\documentclass[a4,11pt]{amsart}
\usepackage{amsmath,amsfonts,amsthm,amssymb,amscd,url}

\newcommand{\quadrec}[2]{({#1}/{#2})}

\DeclareFontFamily{OT1}{rsfs}{}
\DeclareFontShape{OT1}{rsfs}{n}{it}{<-> rsfs10}{}
\DeclareMathAlphabet{\mathscr}{OT1}{rsfs}{n}{it}
\DeclareMathOperator{\sgn}{sgn}

\DeclareMathOperator{\av}{av}

\DeclareMathOperator{\num}{num}
\DeclareMathOperator{\den}{den}

\DeclareMathOperator{\rad}{rad}
\DeclareMathOperator{\mo}{\,mod}
\DeclareMathOperator{\sq}{sq}

\DeclareMathOperator{\Area}{Area}

\DeclareMathOperator{\Res}{Res}
\DeclareMathOperator{\Disc}{Disc}

\DeclareMathOperator{\charac}{char}

\DeclareMathOperator{\rnk}{rank}

\DeclareMathOperator{\Frob}{Frob}
\DeclareMathOperator{\Gal}{Gal}

\newtheorem{prop}{Proposition}[section]
\newtheorem{thm}[prop]{Theorem}
\newtheorem*{main1}{Main Theorem 1}
\newtheorem*{main2}{Main Theorem 2}
\newtheorem{cor}[prop]{Corollary}
\newtheorem{lem}[prop]{Lemma}

\newtheorem*{defn*}{Definition}
\newtheorem*{hypa1}{Hypothesis $\mathscr{A}_1(P)$}
\newtheorem*{hypa2}{Hypothesis $\mathscr{A}_2(P)$}
\newtheorem*{hypb1}{Hypothesis $\mathscr{B}_1(P)$}
\newtheorem*{hypb2}{Hypothesis $\mathscr{B}_2(P)$}
\newenvironment{Rem}{{\bf Remark.}}{}
\numberwithin{equation}{section}
\title{On the behaviour of root numbers in families of elliptic curves}
\author{H. A. Helfgott}
\address{Mathematics Department, University of Bristol, Bristol BS8 1TW,
United Kingdom}
\email{h.andres.helfgott@bristol.ac.uk}
\begin{document}
\begin{abstract}
Let $\mathscr{E}$ be a one-parameter family of elliptic curves over 
$\mathbb{Q}$. 
We prove that the average root number is zero for a large class of families
of elliptic curves of fairly general type. Furthermore, we show that any
family $\mathscr{E}$ with at least one point of multiplicative reduction
over $\mathbb{Q}(T)$ has average root number $0$, provided that
two classical arithmetical conjectures hold
for two polynomials constructed explicitly in terms of $\mathscr{E}$.
The behaviour of the root number in any family $\mathscr{E}$ without
 multiplicative reduction over $\mathbb{Q}(T)$
is shown to be rather regular and non-random; we give expressions for
the average root number in this case. 
\end{abstract}
\maketitle

\setcounter{tocdepth}{1}
\tableofcontents

\section{Introduction}
The precise statement of the main results (listed and discussed in \S \ref{subs:resu}) requires a non-trivial amount of definitions and notation. Let us first
discuss problems and results informally.
\subsection{Summary of results}
Let $E$ be an elliptic curve over $\mathbb{Q}$. Its {\em root number}
$W(E) = \pm 1$ is the sign of its functional equation:
\[\mathcal{N}_E^{(2-s)/2} (2\pi )^{s-2} \Gamma(2-s) L(E,2-s) = 
W(E) \mathcal{N}_E^{s/2} (2\pi)^{-s} \Gamma(s) L(E,s),\]
where $\mathcal{N}_E$ is the conductor of $E$.
Given an elliptic curve\footnote{
We will speak interchangeably about a family of elliptic curves
over a field $K$ and an elliptic curve over $K(T)$. In other
words, we will always understand {\em family} to mean
{\em family fibred over the projective
line}.}
 $\mathscr{E}$ over $\mathbb{Q}(T)$, one may ask
how $W(\mathscr{E}(t))$ varies as $t\in \mathbb{Q}$ varies.

The natural expectation is that $W(\mathscr{E}(t))$ be equidistributed. We
will see that this holds in a strong sense
whenever $\mathscr{E}$ has at least one place
of multiplicative reduction, provided that two standard arithmetical
conjectures hold. Since the conjectures have been proved in some cases,
the result is unconditional for some families. If $\mathscr{E}$ has
no places of multiplicative reduction, the average of $W(\mathscr{E}(t))$
need no longer be zero; we will obtain an expression for the average as
an infinite product.

If the equidistribution of $W(\mathscr{E}(t))$ were proven unconditionally
for any $\mathscr{E}$ for which our results are merely conditional,
a standard and currently
intractable conjecture in number theory
would follow. Thus, we seem to be closing the subject 
of the distribution of the root numbers of elliptic curves
for the time being.

The behaviour of $W(\mathscr{E}(t))$ will be shown to depend strongly
on the arithmetic of $\mathbb{Z}$.
Consider the following examples, 
for which our results are unconditional:\\ \\
\noindent
\begin{tabular}{ll}
$y^2 = x^3 - (t-1)^2 t x/48 -  (t-1)^3 t/864$
&has av. root number $0$ for $t\in \mathbb{Z}$\\
$y^2 = x^3 - x/12 - (11+t)/864$
&has av.\ root number
 $0$ for $t\in \mathbb{Q}$\\
$y^2 = x^3 - x/16 - (2 + 7 t)/864$ &has av.\
root number
 $0$ for $t\in \mathbb{Q}$\\
$y^2 = x^3 - \frac{1}{48} (2 - 4 t + t^2) x -
\frac{1}{864} (3 + 9 t - 6 t^2 + t^3)$ &has av.\ root number $0$ for $t\in \mathbb{Q}$\\
$y^2 = x (x+a) (x+b)$ 
&has av.\ root number $0$ for $a,b\in \mathbb{Z}$.
\end{tabular}\\

(Here ``av.'' stands for ``average''. When taking averages over $\mathbb{Q}$, we order $\mathbb{Q}$ by height;
see \S \ref{subs:preli}.)
We will see that the fact that the first family has average root
number $0$ over every arithmetic progression is equivalent to the
fact that the Liouville function $\lambda$ (or
the M\"obius function $\mu$) averages to zero over every
arithmetic progression. The remaining examples follow from (and
imply) particular cases of the recently solved parity problem for 
degree $3$
(\cite{He3}, \cite{Heirr}; cf.\ \cite{FI}, \cite{HB}, \cite{HBM}).
 The last example
does not require the apparatus developed in the body of this paper,
and will be dealt with briefly (Prop.\ \ref{prop:avrab}). It reduces
to essentially the same problem as the second example.

Consider now the following example:
\[\mathscr{E}:y^2 = x^3 - \frac{1}{48}
f_1(t) f_2(t) (f_1(t)^3 - f_2(t)^3)^2 x
- \frac{1}{864} (f_1(t)^3 + f_2(t)^3) (f_1(t)^3 - f_2(t)^3)^3,\]
where $f_1(t) = - 5 -2 t^2$, $f_2 = 2 + 5 t^2$. We will
prove (unconditionally) that $\mathscr{E}(t)$
has average root number $-0.15294\dotsc$ as $t$ varies over 
$\mathbb{Q}$ (\ref{eq:umi}). 
Notice that the family
has no points of multiplicative reduction as a curve over
$\mathbb{Q}(T)$; only when this is the case can the
average be non-zero over $\mathbb{Q}$. (The exclusion of all other
cases is conditional on the conjectures already mentioned.)

Given any family $\mathscr{E}$ without points of multiplicative reduction
over $\mathbb{Q}(T)$, we will express the average root number of
$\mathscr{E}(t)$ for $t$ varying over $\mathbb{Q}$ or $\mathbb{Z}$
as an infinite product (\S \ref{subs:vogel}, 
(\ref{eq:brah1}) and (\ref{eq:brah2})). 

\begin{center}
* * *
\end{center}

{\em Applications.} Joint work with B. Conrad and K. Conrad (\cite{CCH}), stemming from
the present paper, 
has shown that, over function fields
$\mathbb{F}_q(u)$, there are families $\mathscr{E}$ with 
places of multiplicative reduction and constant root number.
At issue is a deep difference between the arithmetic of $\mathbb{Z}$
and the arithmetic of $\mathbb{F}_q(u)$ (\cite{CCG}).

If the Birch-Swinnerton-Dyer conjecture holds, 
then $\rnk(\mathscr{E}(t)(K))
\geq \rnk(\mathscr{E}(K(T))) + 
\frac{1}{2} (1 - W(\mathscr{E}(t)) (-1)^{\rnk(\mathscr{E}(K(T)))})$ 
for all but finitely many
$t\in K$ (by the specialisation theorem; vd.\ \cite{Sisp}, Thm.\ C).
Thus, the average root number gives a lower bound for
the rank; this fact is used in \cite{CCH} to construct families
with excess rank. The distribution of the root number
is also crucial to issues of the distribution of the zeroes
of $L$-functions of elliptic curves; see \cite{Mi} for
an application of the main result of the present paper.

\subsection{The analytical crux of the problem.}\label{subs:crux}
 It is well-known that there is a
general expression
\begin{equation}\label{eq:tororo} W(E) = - \prod_p W_p(E),\end{equation}
where $W_p$ is a function defined on elliptic curves over $\mathbb{Q}_p$
(the {\em local root number}).
However, we {\em cannot} conclude that
\begin{equation}\label{eq:strauss}
\av_{\mathbb{Z}} W(\mathscr{E}(n)) = - \prod_p \int_{\mathbb{Z}_p} 
W_p(\mathscr{E}(x)) dx
\;\;\;\;\;\;\;\;\;\;\;\;\;\;\text{\bf (FALLACY!)}
\end{equation}
where $\av_{\mathbb{Z}} f(n)$ is the average $\lim_{N\to \infty} 
\frac{1}{N} \sum_{n=1}^N f(n)$. One simply 
cannot do this with infinite products.
 In fact, there are many cases when
(\ref{eq:strauss}) is false, though a different
expression for the average as a product of $p$-adic integrals
is possible (Corollaries \ref{cor:albri} and \ref{cor:albor}
and Propositions \ref{prop:ant1} and \ref{prop:ant2}).

(Even if we did not have an infinite product, we would have to be
a little careful. Consider the function $h_p:\mathbb{Q}_p \to \mathbb{C}$ 
given by
\[h_p(x) = \begin{cases} 1 &\text{if $x\notin \mathbb{Q}$}\\
0 &\text{if $x\in \mathbb{Q}$.}\end{cases}\]
The average of $h_p(x)$ over $\mathbb{Z}$ or $\mathbb{Q}$ has to be zero,
since $h_p(x)$ is identically $0$ on $\mathbb{Q}$. At the same time,
$\int_{\mathbb{Z}_p} h_p(x) dx = 1$, and $0\ne 1$.)


Our main task will be to show that the infinite product
$W(\mathscr{E}(t)) = - \prod_p W_p(\mathscr{E}(t))$ can be rewritten
as a something close to finite product of functions that are locally constant
almost everywhere, times a function of the form $\lambda(P(t))$,
where $\lambda$ is the well-known Liouville function and $P$ is a polynomial. 
(We will actually  be able to do this for all families 
$\mathscr{E}$ over $K(T)$, where $K$ is any global field of characteristic 
$\ne 2,3$. Only Section \ref{sec:avcor} will be specifically about $K=\mathbb{Q}$.) 
Once we have rewritten the root number in this fashion, we will be able
to compute averages correctly.

The fact 
that $\av_{\mathbb{Z}} W(\mathscr{E}(n)) = 0$ in the general case
 will be one of our two main
conditional results; it will be conditional on, and indeed shown to be
equivalent to, a difficult conjecture that remains open in general, viz., that
$\lambda(P(n))$, $P$ a polynomial, averages to zero. 
The other main conditional result will be $\av_{\mathbb{Q}} W(\mathscr{E}(t))
= 0$ -- again, in the general case; there are special families for which
the average is non-zero, and we will classify them and determine their
average root numbers.

Our conditional results ($\av_{\mathbb{Z}} W(\mathscr{E}(n)) = 0$,
$\av_{\mathbb{Q}} W(\mathscr{E}(t))=0$) are unconditional when the polynomial
$P$ we must use is one of those for which the conjecture on $\lambda(P(n))$
has been proven. Some of these cases of the conjecture are due to the
author (\cite{Heirr}, \cite{He3}).

\subsection{Acknowledgements}
The author is most grateful to H. Iwaniec
for bringing the main subject of this paper to his attention
in the correct belief that an analytic problem lurked under the 
surface.
 The advice of B. Conrad and K. Conrad was in general
quite valuable. Thanks are also due to C. Hall, for many 
conversations on elliptic curves, to G. Harcos, for his advice on
formalism, to S. Kobayashi, S. J. Miller, O. Rizzo
and D. Rohrlich, for discussions on their work, and to X. Xarles, who
pointed out errors in the appendices. The opinions of an
anonymous referee also proved very useful. Some computations in the
appendices
were done by SAGE v3.2.1, a free open-source mathematics software system.

\section{Overview}\label{sec:ostrorta}
\subsection{Preliminaries}\label{subs:preli}
\subsubsection{Fields, places and polynomials}\label{subs:fpp}
By a {\em global field} we shall mean either (a) a number field, or 
(b) a function
field over a finite field.

Let $K$ be a global field. 
Given
a place $v$ of $K(T)$, we define the homogeneous polynomial 
$P_v = y^{\deg Q} Q_v(x/y)$ if $v$ is given by a polynomial 
$Q_v\in K \lbrack t\rbrack$, and $P_v = y$ if $v = \deg(\den) -
\deg(\num)$. The choice among all $Q_v$ for a given $v$ is arbitrary; thus,
given a finite, non-empty set $V$ of places of $K$,
we may assume $P_v\in \mathscr{O}_{K,V}\lbrack x,y\rbrack$.

We say that $f\in \mathscr{O}_{K,V}\lbrack x\rbrack$
(resp.\ $f\in \mathscr{O}_{K,V}\lbrack x,y\rbrack$) is {\em square-free}
if there are no polynomials $f_1, f_2\in K\lbrack x\rbrack$, $f_1\notin K$
(resp.\ $f_1, f_2\in K\lbrack x,y\rbrack$, $f_1\notin K$)
such that $f = f_1^2 \cdot f_2$. 

Let $\mathscr{E}$ be an elliptic curve over $K(T)$. Define the
homogeneous polynomials
\begin{equation}\label{eq:M}
M_{\mathscr{E}} = \prod_{\text{$\mathscr{E}$ has mult.\ red.\ at $v$}}
 P_v,\;\;\;\;\;\;\;\;\;\;\;\;\;\;\;\;
B_{\mathscr{E}} = \prod_{\text{$\mathscr{E}$ has q.\  bad red.\ at $v$}}
P_v,
\end{equation}
where we say that $\mathscr{E}$ has quite bad reduction at $v$ if every
quadratic twist of $\mathscr{E}$ has bad reduction at $v$. 
Both $M_{\mathscr{E}}$ and $B_{\mathscr{E}}$ are square-free.
See \S \ref{sss:piu} for criteria for determining the reduction type.

\subsubsection{Averages}\label{subs:pasolini}
We will average functions over the integers, over arithmetic progressions,
over subsets of $\mathbb{Z}^2$,
and over the rationals; we will order $\mathbb{Z}^2$ and the rationals by
height. More precisely, we proceed as follows.
Given a function $f:\mathbb{Z}\to \mathbb{C}$ and an arithmetic progression 
$a + m \mathbb{Z}$, we define
\begin{equation}\label{eq:navy}
\av_{a + m \mathbb{Z}} f = \lim_{N\to \infty}
\frac{1}{N/m} \mathop{\sum_{1\leq n\leq N}}_{n\equiv a \mo m}
 f(n) .\end{equation}
If $\av_{a + m \mathbb{Z}} f = 0$ for all 
$a\in \mathbb{Z}$, $m\in \mathbb{Z}^+$,
we say that $f$ {\em has strong zero average over the integers.}
Given a function $f:\mathbb{Z}^2\to \mathbb{C}$, a lattice coset 
$a + L\subset \mathbb{Z}^2$ and a sector $S\subset \mathbb{R}^2$ (see \S
\ref{subs:seclat}), we define
\begin{equation}\label{eq:pabellon}
\av_{S\cap (a + L)} f = \lim_{N\to \infty}
 \frac{1}{|S\cap (a + L)\cap \lbrack -N, N\rbrack^2|}
 \sum_{(x,y)\in S\cap (a + L)\cap \lbrack -N, N\rbrack^2} f(x,y) .\end{equation}
Here and henceforth, $|A|$ denotes the number of elements of a set $A$.
We say that $f$ {\em has strong zero average over $\mathbb{Z}^2$} if 
$\av_{S\cap (a+L)} f = 0$ for all choices of $S$ and $a+L$.

Given a function $f:\mathbb{Q}\to \mathbb{Z}$, a lattice
coset $a+L\subset \mathbb{Z}^2$ and a sector $S\subset \mathbb{R}^2$, we define
\begin{equation}\label{eq:qavy}\av_{\mathbb{Q},S\cap (a+L)} f =
\lim_{N\to \infty}
 \frac{
 \sum_{(x,y)\in S\cap (a+L)\cap \lbrack -N,N\rbrack^2,\,
\gcd(x,y)=1, y \ne 0} f(x/y) 
}{|\{(x,y)\in S\cap (a+L)\cap \lbrack - N,N\rbrack^2 :
\gcd(x,y)=1\}|} .\end{equation}
 We say that $f$ {\em has strong zero average over the rationals} if
$\av_{\mathbb{Q},S\cap (a+L)} f = 0$ for all choices of $S$ and $a+L$.
We are making our definition of ``strong zero average'' strict 
enough for it to be invariant under fractional linear transformations. 
Moreover, by
letting $S$ be arbitrary, we allow sampling to be restricted to any
open interval in $\mathbb{Q}$. Thus our results will not be imputable
to peculiarities in averaging order or to superficial cancellation.

In all of the above averages, we allow the function $f$ to be undefined for
finitely many possible values of the variables. (If we change $f(7)$ and
only $f(7)$, the average $\av_{\mathbb{Z}} f(n)$ does not change; hence
we may take the average of a function that is undefined at $7$ and only at 
$7$.)


We say that only a {\em zero proportion} of all integers $n$ satisfy 
a given property $\mathbf{P}$ if $|\{1\leq n\leq N : \text{
$\mathbf{P}(n)$ holds}\}|
= o(N)$. Similarly, we say that only a {\em zero proportion} of all pairs
$(x,y)$ of coprime integers satisfy a property $\mathbf{P}$
if $|\{(x,y)\in \lbrack -N,N\rbrack^2: \text{$x$, $y$ are coprime and  
$\mathbf{P}(x,y)$ holds}\}|
= o(N^2)$.

\subsection{Assumptions}\label{subs:htour}
By {\bf Theorem 0.0} ($\mathfrak{X}(P)$, $\mathfrak{Y}(Q)$)
we mean a theorem conditional on hypotheses $\mathfrak{X}$ and $\mathfrak{Y}$
in so far as they concern the objects $P$ and $Q$, respectively. A result
whose statement does not contain parentheses after the numeration should
be understood to be unconditional. 

\begin{hypa1} Let $P\in \mathbb{Z}\lbrack x\rbrack$ be given.
Only for a zero proportion of all integers $n$ do we have
a prime $p> \sqrt{n}$ such that $p^2 | P(n)$.
\end{hypa1} 
\begin{hypa2} Let a homogeneous $P\in \mathbb{Z}\lbrack x,y\rbrack$ be
given. Only for a zero proportion of all pairs of coprime integers
$(x,y)$ do we have a prime $p>\max(x,y)$ such that  $p^2 | P(x,y)$.
\end{hypa2}

Hypotheses $\mathscr{A}_1(P)$ and $\mathscr{A}_2(P)$ are believed to
hold for all square-free $P$. The $abc$-conjecture would imply as much
(\cite{Gran}, Thm.\ 8). We know $\mathscr{A}_1(P)$ (resp.\ $\mathscr{A}_2(P)$)
unconditionally when
$P$ has no irreducible factor of degree larger than $3$ (resp.\ larger than
$6$).
See \cite{Er}, \cite{Es}, \cite{Gr}; vd.\ \cite{He} for sharper bounds.

\begin{hypb1} Let $P\in \mathbb{Z}\lbrack x\rbrack$ be given.
Then $\lambda(P(n))$ has strong zero average over the integers.
\end{hypb1}
\begin{hypb2} Let a homogeneous $P\in \mathbb{Z}\lbrack x,y\rbrack$ be given.
Then $\lambda(P(x,y))$ has strong zero average over $\mathbb{Z}^2$.
\end{hypb2}

(Recall the {\em Liouville function}
$\lambda(n) = \prod_{p|n} (-1)^{v_p(n)}$.)
Hypotheses $\mathscr{B}_1(P)$ and $\mathscr{B}_2(P)$ are believed to
hold for all non-constant, square-free $P$. (This is a conjecture of
Chowla's (\cite{Ch}, p.\ 96), closely related to the Bunyakovsky/Schinzel
conjecture on primes represented by polynomials.)
The prime
number theorem implies
$\mathscr{B}_1(P)$ and $\mathscr{B}_2(P)$ for $\deg P = 1$; essentially
the same analytical techniques
 suffice to prove $\mathscr{B}_2(P)$ 
for $\deg P = 2$.
Hypothesis $\mathscr{B}_2(P)$ has been proved (\cite{Heirr})
for $\deg P = 3$
 by a parity-breaking approach.
 It has also been proved recently
for $\deg P = 4$ and $P$ the product of $4$ distinct linear factors \cite{GT};
it is quite possible that a similar approach will eventually prove
$\mathscr{B}_2(P)$ for $\deg P > 4$, $P$ a product of arbitrarily
many distinct linear factors.

\subsection{Main results}\label{subs:resu}

\begin{main1}[$\mathscr{A}_1(B_\mathscr{E}(t,1))$,
$\mathscr{B}_1(M_\mathscr{E}(t,1))$]
Let $\mathscr{E}$ be a family of elliptic curves over $\mathbb{Q}$. 
Assume that 
$M_\mathscr{E}(t,1)$ is not constant. Then $t\mapsto W(\mathscr{E}(t))$ 
has strong zero average over the integers.
\end{main1}

\begin{main2}[$\mathscr{A}_2(B_\mathscr{E})$,
$\mathscr{B}_2(M_\mathscr{E})$]
Let $\mathscr{E}$ be a family of elliptic curves over $\mathbb{Q}$. 
Assume that 
$M_\mathscr{E}$ is not constant\footnote{In other words,
assume that $\mathscr{E}$ has at least one point of multiplicative
reduction over $\mathbb{Q}(T)$.}. 
Then $t\mapsto W(\mathscr{E}(t))$ has
strong zero average over the rationals.
\end{main2}

We will also be able to prove that we can go in the other direction:
If $\mathscr{A}_1(B_{\mathscr{E}}(t,1))$ 
(resp.\ $\mathscr{A}_2(B_{\mathscr{E}})$) is assumed, then 
$\mathscr{B}_1(M_{\mathscr{E}}(t,1))$ 
(resp.\ $\mathscr{B}_2(M_{\mathscr{E}})$) holds if and only if 
$W(\mathscr{E}(t))$ has strong zero average over the integers (resp.\ over
the rationals). 
Since $\mathscr{B}_1(P)$ 
(resp.\ $\mathscr{B}_2(P)$) does not hold for $P$ constant, the
case $M_{\mathscr{E}}(t,1) = 1$ (or, respectively, $M_{\mathscr{E}} = 1$) is not
covered by Main Theorem 1 (or, respectively, Main Theorem 2). In that case,
we will be able 
to express $\av_{a + m\mathbb{Z}} W(\mathscr{E}(t))$ and $\av_{\mathbb{Q},
S\cap L} W(\mathscr{E}(t))$ as infinite products (Propositions
 \ref{prop:ant1} and \ref{prop:ant2}).


\subsection{Generalisations}\label{subs:atrof}
Let $\mathscr{E}$ now be an elliptic curve over $K(T)$, where $K$ is either a 
number field or a function field. 
The root number $W(\mathscr{E}(t))$ can be described in 
the same way whether or not $K = \mathbb{Q}$ (Thm.\ 
\ref{thm:smet}). Consider first the case of $K$ a number field.
The analogues of $\mathscr{A}_j(P)$ 
(resp.\ $\mathscr{B}_j(P)$) could then probably be proved for $\deg P \leq 6$
(resp.\ $\deg P\leq 3$), just as for $K = \mathbb{Q}$; they are certainly
believed to hold in general.
 It is not clear what the most
natural convention for averaging $W(\mathscr{E}(t))$ is, since
 there will usually be infinitely many elements of $K$ of height below
a given constant; Thm.\ \ref{thm:smet}
ought to yield average zero under any reasonable convention, but a fair amount
of ad-hoc work would be required for any given one. 
Note that, for a generic
number field $K$, the functional equation is still conjectural; the root
number $W(E)$ has to be seen as defined by the product of local root numbers
that is used to compute it.

Let $K$ be a function field. Then  we face a qualitatively different situation:
on the one hand, $\mathscr{A}_1(P)$ and $\mathscr{A}_2(P)$
are known to hold for all square-free $P$ (\cite{Ra}, \cite{Po}); on
the other hand, $\mathscr{B}_1(P)$ and $\mathscr{B}_2(P)$ are 
false for some non-constant, square-free $P$ (\cite{CCG}). There are
 families $\mathscr{E}$ over $K(T)$ with
$W(\mathscr{E}(t))$ constant and $M_{\mathscr{E}}\ne 1$ (\cite{CCH}).

We will express the average of $W(\mathscr{E}(t))$ as an infinite
product when $\mathscr{E}$ is a family without multiplicative
reduction over $K(T)$, where $K = \mathbb{F}_q(u)$ and 
$q$ is not a power of $2$ or $3$. The key intermediate result (Thm.\
\ref{thm:smet}) is valid over any global
field of char.\ $\ne 2,3$. Families of elliptic curves over 
function fields of char.\ $2$ or $3$ present technical difficulties
due to the ubiquity of wild ramification.

\subsection{Relation to the previous literature}
\subsubsection{Algebraic aspect}
The decomposition of the root number $W(E)$ into local root numbers is
a classical result \cite{De}. The local root numbers at places
with residue field characteristic $\ne 2,3$ were made explicit by
Rohrlich (\cite{Ro}, \cite{Ro2}, \cite{Rog}). The case of residue
 field
characteristic $3$ was made explicit by Kobayashi (\cite{Ko}); 
the resulting expressions are complicated enough
 that a similar treatment of residue field characteristic
$2$ is likely
to be nearly unworkable.
({\em Note added during editing:} There
 is some very recent work (over number fields)
due to Dokchitser and Dokchitser \cite{DD}
on the case of residue field characteristic 2. There is also some unpublished
work by D. Whitehouse.)
 Over $\mathbb{Q}$, the root numbers $W_2(E)$,
$W_3(E)$ were described in lengthy tables by 
Halberstadt (\cite{Ha}); no clear pattern
is discernible. As a consequence, much
 work on root numbers up to now (\cite{Riz2}, \cite{Wo}) has had to
involve either very substantial case-work or limitations to certain
convenient congruence classes.

In \S \ref{sec:locroot}, we will show that every local root
number $W(\mathscr{E}(t))$ is locally constant almost everywhere -- not
a surprise -- and, moreover, that it has a very regular kind of
behaviour at the points $t$ where it is not locally constant.
This qualitative characterisation allows us to work 
with the local root numbers at all 
places of a number field 
(or of a function field of $\charac \ne 2,3$) with essentially no case-work.
\subsubsection{Analytic aspect}
Starting with \cite{Ro}, there was a series of papers (\cite{Ma},
\cite{Rog}, \cite{Riz1}, \cite{Riz2}, \cite{Wo}) on average root numbers in
one-parameter families. The cases treated in \cite{Ma}, \cite{Riz1}
and \cite{Riz2} all lack multiplicative reduction, and thus sometimes
have a non-zero average root number, which may be computed with the 
methods developed in the present paper (Propositions \ref{prop:ant1} and
\ref{prop:ant2}). Note that the family treated
in \cite{Riz2} has non-zero average root number over $\mathbb{Z}$,
but not over $\mathbb{Q}$; we will see an example of a non-constant
family with non-zero average root number over $\mathbb{Q}$ in
\S \ref{subs:vogel}. It was intimated in \cite{Ma} that the case
of multiplicative reduction (identified by a certain technical condition,
rather than reduction type) presented essential analytical difficulties.
The likelihood of a relation of some sort
between root number problems and the parity problem was already intuited by
S. Wong \cite{Wo}. One of the main points of this paper is to 
establish a precise relation. 

Our main intermediate result 
(Thm.\ \ref{thm:smet}; see also Prop.\ \ref{prop:ahem}) is valid both over
number fields and function fields. Since, as it was pointed out
by \cite{CCG}, the behaviour of the Liouville function
$\lambda$ on function fields differs
clearly from its behaviour on number fields, Thm.\ \ref{thm:smet}
would lead one to believe that there are families $\mathscr{E}$ over 
$\mathbb{F}_q(u,t)$ for which $W(\mathscr{E}(t))$ cannot
have strong zero average over $\mathbb{F}_q(u)$. This was in fact shown
to be the case in \cite{CCH}, which stemmed in part from an early version of
the present paper.

\section{Notation}
\subsection{Fields and valuations}
For us, a field $K$ will be {\em global} if it is either 
a number field or the function field of a curve over a finite field.
We write $M_K$ for the set of all places of $K$, and $M_{K,\infty}$ for the set
of archimedean places of $K$.
 Given a finite, non-empty set $V$ of places of $K$ containing $M_{K,\infty}$,
we write
 $\mathscr{O}_{K,V}$ for the ring of $V$-integers
of $K$
and $I_{K,V}$ for the semigroup of non-zero ideals of $\mathscr{O}_{K,V}$.
(The ring of $V$-integers $\mathscr{O}_{K,V}$ is defined to be
$\mathscr{O}_{K,V} = \{x\in K : |x|_v\leq 1\; \forall v\in M_K\setminus V\}$.)
By a {\em prime} we shall mean
either a place not in $V$ or the prime ideal of $\mathscr{O}_{K,V}$
corresponding thereto, where $V$ is given.

We denote the algebraic closure
of a field $K$ by $\overline{K}$, and its separable closure by $K_s$. 
Given a global field $K$, a finite Galois extension $L/K$, and a prime
ideal $\mathfrak{p}$ of $K$, 
we denote by $\Frob_{\mathfrak{p}}$ the conjugacy class of Frobenius elements
of $\mathfrak{p}$ in $\Gal(L/K)$.

A field is {\em local} if (a) it is complete with respect to a 
non-trivial discrete 
valuation, and (b) its residue field is finite. Given a local field
$K$ and its valuation $v$, a {\em ball} is a set of the form
$\{x\in K: v(x-x_0)\geq k\}$, where $x_0\in K$ and $-\infty\leq k<\infty$.
(We may also call such a ball an {\em open ball}, as it is open (and closed).)
The valuation $v$ is always normalised so that its range is $\mathbb{Z}$.
The ring of integers $\mathscr{O}_K$ of $K$ is defined to be
$\mathscr{O}_K = \{x \in K : v(x)\geq 0\} = \{x\in K: |x|_v\leq 1\}$.

By a {\em punctured ball} around $x_0$ we mean a set of the form
$U\setminus \{x_0\}$, where $U$ is any ball around $x_0$. In general,
by a {\em punctured neighbourhood} around $x_0$, we mean a neighbourhood
of $x_0$ with $x_0$ itself omitted.

\subsection{Ideals}
Let $R$ be a Dedekind domain. Let $\mathfrak{a}$ be a non-zero ideal of $R$.
We let $\sq(\mathfrak{a}) = \prod_{\mathfrak{p}^2 | \mathfrak{a}}
 \mathfrak{p}^{v_{\mathfrak{p}}(\mathfrak{a}) - 1}$,
$\lambda(\mathfrak{a}) = \prod_{\mathfrak{p} | \mathfrak{a}}
 (-1)^{v_{\mathfrak{p}}(\mathfrak{a})}$. By $\mathfrak{a} | \mathfrak{b}^{\infty}$ we will mean that $\mathfrak{p} | \mathfrak{b}$ for every prime
ideal $\mathfrak{p}$ dividing $\mathfrak{a}$.

\subsection{Polynomials}\label{subs:husb}
Let $R$ be a Dedekind domain; let $K$ be its field of fractions.
As is standard, we 
define the resultant of two polynomials $f, g\in K\lbrack x\rbrack$ 
to be the determinant of the corresponding Sylvester matrix.
Given two homogeneous polynomials $f, g\in K\lbrack x, y\rbrack$,
we define the resultant to be the determinant of the Sylvester matrix
\[\left(\begin{array}{ccccccc}
a_n & \dotsb & \dotsb & a_0 & 0 &\dotsb &0\\
\vdots & \vdots & \vdots & \vdots & \vdots & \vdots & \vdots \\
0 &\dotsb &a_n & a_{n-1} & \dotsb &\dotsb & a_0\\
b_m &\dotsb &b_0 & 0 & 0 &\dotsb & 0\\
\vdots & \vdots & \vdots & \vdots & \vdots & \vdots & \vdots \\
0 & \dotsb & 0 & b_{m} &b_{m-1} &\dotsb & b_0 \end{array}\right)\]
built up from the coefficients $a_j$, $b_j$ of
$f(x,y)=\sum_{j=0}^n a_j x^j y^{n-j}$, 
$g(x,y) = \sum_{j=0}^m b_j x^j y^{n-j}$, where $n= \deg f$, $m=\deg g$.

If $f, g\in R\lbrack x\rbrack$ 
are non-zero and coprime as elements of $K\lbrack x\rbrack$, then $\Res(f,g)$
is an element of $R$ divisible by the ideal $\gcd(f(a),g(a))$
 for every $a\in R$. If $f, g\in R\lbrack x,y\rbrack$ are non-zero,
homogeneous and coprime as elements of $R\lbrack x,y\rbrack$,
 then $\Res(f,g)$ is
an element of $R$ divisible by the ideal $\gcd(f(a,b),g(a,b))$ for
every pair $(a,b)\in R^2$ with $\gcd(a,b)=1$. More generally, 
if $(a,b)\in R^2$ satisfies $\gcd(a,b)|\mathfrak{d}$ for some ideal 
$\mathfrak{d}$ of $R$, then $\gcd(f(a,b),g(a,b))$ divides 
$\mathfrak{d}\cdot \Res(f,g)$. 

\subsection{Elliptic curves}
Let $E$ be an elliptic curve over a field $K$. We write
$E(\lbrack m\rbrack)$ for the set of points of order $m$ on $E$, and
$T_{\ell}(E)$ for the Tate module.
If $K$ is a local field, we denote the 
reduction of $E$ by $\widehat{E}$.
When we say that $E$ has potentially good reduction over $K$, we allow
the possibility that $E$ may actually have good reduction over $K$, as opposed
to additive, potentially good reduction; the same holds with
``multiplicative'' instead of ``good''.
\subsection{Characters}
A {\em character} of an abelian topological
group $G$ is a continuous homomorphism from $G$
to the unit circle
$S^1\subset \mathbb{C}$. Given a global field $K$ and a prime ideal 
$\mathfrak{p}$,
we denote by $(\cdot/\mathfrak{p})$ the quadratic reciprocity
symbol, that is, the character $(\cdot/\mathfrak{p})$ on
$K_{\mathfrak{p}}^*$
given by
\[(x/\mathfrak{p}) = \begin{cases} 1 & \text{if $x = y^2$ for some 
$y\in K_{\mathfrak{p}}^*$,}\\
-1 & \text{otherwise.} \end{cases}\]

\subsection{Sectors and lattices}\label{subs:seclat}
A \emph{lattice} is a subgroup of $\mathbb{Z}^n$ of finite index; 
a \emph{lattice coset} 
 is a coset of such a subgroup. 
The {\em index} 
$\lbrack \mathbb{Z}^n : L\rbrack$
of a lattice $L$ 
is simply its index as a subgroup of $\mathbb{Z}^n$, i.e., the
number of elements of the quotient $\mathbb{Z}^n/L$.
By a \emph{sector} of $\mathbb{R}^n$ we will mean a connected component of a set of the form $%
\mathbb{R}^{n}-(T_{1}\cap T_{2}\cap \dotsb \cap T_{n})$, where $T_{i}$ is a
hyperplane going through the origin.
(We will always have $n=2$.)

\subsection{Complex-valued functions}\label{subs:funcs}
Let $K$ be a field. By a {\em complex-valued function} on $\mathbb{P}^1(K)$ 
we will
mean a function $f:U\to \mathbb{C}$ defined on the complement
$U = \mathbb{P}^1(K) - 
\{x_1,\dotsc,x_n\}$ of a finite set $\{x_1,\dotsc,x_n\}$; by a {\em complex-valued
function} on $\mathbb{A}^2(K)$ we will mean a function $f:U\to \mathbb{C}$ defined
on the complement $U\subset \mathbb{A}^2(K)$ of a finite set of
lines through the origin in $\mathbb{A}^2(K)$. (As is usual,
we write $\mathbb{P}^1(K)$ for the projective
line $K^2/K^*$, and $\mathbb{A}^2(K)$ for the affine plane $K^2$.)
 Two complex-valued functions on
$\mathbb{P}^1(K)$ (or on $\mathbb{A}^2(K)$) will be regarded
as identical if they agree on the intersection of their domains of definition.
We write $C_{\mathbb{P}^1(K)}$ for the set of all complex-valued
functions on $\mathbb{P}^1(K)$, and $C_{\mathbb{A}^2(K)}$ for the
set of all complex-valued functions on $\mathbb{A}^2(K)$. (Formally,
we may define $C_{\mathbb{P}^1(K)}$ as the set of all equivalence
classes of pairs $(U,f)$, where $(U,f) \sim (V,g)$ if $f|_{U\cap V} =
g|_{U\cap V}$; the same may be done for $C_{\mathbb{A}^2}(K)$.)

\section{Local root numbers}\label{sec:locroot}
 Let $K$ be a local field of characteristic $\ne 2,3$. Let $\mathscr{E}$
be an elliptic curve over $K(T)$. We will investigate the local behaviour
of $t\to W(\mathscr{E}(t))$; in particular, we will show that
$t\to W(\mathscr{E}(t))$ is locally constant almost everywhere. We will
also see how $t\to W(\mathscr{E}(t))$ behaves at the points where it
is not locally constant.

While we allow ourselves to assume that $\charac(K)\ne 2,3$, we must include
in our treatment local fields $K$ whose residue fields have 
characteristic $2$ or $3$, as such fields $K$
 arise as localisations of number fields. 
Thus we cannot simply make use of formulae valid only 
when the characteristic of the residue field
is greater than $3$ (Prop.\ \ref{prop:rohr}). 
 
\begin{lem}\label{lem:woggle}
 Let $K$ be a local field. 
Let $m\geq 2$ be an integer prime to the characteristic of
the residue field of $K$. Let a minimal Weierstrass model over $K$
give an elliptic curve $E$ with good reduction.
Then the map
\[E\lbrack m\rbrack \to \widehat{E}\lbrack m\rbrack\]
induced by the reduction of the model is a bijective homomorphism.
\end{lem}
\begin{proof}
The map is injective by \cite{Si}, Ch.\ VII, Prop.\ 3.1(b)
and a homomorphism by \cite{Si}, Ch.\ VII, Prop.\ 2.1.
Since $m$ is prime to the residue characteristic of $K$, the finite
sets $E\lbrack m\rbrack$ and $\widehat{E}\lbrack m\rbrack$ have the
same cardinality, and thus an injective map from one to the other must be
a bijection.
\end{proof}

It has been pointed out by B. Conrad \cite{Co} that the following proposition
can be derived from a general result \cite{Ki} on the local constancy of
local systems over schemes of finite type over local fields. A
representation-theoretical argument is also possible (\cite{Roc}).
We give an elementary proof.
\begin{prop}\label{prop:waggle}
Let $K$ be a local field. Let $a_1,\dots,a_6\in \mathscr{O}_K$ 
be the coefficients of a Weierstrass model for an elliptic curve $E$
with potential good reduction. Then there is an $\epsilon>0$ such that
any $\mathbf{a}_1,\dots,\mathbf{a}_6\in \mathscr{O}_K$ with 
$|\mathbf{a}_j - a_j|<\epsilon$
are the coefficients of a Weierstrass model for an elliptic curve $\mathbf{E}$
 with $W(\mathbf{E}) = W(E)$. 
\end{prop}
\begin{proof}
Let $E$ acquire good reduction over a finite extension $L/K$. 
Then there is a change of variables $x' = u^2 x + r$, 
$y' = u^3 y + u^2 s x + t$ with $u, r, s, t\in \mathscr{O}_L$ sending
the original model $a_1,\dotsc,a_6$ to a minimal model 
$a_1',\dotsc,a_6'$ with good reduction. By Lem.\ \ref{lem:woggle},
any two distinct points $P_1, P_2 \in E\lbrack m\rbrack$ have either
$|x(P_1) - x(P_2)|\geq 1$ or $|y(P_1) - y(P_2)|\geq 1$ in the minimal model,
with $|x(P_1) - x(P_2)|<1$ occurring only if $x(P_1) = x(P_2)$.
Thus, either $|x(P_1) - x(P_2)|\geq |u|^{-2}$ or
$|y(P_1) - y(P_2)|\geq |u|^{-3}$ in the original model.

Let $\mathbf{a}_1,\dotsc,\mathbf{a}_6\in \mathscr{O}_K$ satisfy $|\mathbf{a}_j - a_j| < |u|^7$.
Then the same change of variables as above will send 
$\mathbf{a}_1,\dotsc,\mathbf{a}_6$
to a minimal model $\mathbf{a}_1',\dotsc,\mathbf{a}_6'\in \mathscr{O}_K$,
$|\mathbf{a}_j - a_j|<1$, describing an elliptic curve $\mathbf{E}$ 
with good reduction. By Lem.\ \ref{lem:woggle},
the maps $E\lbrack \ell^n\rbrack \to \widehat{E}\lbrack \ell^n \rbrack$ and
$\mathbf{E}\lbrack \ell^n\rbrack 
\to \widehat{\mathbf{E}}\lbrack \ell^n\rbrack$ 
induced by the reduction of the minimal models are bijective
homomorphisms. There is an isomorphism
$\widehat{E}\lbrack \ell^n\rbrack \sim \widehat{\mathbf{E}}\lbrack \ell^n\rbrack$,
since the minimal models reduce to the same model for
$\widehat{E}$ and $\widehat{\mathbf{E}}$. Thus we have
a bijective homomorphism 
$\iota_n:E\lbrack \ell^n\rbrack \to \mathbf{E}\lbrack \ell^n\rbrack$
with $|x(P) - x(\iota_n(P))|<1$, $|y(P) - y(\iota_n(P))|<1$ in the minimal
models, and 
$|x(P) - x(\iota_n(P))|<|u|^{-2}$, $|y(P) - y(\iota_n(P))|<|u|^{-3}$ 
in the original models.

Now let $\gamma
\in \Gal(\overline{K}/K)$. Then, in the original models,
 $|x(\gamma(P)) - x(\gamma(\iota_n(P)))|
= |\gamma(x(P)) - \gamma(x(\iota_n(P)))|
< |u|^{-2}$, $|y(\gamma(P)) - y(\gamma(\iota_n(P)))|<|u|^{-3}$ 
 (since the norm is Galois-invariant) and
$|x(\gamma(P)) - x(\iota_n(\gamma(P)))|<|u|^{-2}$, 
$|y(\gamma(P)) - y(\iota_n(\gamma(P)))|<|u|^{-3}$.
As we would have $|x(\gamma(\iota_n(P))) - x(\iota_n(\gamma(P)))|
\geq |u|^{-2}$ or $|y(\gamma(\iota_n(P))) - y(\iota_n(\gamma(P)))|
\geq |u|^{-3}$ if $\gamma(\iota_n(P))$ were not equal to $\iota_n(\gamma(P))$,
we can conclude that $\gamma(\iota_n(P)) = \iota_n(\gamma(P))$, i.e.,
$\iota_n$ respects the action of the Galois group on $E\lbrack \ell^n\rbrack$
and $\mathbf{E}\lbrack \ell^n\rbrack$.
Therefore, $\iota_n$ induces an isomorphism of Tate modules 
$T_{\ell}(E)\to T_{\ell}(\mathbf{E})$ as $\Gal(\overline{K}/K)$-modules. 
Since $E$ (and hence $\mathbf{E}$) has potentially good reduction, the root
number $W(E)$ (resp.\ $W(\mathbf{E})$) 
is determined by the representation of the Weil
group $\mathscr{W}(K)\subset \Gal(\overline{K}/K)$ on $T_{\ell}(E)$
(resp.\ $T_{\ell}(\mathbf{E})$). Hence $W(E) = W(\mathbf{E})$.
\end{proof}

\begin{prop}\label{prop:thor}
Let $K$ be a local field of characteristic $\ne 2$. 
Let $a_1,\dotsc,a_6\in K$ induce a tuple
$(c_4,c_6,\Delta)\in K^3$ with $c_4\ne 0$, $\Delta=0$. Then there
is a constant $w\in \{-1,1\}$ and
an $\epsilon>0$ such that any $\mathbf{a}_1,\dotsc,\mathbf{a}_6\in K$ with
$|\mathbf{a}_j-a_j|<\epsilon$ and $\mathbf{\Delta}\ne 0$ 
give a Weierstrass model for an elliptic curve $\mathbf{E}$ with
$W(\mathbf{E}) = w$.
\end{prop}
\begin{proof}
We may assume $a_1,\dotsc,a_6\in \mathscr{O}_K$. We wish to show first that
there is one change of variables that sends any nearby $\mathbf{a}_1,\dotsc,\mathbf{a}_6\in 
\mathscr{O}_K$ with $\mathbf{\Delta}\ne 0$ to a minimal Weierstrass model.
 Whether or not a
 change of variables $x' = u^2 x + r$, $y' = u^3 y + u^2 s x + t$
($u,r,s,t\in \mathscr{O}_K$, $u\ne 0$) 
sends a Weierstrass model over $\mathscr{O}_K$ to another Weierstrass
model over $\mathscr{O}_K$ depends on $r$, $s$ and $t$ only modulo 
$u^6 \mathscr{O}_K$; small changes in $u$ do not matter either.
Since $\{u\in K : 1\leq |u| \leq C\}$ is compact for any 
$C$, it follows that, if we need consider only changes of variables
with $u$ bounded in norm from above, 
we may consider only a finite number of changes of 
variables. This is indeed the case, as we must have $|u|\leq |c_4|^{-1/4}$.
For each of these finitely many changes of variables, there is a neighbourhood
$V$ of $(a_1,\dotsc,a_6)$ such that
 either all $(\mathbf{a}_1,\dotsc,\mathbf{a}_6)\in V$ are sent
to integer tuples $(\mathbf{a}_1',\dotsc,\mathbf{a}_6')$, or none of them are. Let $U$
be the intersection of all such neighbourhoods. Then $U$ is an open set
for which there is a single change of variables sending any
$(\mathbf{a}_1,\dotsc,\mathbf{a}_6)\in \mathscr{O}_K$ with $\mathbf{\Delta}\ne 0$
to a minimal 
Weierstrass model.

We may now let that single change of variables act on $(a_1,\dotsc, a_6)$ 
and its
neighbourhood; thus, we may assume that any $(\mathbf{a}_1,\dotsc,
\mathbf{a}_6)\in U$
is minimal. Let $V\subset U$ be the set of all tuples
$(\mathbf{a}_1,\dotsc,\mathbf{a}_6)\in U$ satisfying
$\mathbf{\Delta}\ne 0$,
$v(\mathbf{a}_j-a_j)\geq 0$ and $|j(\mathbf{a}_1,\dotsc,\mathbf{a}_6)|>1$. 
Then $(\mathbf{a}_1,\dotsc,\mathbf{a}_6)$
will have (a) split multiplicative reduction, (b) 
unsplit multiplicative reduction or (c) additive, potentially multiplicative
 reduction
depending on whether (a) $v(c_4)=0$ and the reduction of
$a_1,\dotsc,a_6$ describes a singular cubic with a split node, (b)
$v(c_4)=0$ and the reduction of $a_1,\dotsc,a_6$ describes a singular
cubic with an unsplit node, (c) $v(c_4)>0$. Thus, the reduction type
is constant in the neighbourhood $\{(\mathbf{a}_1,\dotsc,\mathbf{a}_6) 
\in U: |\mathbf{a}_j - a_j|<1, |j(\mathbf{a}_1,\dotsc,\mathbf{a}_6)|>1\}$ 
of $(a_1,\dotsc,a_6)$. When an elliptic curve has 
multiplicative reduction, the root number
depends only on whether the reduction is split or unsplit;
when an elliptic curve has additive, potentially multiplicative reduction,
the root number depends only on $K$ and the class of $-c_6$ in
$K^*/(K^*)^2$ (vd. \cite{Rog}, Thm.\ 2, (ii); cf. \cite{Con1}, \S 3).
It follows that the root
number is constant for all tuples $(\mathbf{a}_1,\dotsc,
\mathbf{a}_6)$ with $\mathbf{\Delta} \ne 0$ in an open neighbourhood
$V_0\subset V$ of $(a_1,\dotsc,a_6)$.
\end{proof}

Let us now examine each kind of reduction that an elliptic curve
$\mathscr{E}$ over $K(T)$ may have.
\begin{lem}\label{lem:goodloc}
Let $K$ be a local field. Let $\mathscr{E}$ be an elliptic curve over $K(T)$.
Suppose $\mathscr{E}$ has a place of good reduction at $(T-t_0)$,
where $t_0\in K$. Then there is a neighbourhood $U$ of $t_0$ on which
$W(\mathscr{E}(t))$ is constant.
\end{lem}
\begin{proof}
If $\mathscr{E}(t_0)$ has potentially good reduction, apply
 Prop.\ \ref{prop:waggle}. If $\mathscr{E}(t_0)$ has potentially multiplicative
reduction, proceed as in the last paragraph of the proof of Prop.\
\ref{prop:thor}.
\end{proof}
\begin{lem}\label{lem:multloc}
Let $K$ be a local field. Let $\mathscr{E}$ be an elliptic curve over $K(T)$.
Suppose $\mathscr{E}$ has a place of multiplicative reduction at $(T-t_0)$,
where $t_0\in K$. Then there is a punctured neighbourhood $U$ of $t_0$ on which
$W(\mathscr{E}(t))$ is constant.
\end{lem}
\begin{proof}
Let $(c_4,c_6,\Delta)$ be parameters of a Weierstrass equation for
$\mathscr{E}$ minimal over $K(T)$ with respect to the valuation induced
by $(T-t_0)$. Then $(T-t_0)$ divides $\Delta$ but not
$c_4$; moreover, $(T-t_0)$ is absent from the denominators of the parameters
$a_1,\dotsc,a_6$. Apply Proposition \ref{prop:thor}.
\end{proof}

\begin{lem}\label{lem:potmultoc}
Let $K$ be a local field of characteristic $\ne 2,3$. Let $\mathscr{E}$ be an elliptic curve over
$K(T)$. Suppose $\mathscr{E}$ has a place of additive, potentially
multiplicative reduction
at $(T-t_0)$, where $t_0\in K$. Then there is a punctured neighbourhood
$U$ of $t_0$ on which $W(\mathscr{E}(t))$ depends only on 
$(t-t_0) \mo {K^*}^2$.
\end{lem}
\begin{proof}
Set $t_0=0$ for notational simplicity. Let 
$\mathscr{E}$ be given by parameters $c_4, c_6\in K(T)$ minimal with respect to
the valuation $v_T$ induced by $T$.  Let $a\in K^*$. For 
$t\in a {K^*}^2$, we may substitute $t = a (t')^2$ and obtain
new parameters $c_4'(t') = c_4(a(t')^2)$, $c_6'(t') = c_6(a(t')^2)$
describing an elliptic curve $\mathscr{E}'$ over $K(T')$. Since
$v_T(c_4)=4 k + 2$, $v_T(c_6) = 6 k + 3$, we obtain
$v_{T'}(c_4') = 4 (2 k + 1)$, $v_{T'}(c_6') = 6 (2 k + 1)$, and thus
$\mathscr{E}'$ has multiplicative reduction over $K(T')$ 
at $T'$. The statement now
follows by Lem.\ \ref{lem:multloc}.
\end{proof}

\begin{lem}\label{lem:potgoo}
Let $K$ be a local field of characteristic $\ne 2,3$. Let $\mathscr{E}$
be an elliptic curve over $K(T)$. Suppose $\mathscr{E}$ has a place
of potential good reduction at $(T-t_0)$, where $t_0\in K$. Then
there is a punctured neighbourhood $U$ of $t_0$ on which $W(\mathscr{E}(t))$
depends only on $(t-t_0) \mo {K^*}^{12}$.
\end{lem}
\begin{proof}
Set $t_0=0$ for notational simplicity.
Let 
$\mathscr{E}$ be given by parameters $c_4, c_6\in K(T)$ minimal with respect to
$v_T$.
Since $\charac(K)\ne 2,3$, 
the curve $\mathscr{E}$ acquires good reduction over $K(T^{1/12})$.
Let $a\in K^*$. For $t\in a {K^*}^{12}$, we may substitute 
$t = a (t')^{12}$ and obtain new parameters $c_4'(t') = c_4(a (t')^{12})$,
$c_6'(t') = c_6(a(t')^{12})$ describing an elliptic curve $\mathscr{E}'$
over $K(T')$. Since $\mathscr{E}$ has good reduction over $K(T^{1/12})$,
there is a non-negative integer $k$ such that $v_{T'}(c_4') = 4 k$,
$v_{T'}(c_6') = 6 k$, $v_{T'}(\Delta') = 12 k$. 
Thus $\mathscr{E}'$ has good reduction over
$K(T')$ at $T'$. Apply Lem.\ \ref{lem:goodloc}.
\end{proof}

\begin{thm}\label{thm:krone}
Let $K$ be a local field of characteristic $\ne 2,3$. Let
$\mathscr{E}$ be an elliptic curve over $K(T)$. 
Then, for every $t_0\in \mathbb{P}^1(K)$,
there is a punctured neighbourhood $U_{t_0}$ of $t_0$ 
such that, for
$t\in U_{t_0}$, the value of $f(t)$ depends only on 
$(t-t_0) \mo (K^*)^{12}$. (For $t_0=\infty$, read $1/t$ instead of $t-t_0$.)

Then, for every $t_0\in K$,
there is a punctured neighbourhood $U_{t_0}$ of $t_0$ and 
a function
$g_{t_0}:K^*/(K^*)^n\to \mathbb{C}$ such that $f|_{U_{t_0}}(t) = g(t-t_0)$
for all $t\in U_{t_0}$. (For $t_0=\infty$, read $1/t$ instead of $t-t_0$.)
\end{thm}
Saying that the value of $f(t)$ depends only on
$(t-t_0) \mo (K^*)^n$ for $t\in U_{t_0}$ is the same as saying that
$f(t)=g(t)$ for all $t\in U_{t_0}$ and some function $g:K^*/(K^*)^n\to
\mathbb{C}$.
\begin{proof}
Immediate from Lemmas \ref{lem:goodloc}--\ref{lem:potgoo}.
\end{proof}

\begin{center}
* * *
\end{center}

Let us see at what sort of conclusion we have arrived. What does Theorem
\ref{thm:krone} tell us about the local root number $W(\mathscr{E}(t))$?

Let $f$ be any convex-valued function on $\mathbb{P}^1(K)$, where $K$
is a local field. Suppose that, for every $t_0\in \mathbb{P}^1(K)$, there
is a punctured neighbourhood $U_{t_0}$ of $t_0$ such that, for $t\in U_{t_0}$,
the value of $f(t)$ depends only on $(t-t_0) \mo (K^*)^n$. (For $t_0=\infty$, read
$1/t$ instead of $t-t_0$.) 
Here $n$ is
a fixed integer; assume $\charac(K)\nmid n$. 
What can we say about $f$?

First of all, $f$ is locally constant almost everywhere. This can be shown as 
follows. Let the punctured neighbourhoods $U_{t_0}$ be as above. Then
$\{U_{t_0}\cup \{t_0\}\}_{t_0\in \mathbb{P}^1(K)}$ is an open cover of 
$\mathbb{P}^1(K)$. Since $\mathbb{P}^1(K)$ is compact, there is a finite
subcover 
$\{U_{t_0}\cup \{t_0\}\}_{t_0\in S}$, $S\subset \mathbb{P}^1(K)$ finite. The
claim
is that, for every $t_1\notin S$, the function $f(t)$ is constant in a 
neighbourhood of $t_1$. This is easy: the point
$t_1$ must lie in some set $U_{t_0}$, $t_0\in S$,
$t_0\ne t_1$; in that set $U_{t_0}$, $f(t)$ depends only on $(t-t_0)
\mod (K^*)^n$; thus, $f(t)$ is constant on the set
$V = (t_0 + (t_1-t_0) \cdot (K^*)^n) \cap U_{t_0}$, which contains $t_1$.
Since $\charac(K)\nmid n$, the set $(K^*)^n$ is an open subset of $K^*$
(this is a consequence of Hensel's lemma) and so $V$ is an open neighbourhood
of $t_1$ on which $f(t)$ is constant.


More can be said about such functions $f$. We could have proven that
$t\to W(\mathscr{E}(t))$ is locally constant almost everywhere directly
from Lemma \ref{lem:goodloc}.
What we proved in Theorem  \ref{thm:krone} is stronger. Here is a way to
use this stronger kind of conclusion. It will be useful later.

\begin{lem}\label{lem:griep}
Let $K$ be a local field. Let $n$ be a positive integer such that
$\charac(K)\nmid n$.
Let $f$ be a convex-valued function on $\mathbb{P}^1(K)$.
Suppose that, for every $t_0\in K$, there
is a punctured neighbourhood $U_{t_0}$ of $t_0$ such that, for $t\in U_{t_0}$,
the value of $f(t)$ depends only on $(t-t_0) \mo (K^*)^n$. 
Suppose also that $f(t)\in \mathbb{Q}$ for all but finitely many $t$.

Then, for every ball $U = \{x\in K: |x-x_0|<r\}$ with $x_0\in K$ and $r<\infty$,
\begin{equation}\label{eq:gotod}\int_{U} f(t) dt\end{equation}
is a rational number.
\end{lem}
We do not need any assumptions about what happens at infinity because
we are working on a neighbourhood $U$ disjoint from infinity. (We do not
consider an integral over all of $K$ because such an integral would not in
general be finite.)
\begin{proof}
We can assume that each punctured neighbourhood $U_{t_0}$ is a punctured ball
around $t_0$, as we can simply replace each open set $U_{t_0}$ by
a punctured ball around $t_0$ contained in it.
Since $U$ is compact, there is a finite cover of $U$ by sets of
the form $U\cap (U_{t_0}\cup \{t_0\})$. In a local field, given two balls,
either they are disjoint or one contains the other. (This is so because
local fields are non-archimedean.) Hence we actually have a cover of $U$
by a finite number of disjoint balls $U_{t_0}\cup \{t_0\}$ contained in $U$.
It will be enough to show that
\[\int_{U_{t_0}} f(t) dt\]
is rational.

By the statement, the values of $f(t)$ on $U_{t_0}$ are rational (with finitely many
exceptions) and depend only on $(t-t_0) \mo (K^*)^n$. 
The punctured ball $U_{t_0}$ is of the form 
$U_{t_0} = V_0 +t_0$ for some punctured ball $V_0$ around the origin.
We can write
\begin{equation}\label{eq:yotor}\int_{U_{t_0}} f(t) dt = \sum_{g\in K^*/(K^*)^n}
\int_{g (K^*)^n \cap V_0} f(t+t_0) dt =
\sum_{g\in K^*/(K^*)^n} f(g+t_0) \Area(g (K^*)^n \cap V_0),\end{equation}
where we write $f(g)$ 
for the value taken by $f(t+t_0)$ for all $t\in g (K^*)^n$. 
We may assume without loss of generality that
$V_0 = \{x\in K^*: v(x)\geq 0\}$. (This assumption is notationally convenient.)
To get that the integral in (\ref{eq:yotor}) is rational,
it will be enough to show that $K^*/(K^*)^n$ is finite and that
every set of the form $g (K^*)^n\cap \{x\in K^*: v(x)\geq 0\}$ 
has rational area.

The quotient $K^*/(K^*)^n$ is finite by Hensel's lemma. (It is here
that we use the assumption that $\charac(K)\nmid n$.) Every coset of
$(K^*)^n$ can be written in the form $a (K^*)^n$ with $a\in K^*$,
$0\leq v(a) < n$. Now 
\[\begin{aligned}
\Area&(\{x\in K^*: v(x)\geq 0,\; x\in a (K^*)^n\}) \\ &= 
\Area(\{x\in K^*: v(x)= 0,\; x\in a (K^*)^n\}) \cdot (1 + q^{-n} + q^{-2 n} +
\dotsc)\\ &= \frac{1}{1-q^{-n}} \Area(\{x\in K^*: v(x)= 0,\; x\in a (K^*)^n\}),
\end{aligned}\]
where $q$ is the number of elements in the residue field of $K$.
Again by Hensel's lemma, the set $\{x\in K^*: v(x)= 0, x\in a (K^*)^n\}$
is a union of (additive) cosets of $\mathfrak{p}^{\ell}$ for some
fixed $\ell$; thus, its area is a multiple of $q^{-\ell}$, and, in particular,
it is a rational number.
\end{proof}


%

We will need to know integrals such as (\ref{eq:gotod}) for the following 
reason. Say $K$ is the localisation $\mathbb{Q}_p$ of $\mathbb{Q}$ at a prime 
$p$. Then the function $t\to W_p(\mathscr{E}(t))$ from $\mathbb{Q}_p$ to
$\{-1,1\}$ induces a function from $\mathbb{Q}$ to $\{-1,1\}$, simply by
the inclusion $\mathbb{Q}\subset \mathbb{Q}_p$. (We write $W_p$ for the
local root number in the local field $\mathbb{Q}_p$.)
 We are interested in the
average of $W_p(\mathscr{E}(t))$ over $\mathbb{Z}$. It is not hard to show that
\begin{equation}\label{eq:kotok}\av_{\mathbb{Z}} W_p(\mathscr{E}(t)) = \int_{\mathbb{Z}_p} W_p(\mathscr{E}(t))
dt.\end{equation}
(Here we use the fact that $W_p(\mathscr{E}(t))$ is locally constant almost
everywhere.) By Thm.\ \ref{thm:krone}, it follows that 
$\av_{\mathbb{Z}} W_p(\mathscr{E}(t))$ is a rational number.

Notice, moreover, that this is a rational number that can be computed in
finite time given the finite data describing $W_p(\mathscr{E}(t))$. There
should be no concerns about the computability or constructibility of $p$-adic
integrals of this kind. 
One can furthermore show that the averages
$\frac{1}{N} \sum_{1\leq n\leq N} W_p(\mathscr{E}(n))$ not only converge
to $\int_{\mathbb{Z}_p} W_p(\mathscr{E}(t))$, but do so very rapidly. (For
more
on rates of convergence, see Appendix \ref{sec:ratcon}.)

We will not go into more details here, since we will prove more general
results later. In particular, we will see that equalities such as
(\ref{eq:kotok}) hold for all products of finitely many functions such as 
$W_p(\mathscr{E}(t))$
(as opposed to just one individual function $W_p(\mathscr{E}(t))$.

The key thing to realise here is that, as we saw in the introduction 
(\S \ref{subs:crux}), such equalities are not in general true for infinite 
products. Now, the global root number -- our central object of study in this
paper -- is an infinite product 
$W(\mathscr{E}(t)) = \prod_v W_v(\mathscr{E}(t))$. Our task is then
to re-express it as a finite product, or -- when that is not possible --
 as the product of a finite product and an infinite product of a special kind.

It is no wonder that reciprocity will play a role here. Indeed, the law
of quadratic reciprocity -- just like higher reciprocity laws -- can
be expressed as the fact that $\prod_v \left(\frac{a,b}{v}\right)=1$,
where $v$ goes over all places of a global field $K$ and
$\left(\frac{a,b}{v}\right)$ is the Hilbert symbol. (More about this later.)
In other words, a product of Hilbert symbols over all but finitely many
places of $K$ is equal to (the multiplicative inverse of) a finite product
of Hilbert symbols. This is what lies at the core of the following section.
\section{Reciprocity and polynomials}\label{eq:ostorma}
\subsection{Introduction}
Characters to fixed moduli are harmless; characters to variable moduli are
not. Consider, for instance, the Jacobi symbol $\left(\frac{\cdot}{\cdot}\right)$, which
is a character modulo $b$. Suppose we know how to estimate sums of the form
\begin{equation}\label{eq:inplus}
\mathop{\sum_{x,y\leq N}}_{x\equiv x_0,\; y\equiv y_0 \mo m} f(x,y)\end{equation}
for all $m\in \mathbb{Z}^+$, $x_0,y_0\in \mathbb{Z}/m\mathbb{Z}$, where
$f:\mathbb{Z}^2\to\mathbb{C}$ is some function. Then we know how to estimate
\begin{equation}\label{eq:incirc}
\sum_{x,y\leq N} \left(\frac{x}{5}\right) f(x,y)
\end{equation}
but not how to estimate
\begin{equation}\label{eq:inast}
\sum_{x,y\leq N} \left(\frac{x}{y}\right) f(x,y) .
\end{equation}
We may be able to estimate a sum such as (\ref{eq:inast}) in some particular
cases (see, for instance, Appendix \ref{sec:goron}); however,
in general, a factor of $\left(\frac{x}{y}\right)$ is a very great bother.

When we start to compute the averages of the root number $W(\mathscr{E}(t))$,
we shall have to deal with sums such as
\begin{equation}\label{eq:brascho}
\sum_{x,y\leq N} \left(\frac{2 x + 7 y}{x^2 + 3 x y + 5 y^2}\right) f(x,y) .
\end{equation}
This looks at least as bad as (\ref{eq:inast}). In fact, this sum is as nice
as (\ref{eq:incirc}), in that it can be computed given (\ref{eq:inplus}).
The following observation is crucial:
\begin{equation}\label{eq:orf}\begin{aligned}
\left(\frac{2 x + 7 y}{x^2 + 3 x y + 5 y^2}\right) &= (-1)^{(2 x + 7 y - 1)
(x^2 + 3 x y + 5 y^2)/4} \left(\frac{x^2 + 3 x y + 5 y^2}{2 x + 7 y}\right)\\
&= (-1)^{(2 x + 7 y - 1) (x^2 + 3 x y + 5 y^2)/4} 
\left(\frac{2}{2 x + 7 y}\right)
\left(\frac{2 x^2 + 6 x y + 10 y^2}{2 x + 7 y}\right)\\
&= (-1)^{P_1(x,y)}
\left(\frac{(2 x + 7 y)\cdot x  - x y + 10 y^2}{2 x + 7 y}\right)\\
&= (-1)^{P_1(x,y)} \cdot \left(\frac{- x y + 10 y^2}{2 x + 7 y}\right)\\
&= (-1)^{P_2(x,y)} \cdot \left(\frac{y}{2 x + 7 y}\right) 
\left(\frac{2}{2 x + 7 y}\right) \left(\frac{- 2 x + 20 y}{2 x + 7 y}\right)\\
&= (-1)^{P_3(x,y)} \cdot \left(\frac{2 x + 7 y}{y}\right) 
\left(\frac{2 x + 7 y}{2}\right) \left(\frac{27 y}{2 x + 7 y}\right)\\
&= (-1)^{P_4(x,y)} \left(\frac{2 x}{y}\right) \left(\frac{y}{2}\right)
\left(\frac{2 x + 7 y}{27 y}\right)\\ &= (-1)^{P_5(x,y)}
\left(\frac{x}{y}\right) \left(\frac{x}{y}\right) = (-1)^{P_5(x,y)},
\end{aligned}\end{equation}
where $P_1(x,y) = (2 x + 7 y - 1)
(x^2 + 3 x y + 5 y^2)/4 + ((2 x + 7 y)^2 - 1)/8$, $P_2(x,y)$,\dots,
$P_5(x,y)$ stand for homogeneous polynomials on $x$ and $y$. 
The factors $\left(\frac{x}{y}\right)$ have cancelled out.
Notice that the procedure is basically the Euclidean algorithm.

(Of course, there are parts of (\ref{eq:orf}) that are not quite right -- 
we proceeded as if the expressions inside a Jacobi symbol were always positive
and odd, for instance. We can be careless because this is an introductory
sketch.) 

The crucial thing here is to realise that an expression of the form
 $(-1)^{P_5(x,y)}$ (where $P_5$ is a polynomial) is much better to average against than
$\left(\frac{x}{y}\right)$ is. If we know (\ref{eq:inplus}), then
we can average
\[\sum_{x,y\leq N} (-1)^{P_5(x,y)} f(x,y) ,
\]
and thus we can average (\ref{eq:brascho}), which we have just shown to be the 
same.

In this section, we shall show that we can simplify matters as we just did
whenever we have a Jacobi symbol (or any reciprocity symbol like it) 
of the form $\left(\frac{f(x,y)}{g(x,y)}\right)$, where $f$ and $g$
are homogeneous polynomials and either $\deg(f)$ or $\deg(g)$ is even.
Happily enough, either $\deg(f)$ or $\deg(g)$ will be even in all of our
later applications.

\subsection{Symbols and polynomials}\label{subs:recsec}
We will find it convenient to work abstractly; that is, we will work for
a while with an abstract symbol $(a|b)$ that is assumed to have all of
the properties that a suitably generalised quadratic reciprocity symbol
must have. (It will actually be better to work with a family of
symbols $(a|b)_{\mathfrak{d}}$, rather than a single symbol $(a|b)$.)

Our aim is to show that expressions of the form $(f(x,y)|g(x,y))$ can be
simplified. Here $f$ and $g$ are homogeneous polynomials.

Let $K$ be a global field.
Let $V\subset M_K$ 
be a finite set of places including all the archimedean ones.
Consider a function
\begin{equation}\label{eq:bdb}
(\cdot|\cdot)_{\mathfrak{d}} : \{a,b\in \mathscr{O}_{K,V} : 
b\ne 0, \gcd(a,b) | \mathfrak{d}^{\infty}\} \to \mathbb{C}\end{equation}
for every non-zero ideal $\mathfrak{d}$ of $\mathscr{O}_{K,V}$. Assume
that $(\cdot|\cdot)_{\mathfrak{d}}$ satisfies the following conditions
for some set $\mathscr{C}$ of functions from $\mathscr{O}_{K,V}$ to
$\mathbb{C}$, some set of functions $\mathscr{C}'$ from
$K^* \times K^*$ to $\mathbb{C}$ and all non-zero $a,b,c\in 
\mathscr{O}_{K,V}$:
\begin{enumerate}
\item\label{en:3} $( a b | c ) _{\mathfrak{d}} = (a | c)_{\mathfrak{d}} \cdot 
(b | c)_{\mathfrak{d}},$

\item\label{en:4} $(a | b c)_{\mathfrak{d}} = (a | b)_{\mathfrak{d}} \cdot 
(a | c)_{\mathfrak{d}},$

\item\label{en:5} $(a | b)_{\mathfrak{d}} = (a + b c | b)_{\mathfrak{d}}$,

\item\label{en:7} $(a | b)_{\mathfrak{d}} = f_{\mathfrak{d},b}(a)
\text{\:for $a,b\ne 0$, \;where $f_{\mathfrak{d},b} \in \mathscr{C}$,}$

\item\label{en:8} $(a | b)_{\mathfrak{d}_1} = 
 f_{\mathfrak{d}_1,\mathfrak{d}_2}(a,b) (a | b)_{\mathfrak{d}_2} \text{%
\:for $a,b\ne 0$, $\mathfrak{d}_1|\mathfrak{d}_2$ with
$\gcd(a,b)|\mathfrak{d}_1^{\infty}$, \;where 
$f_{\mathfrak{d}_1,\mathfrak{d}_2} \in \mathscr{C}'$,}$

\item\label{en:8b} $(a | b)_{\mathfrak{d}_1} = 
 g_{\mathfrak{d}_1,\mathfrak{d}_2}(a,b) (a | b) _{\mathfrak{d}_2} \text{%
\:for $a,b\ne 0$, $\mathfrak{d}_2|\mathfrak{d}_1$ with 
$\gcd(a,b)|\mathfrak{d}_2^{\infty}$, \;where 
$g_{\mathfrak{d}_1,\mathfrak{d}_2} \in \mathscr{C}'$,}$

\item\label{en:6} $(a | b)_{\mathfrak{d}} = f_{\mathfrak{d}}(a,b)\cdot 
    (b|a)_{\mathfrak{d}}$ for $a, b\ne 0$,
where $f_{\mathfrak{d}} \in \mathscr{C}'$.
\end{enumerate}
(Each of these equations is meant to hold whenever all symbols on the
right side are defined. For example, (\ref{en:3}) is applicable when
$\gcd(a,c)|\mathfrak{d}^{\infty}$, $\gcd(b,c)|\mathfrak{d}^{\infty}$ and
$c=0$: 
by (\ref{eq:bdb}), this is what must hold for $(a|c)_{\mathfrak{d}}$
and $(b|c)_{\mathfrak{d}}$ to be both defined.)

We will now define a set of functions $\mathscr{D}$ that we will use
in our next result. (We will later be able to show that, 
in the case we care about, all functions
in $\mathscr{D}$ have some rather nice properties.)
Let $\mathscr{D}$ be the set of all complex-valued functions
$f$ on $\mathbb{A}^2(K)$ 
of the form \[\begin{aligned}
f(x,y) &= F_{1,1}(P_1(x,y)) \cdot F_{1,2}(P_2(x,y)) \dotsb
F_{1,k}(P_k(x,y))\\ &\cdot F_{2,1}(Q_1(x,y),R_1(x,y)) \cdot
F_{2,2}(Q_2(x,y),R_2(x,y)) \dotsb F_{2,k'}(Q_{k'}(x,y),R_{k'}(x,y))
,\end{aligned}\]
where $k$, $k'$ are non-negative integers, $F_{1,j}\in \mathscr{C}$,
$F_{2,j}\in \mathscr{C}'$, $P_j\in \mathscr{O}_{K,V}\lbrack x,y\rbrack$
is non-zero and
homogeneous, and, for each $j\leq k'$,
$Q_j, R_j \in \mathscr{O}_{K,V}\lbrack x,y\rbrack$ are non-zero, homogeneous
 and coprime as elements of $K\lbrack x,y\rbrack$.

Since $F_{2,j}(a,b)$ is undefined when $a=0$ or $b=0$, a function
$f(x,y)$ in $\mathscr{D}$ may be undefined when $x/y$ takes a value
within a finite set (namely, the set of roots of $F_{2,j}(t,1)=0$).

We recall that $\Res(f,g)$ stands for the resultant of $f$ and $g$
(\S \ref{subs:husb}).
\begin{prop}
\label{prop:wichtig} Let $K$ be a global field; let $V$ be a finite set of
places of $K$ including all archimedean ones.
 Let $\mathscr{D}$ and $(\cdot |\cdot)$ be as above.
Let $f, g\in \mathscr{O}_{K,V}\lbrack x,y\rbrack$ be non-zero,
homogeneous and coprime as elements of $K\lbrack x,y\rbrack$. 
Let $\mathfrak{d}$ be an ideal of $\mathscr{O}_{K,V}$ such that
$\Res(f,g) | \mathfrak{d}^{\infty}$. Let $\mathfrak{d}_0 |
\mathfrak{d}^{\infty}$.
Then there is a function $h$ in $\mathscr{D}$ such that 
\begin{equation*}
(f(x,y) | g(x,y))_{\mathfrak{d}} = h(x,y) 
  (x|y)_{\mathfrak{d}_0}^{(\deg f) (\deg g)}
\end{equation*}
for all $x,y\in \mathscr{O}_{K,V}$ with $\gcd(x,y)$ dividing $\mathfrak{d}^{\infty}$ and
$x/y$ outside a finite subset of $\mathbb{P}^1(K)$.
\end{prop}
The proof of the proposition mimics the Euclidean algorithm.
\begin{proof}
We may assume $x,y\ne 0$, $f(x,y)\ne 0$, $g(x,y)\ne 0$ throughout, as these
conditions hold for all $x$, $y$ with $x/y$ outside a finite subset
of $\mathbb{P}^1(K)$.
If $\deg(g) = 0$, the result follows from condition (\ref{en:7}). 
If $\deg(f)=0$, the
result follows from (\ref{en:6}) and (\ref{en:7}). 
If $f$ or $g$ is reducible, the statement
follows by (\ref{en:3}) or (\ref{en:4}) 
from cases with lower $\deg(f) + \deg(g)$. If $f$ is
irreducible and $g = c x$, $c$ a non-zero element of $\mathscr{O}_{K,V}$, 
then, by (\ref{en:3}), (\ref{en:4}), (\ref{en:5}), (\ref{en:8b})
 and (\ref{en:6}), 
\begin{equation*}
\begin{aligned} (f(x,y)|g(x,y))_{\mathfrak{d}} &=
(a_0 x^k + a_1 x^{k-1} y + \dotsb + a_k y^k | c x)_{\mathfrak{d}} \\
&=
(f(x,y)|c)_{\mathfrak{d}} \cdot (a_k y^k | x)_{\mathfrak{d}} \\ &=
(f(x,y)|c)_{\mathfrak{d}} \cdot (a_k|x)_{\mathfrak{d}}\cdot 
(y|x)_{\mathfrak{d}}^k \\ &=
(f(x,y)|c)_{\mathfrak{d}} \cdot (a_k|x)_{\mathfrak{d}}
\cdot f_{\mathfrak{d}}^k(y,x) g_{\mathfrak{d},\mathfrak{d}_0}^{k}(x,y)
(x|y)_{\mathfrak{d}_0}^k
\end{aligned}
\end{equation*}
for some $f_{\mathfrak{d}}, g_{\mathfrak{d},\mathfrak{d}_0}\in \mathscr{C}'$.
The result then follows from (\ref{en:7}), the definition
of $\mathscr{D}$ and the case of $(\text{constant} |
x)_{\mathfrak{d}}$, which we have already treated.
 The same works for $f$ irreducible, $g = c y$. The case of $g$
irreducible, $f = \text{$c x$ or $c y$}$ follows from (\ref{en:6}) 
and the above.
For $f$, $g$ irreducible, $\deg(f)<\deg(g)$, we apply (\ref{en:6}).
 We are left with
the case of $f$, $g$ irreducible, $f,g \ne cx, cy$, $\deg(f)\geq \deg(g)$.
Write $f = a_0 x^k + \dotsb + a_k y^k$, $g = b_0 x^l + b_1 x^{l-1} y +
\dotsb + b_l y^l$. Then 
\begin{equation*}
\begin{aligned}
(f(x,y)|g(x,y))_{\mathfrak{d}} &= f_{\mathfrak{d},b_0 \mathfrak{d}}(x,y) 
(f(x,y)|g(x,y))_{b_0 \mathfrak{d}} \\&=
f_{\mathfrak{d},b_0 \mathfrak{d}}(x,y) 
(b_0|g(x,y))_{b_0 \mathfrak{d}} 
(b_0 f(x,y)|g(x,y))_{b_0 \mathfrak{d}} \\ &=
f_{\mathfrak{d},b_0 \mathfrak{d}}(x,y) (b_0|g(x,y))_{b_0 \mathfrak{d}} 
(b_0 f(x,y) - a_0 g(x,y)|g(x,y))_{b_0 \mathfrak{d}}
\end{aligned}
\end{equation*}
for all $(x,y)\in \mathscr{O}_{K,V}^2$ such that
 $\gcd(x,y) | \mathfrak{d}_0$ and $b_0 f(x,y) - a_0 g(x,y)\ne 0$.
The coefficient of $x^k$ in $b_0 f(x,y)
- a_0 g(x,y)$ is zero. Hence $b_0 f(x,y) - a_0 g(x,y)$ is a multiple of $y$.
Since $f$ and $g$ are coprime as elements of $K\lbrack x,y\rbrack$,
the linear combination $b_0 f(x,y) - a_0 g(x,y)$ cannot be the zero
polynomial; either it is a reducible polynomial or it is a non-zero
constant times $y$. 
The former case has already been treated. 
By (\ref{en:6}), the latter case
reduces to $f$ irreducible, $g = c y$, which is a situation we have
already considered.
\end{proof}

\subsection{Conclusions on quadratic reciprocity}\label{subs:conc}
Let $K$ be a global field.
Assume that the characteristic of $K$ is $\ne 2$. Let $V\subset M_K$ be a finite set of places
including all archimedean ones. Define
\begin{equation}  \label{eq:brack}
(a|b)_{\mathfrak{d}} = 
\mathop{\prod_{\mathfrak{p}\in M_K\setminus V}}_{\mathfrak{p}\nmid 2 \mathfrak{d}}
(a / \mathfrak{p})^{v_\mathfrak{p}(b)}, 
\end{equation}
for all non-zero $a,b\in K^*$. For $a=0$, $b\in K^*$, set
$(a|b)_{\mathfrak{d}}=1$.  (For further reference, note that the condition
$\mathfrak{p}\nmid 2 \mathfrak{d}$ in (\ref{eq:brack}) is equivalent to
$\mathfrak{p}\nmid \mathfrak{d}$ when $2|\mathfrak{d}$.)

Here $(a/\mathfrak{p})$ is simply the quadratic
reciprocity symbol, meaning a non-trivial character $(\cdot/\mathfrak{p})$ on
$K_{v_{\mathfrak{p}}}^*$
given by
\begin{equation}\label{eq:dostor}(x/\mathfrak{p}) = \begin{cases} 1 & \text{if $x = y^2$ for some 
$y\in K_{v_{\mathfrak{p}}}^*$,}\\
-1 & \text{otherwise.} \end{cases}\end{equation}

We want to prove that conditions (\ref{en:3})--(\ref{en:6}) in 
subsection \ref{subs:recsec} hold. Let us define the sets
$\mathscr{C}$ and $\mathscr{C}'$ first, since the statements of
conditions (\ref{en:3})--(\ref{en:6}) involve these sets.

We define $\mathscr{C}$ to be
the set of all functions $f:\mathscr{O}_{K,V}\to \mathbb{C}$ of the form
\[f(x) = (x/\mathfrak{p}_1) \cdot (x/\mathfrak{p}_2) \dotsb 
(x/\mathfrak{p}_k),\]
where $k$ is any non-negative integer and $\mathfrak{p}_1$, $\mathfrak{p}_2$,
\dots, $\mathfrak{p}_k$ are finite places of $K$ not in $V$.
We define $\mathscr{C}'$ to be the set of all functions
$g:K^*\times K^*\to \{-1,1\}$ of the form
\begin{equation}\label{eq:kot}
g(x,y) = g_1(x,y) \cdot g_2(x,y) \dotsb g_{k}(x,y),\end{equation}
where, for each $j\leq k$, $v_j$ is a place of $K$
and
 $g_j$ is a function from $K_{v_j}^*/(K_{v_j}^*)^2 \times K_{v_j}^*/(K_{v_j}^*)^2$
to $\{0,-1,1\}$.
\begin{cor}[to Proposition \ref{prop:wichtig}]
\label{cor:important} 
Let $K$ be a global field with $\charac(K)\ne 2$; let $V$ be a finite,
non-empty
 set of places of $K$ including all the archimedean places.
 Let $f,g\in \mathscr{O}_{K,V}\lbrack
x,y\rbrack$ be non-zero, homogeneous, and coprime as elements
of $K\lbrack x,y\rbrack$.
 Let $\mathfrak{d}$ be a 
non-zero ideal of $\mathscr{O}_{K,V}$ such that $\Res(f,g) | \mathfrak{d}^{\infty}$.
Let $\mathfrak{d}_0 | \mathfrak{d}^{\infty}$.
Let $(\cdot | \cdot )_{\mathfrak{d}}$ be as in (\ref{eq:brack}). 
Then, 
Then, for all $x,y \in \mathscr{O}_{K,V}$ with 
$x/y$ outside a finite subset of $\mathbb{P}^1(K)$
and $\gcd(x,y)|\mathfrak{d}^{\infty}$,
\begin{equation}\label{eq:coton}(f(x,y)|g(x,y))_{\mathfrak{d}} = \begin{cases}
h(x,y) \cdot (x|y)_{\mathfrak{d}_0} &\text{if $\deg f$, $\deg g$ are both
odd} \\
h(x,y)&\text{otherwise}\end{cases}\end{equation}
for some complex-valued function $h$ on $\mathbb{A}^2(K)$ of the form
\begin{equation}\label{eq:sdb}
h(x,y) = 
g_1(P_1(x,y),Q_1(x,y)) \cdot g_2(P_2(x,y),Q_2(x,y))\cdot \dotsb \cdot
g_{k}(P_{k}(x,y),Q_{k}(x,y)),
\end{equation}
where, for each $j\leq k$, $v_j$ is a place of $K$,
$g_j$ is a function
from 
$K_{v_j}^*/(K_{v_j}^*)^2 \times K_{v_j}^*/(K_{v_j}^*)^2$ to
$\{-1,1\}$, and
$P_j, Q_j \in \mathscr{O}_{K,V}\lbrack x,y\rbrack$ are non-zero, homogeneous
 and coprime as elements of $K\lbrack x,y\rbrack$.
\end{cor}
\begin{proof}
It is clear that, given our current definitions of $\mathscr{C}$
and $\mathscr{C}'$, the functions $h(x,y)$ of the form above are
simply the elements of the set $\mathscr{D}$ defined in
\S \ref{subs:recsec}. (A function of the form $(P(x,y)/\mathfrak{p})$
can be written in the form $g(P(x,y),Q(x,y))$, where $Q$ is identically
one and $g(x,y) = (x/\mathfrak{p})$.)
For the corollary to follow from Prop.\ 
\ref{prop:wichtig}, it remains only to check that the symbols 
$(a|b)_{\mathfrak{d}}$
defined by (\ref{eq:brack}) satisfy the properties (\ref{en:3})--(\ref{en:6})
listed in \S \ref{subs:recsec}
for all
$a,b\in \mathscr{O}_{K,V}$ with $b\ne 0$, $\gcd(a,b) | \mathfrak{d}^{\infty}$.

(\ref{en:3}): This is true because $a\mapsto (a/\mathfrak{p})$ is
a multiplicative function.

(\ref{en:4}): True by definition (\ref{eq:brack}).

(\ref{en:5}): This is true because $(a/\mathfrak{p}) = (a+x/\mathfrak{p})$
for all $a\notin \mathfrak{p}$, $x\in \mathfrak{p}$, provided that $\mathfrak{p}\nmid 2$. (This
is why we omitted primes $\mathfrak{p}$ dividing
$2$ from the product in (\ref{eq:brack}).)

(\ref{en:7}): This holds because $(a/\mathfrak{p})$ is a character of
the multiplicative group of the residue field of $K_{v_{\mathfrak{p}}}$.

(\ref{en:8}): By (\ref{eq:brack}),
\[\frac{(a|b)_{\mathfrak{d}_1}}{
(a|b)_{\mathfrak{d}_2}} = 
\mathop{\mathop{\prod_{\mathfrak{p}\in M_K \setminus V}}_{\mathfrak{p}|
    2 \mathfrak{d}_2}}_{\mathfrak{p}\nmid 2 \mathfrak{d}_1}
(a/\mathfrak{p})^{v_{\mathfrak{p}}(b)}.\]
Since $(a/\mathfrak{p})$ is the quadratic reciprocity symbol, it
depends on $a\in \mathscr{O}_{K,V}\setminus \{0\}$ 
only as an element of $K_{v_{\mathfrak{p}}}^*/
(K_{v_{\mathfrak{p}}}^*)^2$. 
The valuation $v_{\mathfrak{p}}(b) \mo 2$ depends on $b$ only as an element of
$K_{v_{\mathfrak{p}}}^*/
(K_{v_{\mathfrak{p}}}^*)^2$ (since the valuation of a square is even, i.e.,
$\equiv 0 \mo 2$), and so we are done.

(\ref{en:8b}): Same as (\ref{en:8}).

(\ref{en:6}): This is essentially the quadratic reciprocity law, and thus
it will take the most work to check.
Let $\left(\frac{a,b}{\mathfrak{p}}\right)$ be the quadratic Hilbert
symbol, defined by
\begin{equation}\label{eq:utur}\left(\frac{a,b}{v}\right) = \begin{cases}
1 & \text{if $z^2 = a x^2 + b y^2$ has a non-zero solution $(x,y,z)\in
K_{v} \times K_{v}\times K_{v}$}\\ -1 &\text{otherwise}
\end{cases}\end{equation}
for $a, b\in K_{v}^*$, $v$ a place of $K$. For any
$a,b \in \mathscr{O}_{K,V}$ with $b\ne 0$, $\gcd(a,b)|\mathfrak{d}^{\infty}$,
\[(a|b)_{\mathfrak{d}} 
= \mathop{\prod_{\mathfrak{p}\in M_K\setminus V}}_{\mathfrak{p}\nmid 2 \mathfrak{d}}
(a/\mathfrak{p})^{v_\mathfrak{p}(b)}
= 
\mathop{
\mathop{\prod_{\mathfrak{p}\in M_K\setminus V}}_{\mathfrak{p}\nmid 2
  \mathfrak{d}}}_{
\mathfrak{p}|b}
\left(\frac{b,a}{\mathfrak{p}}\right),\]
since $\left(\frac{b,a}{\mathfrak{p}}\right) = (a/\mathfrak{p})$
for $\mathfrak{p}\nmid a$, $v_{\mathfrak{p}}(b)$ odd, and
$\left(\frac{b,a}{\mathfrak{p}}\right) = 1$ for
$\mathfrak{p}\nmid a$, $v_{\mathfrak{p}}(b)$ even
(see, e.g., \cite[Prop.\ V.3.4]{Ne}).
Similarly,
$(b|a)_{\mathfrak{d}} = 
\prod_{\mathfrak{p}\in M_K\setminus V: \mathfrak{p}\nmid 2 \mathfrak{d}} (b/\mathfrak{p})^{v_{\mathfrak{p}}(a)}
= \prod_{\mathfrak{p}\in M_K\setminus V: \mathfrak{p}\nmid 2 \mathfrak{d},
\mathfrak{p}|a} 
\left(\frac{a,b}{\mathfrak{p}}\right)$.

Thus
\[\begin{aligned}
\frac{(a|b)_{\mathfrak{d}}}{(b|a)_{\mathfrak{d}}}
&= \mathop{\prod_{\mathfrak{p}\in M_K\setminus V}}_{\mathfrak{p}\nmid 2 \mathfrak{d}} (a/\mathfrak{p})^{v_{\mathfrak{p}}(b)}
 \mathop{\prod_{\mathfrak{p}\in M_K\setminus V}}_{\mathfrak{p}\nmid 2
   \mathfrak{d}} 
(b/\mathfrak{p})^{-v_{\mathfrak{p}}(a)} 
= 
\mathop{\mathop{\prod_{\mathfrak{p}\in M_K\setminus V}}_{\mathfrak{p}\nmid 2 \mathfrak{d}}}_{\mathfrak{p}| b} \left(\frac{b,a}{\mathfrak{p}}\right) 
\mathop{\mathop{\prod_{\mathfrak{p}\in M_K\setminus V}}_{\mathfrak{p}\nmid 2 \mathfrak{d}}}_{\mathfrak{p}| a}
\left(\frac{a,b}{\mathfrak{p}}\right)^{-1}\\ 
&=
\mathop{\mathop{\prod_{\mathfrak{p}\in M_K\setminus V}}_{\mathfrak{p}\nmid 2
    \mathfrak{d}}}_{\mathfrak{p}| a b} \left(\frac{a,b}{\mathfrak{p}}\right)^{-1} 
,\end{aligned}\]
where we use the fact that $\left(\frac{a,b}{p}\right)^{-1} = 
\left(\frac{b,a}{p}\right)$. (This is true for all Hilbert symbols
(\cite[Prop.\ V.3.2 (iv)]{Ne}); in the case of the quadratic Hilbert symbol,
which is the one we are treating, this follows easily from
$\left(\frac{a,b}{p}\right) = 
\left(\frac{b,a}{p}\right)$ (directly from the definition)
and $\left(\frac{b,a}{p}\right) = \pm 1$.)

Now, since $\left(\frac{a,b}{p}\right) = 1$ for $\mathfrak{p}$ finite
and $\mathfrak{p}\nmid a,b$ (see, e.g., \cite[Prop.\ V.3.4]{Ne}),
\[
\mathop{\mathop{\prod_{\mathfrak{p}\in M_K\setminus V}}_{\mathfrak{p}\nmid 2
    \mathfrak{d}}}_{\mathfrak{p}| a b} \left(\frac{a,b}{\mathfrak{p}}\right)^{-1} =
\mathop{\prod_{\mathfrak{p}\in M_K\setminus V}}_{\mathfrak{p}\nmid 2
    \mathfrak{d}} \left(\frac{a,b}{\mathfrak{p}}\right)^{-1}.\]

The above is mostly language. Now we will use the product formula
\begin{equation}\label{eq:stski}
\prod_{v\in M_K} \left(\frac{a,b}{v}\right) = 1
\end{equation}
(vd.\ \cite[Thm.\ VI.8.1]{Ne}); this is in some sense the core of 
reciprocity (both quadratic and higher). By (\ref{eq:stski}),
we now have
\begin{equation}\label{eq:ednes}\frac{(a|b)_{\mathfrak{d}}}{(b|a)_{\mathfrak{d}}} = 
\mathop{\prod_{\mathfrak{p}\in M_K\setminus V}}_{\mathfrak{p}\nmid 2
    \mathfrak{d}} \left(\frac{a,b}{\mathfrak{p}}\right)^{-1} =
\prod_{v\in V \cup \{\mathfrak{p}: \mathfrak{p}|2 \mathfrak{d}\}}
\left(\frac{a,b}{\mathfrak{p}}\right)^{-1}.\end{equation}
For every place $v$, the map
$(a,b)\mapsto \left(\frac{a,b}{v}\right)$ is in fact a function
on $K_v^*/{K_v^*}^2 \times K_v^*/{K_v^*}^2$, i.e., it depends on $a$ and
$b$ only modulo $(K_v^*)^2$; this follows from 
the definition (\ref{eq:utur}) of the quadratic Hilbert symbol\footnote{This is
  also true by definition in greater generality, with
$K_v^*/{K_v^*}^2$ replaced by $K_v^*/{K_v^*}^n$. See the definition in
\cite[p.\ 333]{Ne}.}. 
Hence $\frac{(a|b)_{\mathfrak{d}}}{(b|a)_{\mathfrak{d}}}$ is of the form
(\ref{eq:kot}).
\end{proof}

\begin{Rem}
It should be clear from the proof that one should be able to prove similar
statements for higher reciprocity symbols defined in terms of the general
Hilbert symbol for arbitrary $n$ (defined as in \cite{Ne}, Prop.\ V.3).
We shall not need this, and will not elaborate on this further.
\end{Rem}

\begin{Rem}
This is as good a place as any to answer the following question: why
do we insist in statements of the form ``for all
 $(x,y)\in K\times K$ satisfying
 $\gcd(x,y) | \mathfrak{d}_0$'' rather than simply requiring
$x$ and $y$ to be coprime (that is, $\gcd(x,y)=1$)?
The answer is that, if $\mathscr{O}_{K,V}$ is not a unique factorisation
domain, there may be elements
$t\in K$ that cannot be expressed as quotients $t=x/y$ with $x$, $y$
coprime.
\end{Rem}

We will find the following consequence of Corollary \ref{cor:important} to be
particularly useful. We recall that $P_w\in \mathscr{O}_{K,V}\lbrack x,y
\rbrack$ is a homogeneous polynomial associated to a place $w$ of $K(T)$
(see \S \ref{subs:fpp}).
\begin{cor}\label{cor:macizo}
Let $K$ be a global field with $\charac(K)\ne 2$. Let $w$ be a place
of $K(T)$. Let $f\in (K(T))^*$ be a rational function with $w(f)$ even.

Let $V$ be a finite, non-empty set of places of $K$ including all archimedean places.
Let $\mathfrak{d}_0$ be a non-zero ideal of $\mathscr{O}_{K,V}$.
Then, for all $x,y \in \mathscr{O}_{K,V}$ with 
$x/y$ outside a finite subset of $\mathbb{P}^1(K)$ 
and $\gcd(x,y)|\mathfrak{d}^{\infty}$,
the expression $(f(x/y)|P_w(x,y))_{\mathfrak{d}}$ can be written
in the form
\[\begin{aligned}(f(x/y)&|P_w(x,y))_{\mathfrak{d}_0}\\
&= g_1(P_1(x,y),Q_1(x,y)) \cdot g_2(P_2(x,y),Q_2(x,y))\cdot \dotsb 
g_{k}(P_{k}(x,y),Q_{k}(x,y)),\end{aligned}\]
where $v_j$ is a place of $K$,
$g_j$ is a function from 
$K_{v_j}^*/(K_{v_j}^*)^2 \times K_{v_j}^*/(K_{v_j}^*)^2$ to
$\{-1,1\}$, and
and
$P_j, Q_j \in \mathscr{O}_{K,V}\lbrack x,y\rbrack$ are non-zero, homogeneous
 and coprime as elements of $K\lbrack x,y\rbrack$.
\end{cor}
\begin{proof}
We can write 
\[f(x/y) = c\cdot \prod_{w': w'(f)\ne 0} P_{w'}(x,y)^{w'(f)},\]
where $c\in K^*$. 
Write $c = c_0/c_1$, where $c_0, c_1\in \mathscr{O}_{K,V}$.
Let \[\mathfrak{d} = c_0 \cdot c_1 \cdot \mathfrak{d}_0 \cdot
\prod_{w': w'(f)\ne 0} \Res(P_w,P_{w'}).\]

Since $(a|b)_{\mathfrak{d}}$ is multiplicative on $1$,
and equal to $1$ when $a\in (K^*)^2$, 
\[(f(x/y)|P_w(x,y))_{\mathfrak{d}} = 
(c_0|P_w(x,y))_{\mathfrak{d}} (c_1|P_w(x,y))_{\mathfrak{d}}^{-1}
\prod_{w': \text{$w'(f)$ odd}} (P_{w'}(x,y)|P_w(x,y))_{\mathfrak{d}}.\]
Since we are assuming that $w(f)$ is even, every place $w'$ appearing
in $ \prod_{w': \text{$w'(f)$ odd}}$ satisfies $w'\ne w$.

By Corollary \ref{cor:important}, each factor $(P_{w'}(x,y)|P_w(x,y))_{
\mathfrak{d}}$ is of the form
\[
(P_{w'}(x,y)|P_w(x,y))_{
\mathfrak{d}} = 
\begin{cases}
h_{w'}(x,y) \cdot (x|y)_{\mathfrak{d}_0} &\text{if $\deg(w')$, $\deg(w)$ are both
odd} \\
h_{w'}(x,y)&\text{otherwise,}\end{cases}\]
where $h_{w'}$ is as in (\ref{eq:sdb}). Now, since $f$ is an element
of $K(T)$, we know that \[\prod_{w'} w'(f) \deg(w') = 0,\] and so
$\prod_{w': \text{$w'(f)$ odd}} \deg(w')$ is even. Hence, the factors
$(x|y)_{\mathfrak{d_0}}$ that arise in the terms 
$(P_{w'}(x,y)|P_w(x,y))_{\mathfrak{d}}$ cancel out when the product
is taken over all $w'$ with $w'(f)$ odd.

Again by Corollary \ref{cor:important}, $(c_0|P_w(x,y))_{\mathfrak{d}}$
and $(c_1|P_w(x,y))_{\mathfrak{d}}$ are themselves of the form (\ref{eq:sdb}).

It remains to show that 
$(f(x/y)|P_w(x,y))_{\mathfrak{d}_0}/(f(x/y)|P_w(x,y))_{\mathfrak{d}}$ is 
of the form
(\ref{eq:sdb}). Now
\[\frac{(f(x/y)|P_w(x,y))_{\mathfrak{d}_0}}{
(f(x/y)|P_w(x,y))_{\mathfrak{d}}} = 
\mathop{\prod_{\mathfrak{p}\in M_K\setminus V}}_{\mathfrak{p}\nmid 2
\mathfrak{d_0} \wedge \mathfrak{p}|2 \mathfrak{d}}
(f(x/y)/\mathfrak{p})^{v_{\mathfrak{p}}(P_w(x,y))}.\]
This is clearly of the form (\ref{eq:sdb}) with
$P_j(x,y) = P_w(x,y)$ and $Q_j(x,y) = f(x/y)\cdot y^{2\lceil \deg(f)\rceil}$.
\end{proof}

\subsection{Integrals}
We still need to show that expressions such as (\ref{eq:sdb}) are good
to work with. In other words, we need to show that, in Corollary
\ref{cor:important} and \ref{cor:macizo}, we reduced expressions
involving reciprocity symbols to something simpler and more useful,
not to something
more complicated. It is clear that (\ref{eq:sdb}) is a product
over finitely many places, whereas (global) reciprocity symbols are products
over all places of $K$. We must still examine the individual terms of
(\ref{eq:sdb}) and show we can work with them.

Each place $v$ that contributes to (\ref{eq:sdb}) makes a contribution
of the form
\[g_1(P_1(x,y),Q_1(x,y)) \cdot g_2(P_2(x,y),Q_2(x,y))\cdot \dotsb \cdot
g_{r}(P_{r}(x,y),Q_{r}(x,y)),\]
where each $g_j$ is a function from 
$K_{v}^*/(K_{v}^*)^2 \times K_{v}^*/(K_{v}^*)^2$ to
$\{-1,1\}$, and
$P_j, Q_j \in \mathscr{O}_{K,V}\lbrack x,y\rbrack$ are non-zero, homogeneous
 and coprime as elements of $K\lbrack x,y\rbrack$.

\begin{lem}\label{lem:hussite}
Let $K_v$ be a local field. Let $n$ be a positive integer; assume either
$\charac(K_v)=0$ or $n\nmid \charac(K_v)$. Let $g$ be a complex-valued
function on $\mathbb{A}^2(K_v)$ of the form
\[g(x,y) = g_1(P_1(x,y),Q_1(x,y)) \cdot g_2(P_2(x,y),Q_2(x,y))\cdot \dotsb \cdot
g_{r}(P_{r}(x,y),Q_{r}(x,y)),\]
where each $g_j$ is a function from 
$K_{v}^*/(K_{v}^*)^n \times K_{v}^*/(K_{v}^*)^n$ to
$\{-1,1\}$, and
$P_j, Q_j \in K_v\lbrack x,y\rbrack$ are non-zero and homogeneous.

Then $g$ is locally constant at all points where it is defined, i.e.,
it is locally constant on the complement of finitely many lines through
the origin.
\end{lem}
\begin{proof}
 Since 
(by Hensel's lemma) $(K_v^*)^n$ is an open subset of $K_v^*$,
every function from $(K_v^*)/(K_v^*)^n \times
(K_v^*)/(K_v^*)^n$  
to $\mathbb{C}$ induces a function from $K_v^*\times K_v^*$ to 
$\mathbb{C}$ that is everywhere locally constant. A pair of polynomials 
$P,Q\in K_v\lbrack x,y\rbrack$ induce a function 
$(x,y)\mapsto (P(x,y),Q(x,y))$ from
$K_v\times K_v$ to $K_v^*\times K_v^*$ that is defined and locally
constant outside a finite number of lines through the origin (given
by $y= r x$, where $r$ runs through the set of roots to $P(1,r)=0$
and the set of roots to $Q(1,r)=0$). Hence the composition
$(x,y)\mapsto g(P(x,y),Q(x,y))$ is locally constant outside
a finite number of lines through the origin. A finite product of such
functions will itself be locally constant outside a finite number of
lines through the origin.
\end{proof}

We can go further. Let us examine the one-variable case first. The main
results in what follows are Cor.\ \ref{cor:pathan} and Cor.\ \ref{cor:zobra}.

\subsubsection{Lemmas on integration in one variable}

\begin{lem}\label{lem:rotwos}
Let $K_v$ be a local field.  Let $n$ be a positive integer; assume either
$\charac(K_v)=0$ or $n\nmid \charac(K_v)$. Let $P\in K_v\lbrack t\rbrack$
be a non-zero polynomial. Then, for every $t_0\in \mathbb{P}^1(K)$, there
is a punctured neighbourhood
$U_{t_0}$ of $t_0$
such that, for all $t\in U_{t_0}$, the value of $P(t) \mo (K_v^*)^n$
depends only on $(t-t_0)\mo (K_v^*)^n$. (If $t_0 = \infty$, read
$1/t$ instead of $t-t_0$.)
\end{lem}
\begin{proof}
If two polynomials $P_1$ and $P_2$ satisfy the conclusion, so does their
product $P_1\cdot P_2$. We can hence assume that $P$ is irreducible.

We can also assume that we are looking at 
points $t_0$ with $v(t_0)\geq 0$: to look at $t_0$ with $v(t_0)<0$ 
(or at $t_0=\infty$), replace $P(t)$
by $P(1/t) \cdot t^{r n}$, where $r$ is the least integer such that
$r n \geq \deg(P)$.

If $P$ is a constant, then what we seek to prove is trivially true.

Assume now that $P\in K_v\lbrack t\rbrack$ is linear, i.e., of degree $1$. 
Then $P(t) = a t + b$, $a,b\in K_v$. The value of $P(t) \mo\; (K_v^*)^n$ will
depend only on $t - t_0 \mo\; (K_v^*)^n$ for $t_0 = -b/a$. At all points other
than $t_0=-b/a$, the function $t\to P(t) \mo (K_v^*)^n$ is locally constant.

Assume, lastly, that $P\in K_v\lbrack t\rbrack$ is an irreducible
polynomial of degree $2$ or more. By Hensel's lemma, there is an integer $k$
such that $v(P(t))\leq k$ for all $t\in K_v$ with $v(t)\geq 0$. (If
$v(P(t))$ is too large, Hensel's lemma states that there is a root of $P(t)$
in $K_v$ near $t$, and so $P$ would not be irreducible over $K_v$.)
Again by Hensel's lemma, there is an integer $\ell$ such that
$1+x\in (K^*)^n$ for all $x\in K_v$ with $v(x)\geq \ell$. 
Let $t_0\in K_v$ satisfy $v(t_0)\geq 0$. Since $t\mapsto
P(t)$ is a continuous map, there is a 
neighbourhood $U_{t_0}$ of $t_0$ such that $v(t-t_0)\geq \ell+k$ for all
$t\in U_{t_0}$. By the preceding discussion, $P(t) \mo (K_v^*)^n$ is constant
on $U_{t_0}$.
\end{proof}

\begin{cor}\label{cor:marneg}
Let $K_v$ be a local field. Let $n$ be a positive integer; assume either
$\charac(K_v)=0$ or $n\nmid \charac(K_v)$. Let $g$ be a complex-valued
function on $\mathbb{P}^1(K_v)$ of the form
\[g(t) = g_1(P_1(t),Q_1(t)) \cdot g_2(P_2(t),Q_2(t))\cdot \dotsb \cdot
g_{r}(P_{r}(t),Q_{r}(t)),\]
where each $g_j$ is a function from 
$K_{v}^*/(K_{v}^*)^n \times K_{v}^*/(K_{v}^*)^n$ to
$\mathbb{C}$, and
$P_j, Q_j \in \mathscr{O}_{K,V}\lbrack t\rbrack$ are non-zero.

Then, around each $t_0\in \mathbb{P}^1(K_v)$, there is a punctured
neighbourhood $U_{t_0}$ such that, for all $t\in t_0$,
$g(t)$ depends only on $(t-t_0) \mo (K^*)^n$. (If
$t_0=\infty=\mathbb{P}^1(K_v)$, read $1/t$ instead of $t-t_0$.)
\end{cor}
\begin{proof}
Apply Lemma \ref{lem:rotwos} to each $P_j$ and each $Q_j$.
\end{proof}

\begin{cor}\label{cor:pathan}
Let $K_v$ be a local field. Let $n$ be a positive integer; assume either
$\charac(K_v)=0$ or $n\nmid \charac(K_v)$. Let $g$ be a complex-valued
function on $\mathbb{P}^1(K_v)$ of the form
\[g(t) = g_1(P_1(t),Q_1(t)) \cdot g_2(P_2(t),Q_2(t))\cdot \dotsb \cdot
g_{r}(P_{r}(t),Q_{r}(t)),\]
where each $g_j$ is a function from 
$K_{v}^*/(K_{v}^*)^n \times K_{v}^*/(K_{v}^*)^n$ to
$\mathbb{Q}$, and
$P_j, Q_j \in \mathscr{O}_{K,V}\lbrack t\rbrack$ are non-zero.

Then, for every ball $U$ in $K_v$,
\[\int_U g(t) dt\]
is a rational number.
\end{cor}
\begin{proof}
Immediate from Corollary \ref{cor:marneg} and Lemma \ref{lem:griep}.
\end{proof}



\subsubsection{Lemmas on integration in two variables}
\begin{lem}\label{lem:wour}
Let $K_v$ be a local field.  Let $n$ be a positive integer; assume either
$\charac(K_v)=0$ or $n\nmid \charac(K_v)$. Let $P\in K_v\lbrack x,y
\rbrack$ be a non-zero homogeneous polynomial. 
Let $R = \{(x,y)\in K_v\times K_v: \min(v(x),v(y))=0\}$. Let
$\pi_r:R\to K_v$ be the map given by $(x,y)\mapsto x/y$.

Then, for every $t_0\in \mathbb{P}^1(K_v)$, there is a punctured
neighbourhood $U_{t_0}$ of $t_0$ such that, for all
$(x,y)\in \pi_r^{-1}(U_{t_0})$, the value of $P(x,y) \mo (K_v^*)^n$ depends only
on $(x - t_0 y) \mo (K_v^*)^n$ and $y\mo (K_v^*)^n$
(or $x\mo (K_v^*)^n$ and $y \mo (K_v^*)^n$ if $t_0 = \infty$).
\end{lem}
\begin{proof}
Since $P$ is homogeneous, $Q(x/y) = P(x,y) \cdot y^{-\deg(P)}$
is a polynomial on the variable $x/y$. Let $t_0\in \mathbb{P}^1(K_v)$.
 By Lemma \ref{lem:rotwos}, there is a punctured neighbourhood $U_{t_0}$
of $t_0$ such that, for all $(x,y)$ with $\pi_r(x,y)=x/y$ inside
$U_{t_0}$, the value of $Q(x/y) \mo (K^*)^n$ depends only on 
$(x/y - t_0) \mo (K^*)^n$ (or $y/x \mo (K^*)^n$ if $t_0=\infty$).
Hence the value of $P(x,y) \mo (K^*)^n = Q(x/y) y^{-\deg(P)} \mo (K^*)^n$
depends only on $(x/y - y_0) \mo (K^*)^n$ and $y \mo (K^*)^n$
(or $y/x \mo (K^*)^n$ and $y \mo (K^*)^n$ if $t_0=\infty$).
\end{proof}

\begin{cor}\label{cor:marst}
Let $K_v$ be a local field. Let $n$ be a positive integer; assume either
$\charac(K_v)=0$ or $n\nmid \charac(K_v)$. Let $g$ be a complex-valued
function on $\mathbb{A}^2(K_v)$ of the form
\[g(x,y) = g_1(P_1(x,y),Q_1(x,y)) \cdot g_2(P_2(x,y),Q_2(x,y))\cdot \dotsb \cdot
g_{r}(P_{r}(x,y),Q_{r}(x,y)),\]
where each $g_j$ is a function from 
$K_{v}^*/(K_{v}^*)^n \times K_{v}^*/(K_{v}^*)^n$ to
$\mathbb{C}$, and
$P_j, Q_j \in K_v\lbrack x,y\rbrack$ are non-zero and homogeneous.
Let
$\pi_r:R\to K_v$ be the map given by $(x,y)\mapsto x/y$.

Then, for every $t_0\in \mathbb{P}^1(K_v)$, there is a punctured
neighbourhood $U_{t_0}$ of $t_0$ such that, for all
$(x,y)\in \pi_r^{-1}(U_{t_0})$, the value of $P(x,y) \mo (K_v^*)^n$ depends only
on $(x - t_0 y) \mo (K_v^*)^n$ and $y\mo (K_v^*)^n$
(or $x\mo (K_v^*)^n$ and $y \mo (K_v^*)^n$ if $t_0 = \infty$).
\end{cor}
\begin{proof}
Apply Lemma \ref{lem:wour} to each $P_j$ and each $Q_j$.
\end{proof}

\begin{cor}\label{cor:zobra}
Let $K_v$ be a local field. Let $n$ be a positive integer; assume either
$\charac(K_v)=0$ or $n\nmid \charac(K_v)$. Let $g$ be a complex-valued
function on $\mathbb{A}^2(K_v)$ of the form
\[g(x,y) = g_1(P_1(x,y),Q_1(x,y)) \cdot g_2(P_2(x,y),Q_2(x,y))\cdot \dotsb \cdot
g_{r}(P_{r}(x,y),Q_{r}(x,y)),\]
where each $g_j$ is a function from 
$K_{v}^*/(K_{v}^*)^n \times K_{v}^*/(K_{v}^*)^n$ to
$\mathbb{Q}$, and
$P_j, Q_j \in K_v\lbrack x,y\rbrack$ are non-zero and homogeneous.

Let $R = \{(x,y)\in K_v\times K_v: \min(v(x),v(y))=0\}$. Then,
for any balls\, $U_1, U_2 \subset K_v$,
\[\int_{R\cap (U_1\times U_2)} f(x,y)\; dx\, dy\]
is a rational number.
\end{cor}
\begin{proof}
Let
$\pi_r:R\to K_v$ be the map given by $(x,y)\mapsto x/y$.
We can cover $\mathbb{P}^1(K_v)$ by neighbourhoods $U_{t_0} \cup \{t_0\}$,
where $U_{t_0}$ is as in Cor.\ \ref{cor:marst}, and then refine this
covering into a partition of $\mathbb{P}^1(K_v)$ into a finite number
of balls $U_{t_0} \cup \{t_0\}$. It is enough to show that
\begin{equation}\label{eq:yows}
\int_{\pi_r^{-1}(U_{t_0})\cap (U_1\times U_2)} f(x,y) dx dy\end{equation}
is rational for one such ball $U_{t_0}$.
We can assume without loss of generality that $v(t_0)\geq 0$,
and, in particular, that $t_0\ne \infty$. (If
$v(t_0)<0$, switch the variables $x$ and $y$.) We can also assume
that $v(t)=v(t_0)$ for all $t\in U_{t_0}$. (Make all neighbourhoods
$U_{t_0}$ small enough at the beginning.)

As in the proof of Lemma \ref{lem:griep}, it is enough to show that
the area of each set of the form $\pi_r^{-1}((t_0 + g (K^*)^n)\cap U_{t_0})\cap (U_1\times
U_2)$ is rational. We can assume that $U_{t_0}$ is small enough for
$U_{t_0} \cdot U_1$ to be equal to $t_0\cdot U_1$. (For this to be true, it
is enough that the radius of $U_{t_0}$ be smaller than the radius of $U_1$.)
Either the ball $t_0\cdot U_1$ is disjoint from $U_2$ (and we are done),
or $t_0\cdot U_1$ is contained in $U_2$ (and we replace $U_1$ by $t_0^{-1}
U_2$) or $t_0^{-1} \cdot U_2$ is contained in $U_1$ (and we replace
$U_2$ by $t_0 U_2$). In any of these cases, 
we have that it is enough to show that
\[\Area(\pi_r^{-1}((t_0 + g (K^*)^n)\cap U_{t_0})\cap (U_1\times K_v))\]
is rational.

If $v(t_0)>1$,
\[\begin{aligned}
\Area&(\pi_r^{-1}((t_0 + g (K^*)^n)\cap U_{t_0})\cap (U_1\times K_v))\\
&= \Area((t_0 + g (K^*)^n) \cap U_{t_0}) \cdot \Area(U_1 \cap \{x\in K_v:
v(x)\geq 0\}).\end{aligned}\]
If 
\[\begin{aligned}
\Area&(\pi_r^{-1}((t_0 + g (K^*)^n)\cap U_{t_0})\cap (U_1\times K_v))\\
&= \Area((t_0 + g (K^*)^n) \cap U_{t_0}) \cdot \Area(U_1 \cap \{x\in K_v:
v(x)>0\}).\end{aligned}\]
Both $U_1 \cap \{x\in K_v:
v(x)\geq 0\}$ and $U_1\cap \{x\in K_v: v(x)>0\}$ are unions of finitely
many disjoint balls, and thus they have rational areas. We show
as in the proof of 
Lem.\ \ref{lem:griep} that $\Area((t_0 + g (K^*)^n) \cap U_{t_0})$
is rational, and we are done.
\end{proof}

Corollary \ref{cor:pathan} and \ref{cor:zobra}  will all be useful when
we determine the average root number for those rare families where it is
not zero.
\section{The shape of the global root number}\label{sec:globroo}
\subsection{Outline}
Let $\mathscr{E}$ be an elliptic curve over $K(T)$, where $\mathscr{E}$ is
a global field. For simplicity, say $K = \mathbb{Q}$. Let $x,y\in \mathbb{Z}$.
We can write the root number of $\mathscr{E}(x/y)$ as follows:
\[W(\mathscr{E}(x/y)) = - \prod_p W_p(\mathscr{E}(x/y)),\]
where $W_p(\mathscr{E}(x/y))$ is the local root number at $p$. 

Local root
numbers will be described explicitly in Prop.\ \ref{prop:rohr}. As one
can see there, the most interesting case happens when $\mathscr{E}(x/y)$
multiplicative reduction at $p$; the local root number at $p$ is then
\[W_p(\mathscr{E}(x/y)) = - \left(\frac{-c_6(x/y)}{p}\right),\]
where $c_6$ is one of the parameters describing $\mathscr{E}$ in the
standard fashion. (See \S \ref{sss:piu}.) The local root numbers at
places of non-multiplicative reduction are dull in comparison.

Thus, $W(\mathscr{E}(x/y))$ equals something dull times
\begin{equation}\label{eq:litstar}
\mathop{\mathop{\prod_p}_{\text{$\mathscr{E}(x/y)$ has mult.\ }}}_{
\text{reduction at $p$}} W_p(\mathscr{E}(x/y)) =
\mathop{\mathop{\prod_p}_{\text{$\mathscr{E}(x/y)$ has mult.\ }}}_{
\text{reduction at $p$}}  \left(- \left(\frac{-c_6(x/y)}{p}\right)\right),
\end{equation}

Now, the set of primes at which $\mathscr{E}(x/y)$ has multiplicative 
reduction is more or less the same (though not quite the same!) as the set
of primes dividing $M_{\mathscr{E}}(x,y)$, where $M_{\mathscr{E}}\in
\mathbb{Z}\lbrack X,Y\rbrack$ is the product of the polynomials corresponding
to the places $v$ of $K(T)$ where $\mathscr{E}$ has multiplicative
reduction. Hence (\ref{eq:litstar}) equals
\begin{equation}\label{eq:orgo}\begin{aligned}
\text{dull} \cdot \prod_{p|M_{\mathscr{E}}(x,y)}
\left(- \left(\frac{-c_6(x/y)}{p}\right)\right) &=
\text{dull} \cdot \lambda(M_{\mathscr{E}}(x,y)) \cdot
\prod_{p|M_{\mathscr{E}}(x,y)}\left(\frac{-c_6(x/y)}{p}\right)\\
&= \text{dull} \cdot \lambda(M_{\mathscr{E}}(x,y)) \cdot
\left(\frac{-c_6(x/y)}{M_{\mathscr{E}}(x,y)}\right),
\end{aligned}\end{equation}
where $\lambda(n) = \prod_{p|n} (-1)^{v_p(n)}$ is the Liouville function.

We proved in \S \ref{eq:ostorma} that quadratic reciprocity symbols
are dull. Thus, 
\[\text{dull} \cdot \lambda(M_{\mathscr{E}}(x,y)) \cdot
\left(\frac{-c_6(x/y)}{M_{\mathscr{E}}(x,y)}\right) = \text{dull} \cdot \lambda(M_{\mathscr{E}}(x,y)).\]
Therefore
\[W(\mathscr{E}(x/y) = \text{dull} \cdot \lambda(M_{\mathscr{E}}(x,y)).\]
This is the gist of Theorem \ref{thm:smet}, the main result in this section.

\subsection{Preliminaries}

\subsubsection{Valuative criteria for reduction type}\label{sss:piu}
Let $K_v$ be a Henselian
field
 whose residue field has
 characteristic neither $2$ nor $3$. Let $E$ be an elliptic curve over 
$K_v$. Let $c_4, c_6\in K_v$
be a set of parameters describing $E$. As is usual, define the
{\em discriminant} $\Delta$ of (the given model of) $E$ to be
\[\Delta = \frac{c_4^3 - c_6^2}{1728} .\]
Then the reduction of $E$
at $v$ will be
\begin{itemize}
\item {\em good} if $v(c_4) \geq 4 k$, $v(c_6) \geq 6 k$, $v(\Delta) = 12 k$ for some integer $k$;
\item {\em multiplicative} if $v(c_4) = 4 k$, $v(c_6) = 6 k$, 
$v(\Delta) > 12 k$ for some integer $k$;
\item {\em additive} and {\em potentially multiplicative} if 
$v(c_4) = 4 k + 2$, $v(c_6) = 6 k + 3$ and $v(\Delta) > 12 k + 6$ for
some integer $k$;
\item {\em additive} and {\em potentially good} in all remaining cases.
\end{itemize}

It is easy to derive these criteria for reduction type from, say, \cite{Si},
Prop.\ VII.5.1; all we need is a minimal model for $E$. We can find such a 
model as follows. Let $\pi\in \mathscr{O}_{K_v}$ be a uniformiser for $K_v$,
i.e., an element of $K_v$ such that $v(\pi)=1$. Let $k = \min\left(
\left\lfloor \frac{v(c_4)}{4}\right\rfloor,
\left\lfloor \frac{v(c_6)}{6}\right\rfloor\right)$. Since the characteristic
of $K_v$ is $\ne 2,3$, the equation $y^2 = x^3 - 27 c_4 \pi^{-4 k} x - 54 c_6
\pi^{-6k}$ describes an elliptic curve with $c_4$- and $c_6$-parameters equal
to $6^4 \pi^{-4 k} c_4$ and $6^6 \pi^{-6 k} c_6$, respectively. Since the
residue field characteristic of $K_v$ is also $\ne 2,3$, the valuations
of these parameters
are $v(6^4 \pi^{-4k} c_4) = v(c_4) - 4 k$ and $v(6^6 \pi^{-6k} c_6) = 
v(c_6) - 6 k$. By the definition of $k$, either $v(c_4) - 4 k < 4$ or
$v(c_6) - 6 k< 6$. Thus, the equation $y^2 = x^3 - 27 c_4 \pi^{- 4 k}
- 54 c_6 \pi^{-6 k}$ provides a minimal model for $E$. Now apply
\cite{Si},
Prop.\ VII.5.1.

We will say that a curve with bad reduction has {\em half bad} reduction
if a quadratic twist of the curve in question
has good reduction; we will say it has {\em quite bad} reduction 
if no quadratic twist has good 
reduction. It is simple to see that the reduction is {\em half bad} exactly when
$v(c_4) \geq 4 k +2$, $v(c_6) \geq 6 k +3$ and
$v(\Delta)=12 k + 6$ for some integer $k$.

\subsubsection{Local root numbers for $\mathfrak{p}\nmid 2, 3$.
Local root numbers at infinite places}
Recall that we write $(\cdot/\mathfrak{p})$ for the quadratic reciprocity
symbol on $(K_{\mathfrak{p}}^*$; see (\ref{eq:dostor}).
\begin{prop}\label{prop:rohr}
Let $K$ be a global field. Let $\mathfrak{p}$ be a prime of $K$ such that
$\charac(K_{\mathfrak{p}}/\mathfrak{p} K_{\mathfrak{p}})\ne 2,3$.
Let $E$ be an elliptic curve over $K$. If the reduction of $E$ at
$\mathfrak{p}$ is good, then $W_{\mathfrak{p}}(E) = 1$.
 If the reduction 
of $E$ at
$\mathfrak{p}$ is additive and potentially good, then
\begin{enumerate}
\item $W_{\mathfrak{p}}(E) = (-1 / \mathfrak{p})$ if 
$v_{\mathfrak{p}}(\Delta(E))$ is even but not divisible by four,
\item $W_{\mathfrak{p}}(E) = \quadrec{-2}{\mathfrak{p}}$ if 
$v_{\mathfrak{p}}(\Delta(E))$ is odd and divisible by three,
\item $W_{\mathfrak{p}}(E) = \quadrec{-3}{\mathfrak{p}}$ if 
$v_{\mathfrak{p}}(\Delta(E))$ is divisible by four but not by three.
\end{enumerate}

If the reduction of $E$ at $\mathfrak{p}$ is additive and potentially
multiplicative, then 
\begin{equation}\label{eq:totoru}
W_{\mathfrak{p}}(E) = \quadrec{-1}{\mathfrak{p}}.\end{equation}
If the reduction is multiplicative and $c_6$ is the parameter from any
Weierstrass model of $E$, then $W_{\mathfrak{p}}(E) = -1$ when
$-c_6 \in (K_{\mathfrak{p}}^*)^2$ and $W_{\mathfrak{p}}(E) = 1$ when
$-c_6 \notin (K_{\mathfrak{p}}^*)^2$. 
In other words,
when the reduction is multiplicative,
\begin{equation}\label{eq:shackle}
W_{\mathfrak{p}}(E) = - \quadrec{-c_6}{\mathfrak{p}}.\end{equation}
\end{prop}
\begin{proof}
The formulae for good and multiplicative reduction are classical; see, e.g.,
parts (i) and (ii) of the Proposition in \cite{Ro2}, \S 19. 
The case of multiplicative reduction was put in the form (\ref{eq:shackle})
in \cite[Lemma 2.2]{CCH} (but similar expressions were known before).

The formulae
for additive reduction are proved 
in \cite[Prop.\ 2--3]{Ro}, when $K=\mathbb{Q}$,
and in \cite[Prop.\ 2(iii)]{Rog}
for $K$ a number field. (Apply \cite[Prop.\ 2(iii)]{Rog} with
$\tau$ equal to the trivial (one-dimensional) complex representation.) 
Rohrlich's arguments rest on
older work (see \cite{De}) valid for local fields of any
characteristic, and his proofs are general enough to carry over to all
local fields $K_{\mathfrak{p}}$ with residue characteristic $\ne 2,3$.
See the discussion in the proof of \cite{CCH}, Thm.\ 3.1.
\end{proof}

Recall that, for us, a global field is either a number field or a function field
over a finite field.
In the following proposition,
 we assume that $K$ is a number field, as opposed to any global field,
simply because function fields have no infinite places. 
\begin{prop}\label{prop:inflace}
Let $K$ be a number field. Let $E$ be an elliptic curve over $K$.
Let $v$ be an infinite place of $K$.
Then
\[W_v(E) = -1.\]
\end{prop}
\begin{proof}
See \cite[\S 20]{Ro2}.
\end{proof}

\subsubsection{Notation}\label{sss:nota}
Henceforth $K$ will be a global field of characteristic neither $2$ nor
$3$. Fix a finite, non-empty set
of places $V$ of $K$ including all archimedean places.
Let $\mathscr{E}$ be an elliptic curve over $K(T)$ 
given by $c_4, c_6\in K(T)$. The parameters $c_4, c_6 \in K(T)$ give us
the discriminant $\Delta  = \frac{1}{1728} (c_4^3 - c_6^2) \in K(T)$.

Define 
$M_{\mathscr{E}}, B_{\mathscr{E}}\in \mathscr{O}_{K,V}\lbrack x,y\rbrack$ 
as in (\ref{eq:M}). 
(They are the products of the polynomials corresponding to the places
of $K(T)$ where $\mathscr{E}$ has multiplicative or quite bad reduction,
respectively.)
We will henceforth let $(a|b)_{\mathfrak{d}}$ be as in (\ref{eq:brack}). 
(Essentially, $(a|b)_{\mathfrak{d}}$ is a quadratic-reciprocity symbol
that ignores the contributions of the prime ideals dividing $\mathfrak{d}$.)

Let $\mathscr{P}$ be the set of all places $w$ of $K(T)$ such that
$w(c_4)\ne 0$, $w(c_6)\ne 0$ or $w(\Delta)\ne 0$.
We can write
\begin{equation}\label{eq:bahc}\begin{aligned}
c_4(x/y) &= \frac{d_0}{d_1} \prod_{w: w(c_4)\ne 0} (P_w(x,y))^{w(\Delta)}\\
c_6(x/y) &= \frac{d_2}{d_3} \prod_{w: w(c_6)\ne 0} (P_w(x,y))^{w(\Delta)}\\
\Delta(x/y) &= \frac{d_4}{d_5} \prod_{w: w(\Delta)\ne 0} (P_w(x,y))^{w(\Delta)},
\end{aligned}\end{equation}
for some $d_0, d_1, d_2, d_3, d_4, d_5\in \mathscr{O}_{K,V}$, where $P_w$ is the homogeneous
polynomial corresponding to $w$ (see \S \ref{subs:fpp}). 
Let $\mathfrak{d}_{\mathscr{E}}\in I_{K,V}$ be the principal ideal generated by
\begin{equation}\label{eq:toolong}
2 \cdot 3  \cdot d_0 d_1 d_2 d_3 d_4 d_5 \cdot
\mathop{\mathop{\prod_{w_1,w_2\in \mathscr{P}}}_{w_1\ne w_2}}_{
\text{$\mathscr{E}$ has bad red.\ at $w_1$, $w_2$}}
\Res(P_{w_1},P_{w_2}) .\end{equation} 
One may say that $\mathscr{P}$ is the set of interesting places of $K(T)$,
whereas
 $\mathfrak{d}_{\mathscr{E}}$ is the product of the prime ideals 
corresponding to the interesting places of $K$.

\subsection{The contribution of a place of $K(T)$ to $W(\mathscr{E}(t))$}
Let $w$ be a place of $K(T)$.
It is natural to define $W_{\mathscr{E},\mathfrak{d},w}(x,y) =
\prod_{p\nmid \mathfrak{d}: \mathfrak{p}|P_w(x,y)} 
W_{\mathfrak{p}}(\mathscr{E}(x/y))$, which may be thought of as
the contribution to
$W(\mathscr{E}(x/y)) = \prod_{v} W_{v}(\mathscr{E}(x/y))$
coming from the place $w$ of $K(T)$. 
(Here and henceforth, the variable $\mathfrak{p}$ in a product 
$\prod_{\mathfrak{p}}$ is understood to range over the prime ideals
of $\mathscr{O}_{K,V}$.)

\begin{prop}\label{prop:rijk}
Let $K$ be a global field of characteristic $\ne 2,3$; fix a finite, non-empty
set of place $V$ of $K$ including all of the infinite places. Let
$\mathscr{E}$ be an elliptic curve over $K(T)$. Let $w$ be a place of $K(T)$. 
Let $\mathfrak{d}\in I_{K,V}$
be any ideal divisible by $\mathfrak{d}_{\mathscr{E}}$.

Then, for all $x,y \in \mathscr{O}_{K,V}$ with 
$x/y$ outside a finite set of values
and $\gcd(x,y)|\mathfrak{d}^{\infty}$,
\begin{equation}\label{eq:rocco}W_{\mathscr{E}, \mathfrak{d},w}(x,y) = 
\begin{cases} g(x,y) \cdot h(x,y) \cdot \lambda(P_w(x,y))
&\text{if $\mathscr{E}$ has multiplicative reduction at $w$,}\\
g(x,y) \cdot h(x,y) &\text{otherwise,}\end{cases}\end{equation}
where
\begin{enumerate}
\item $\lambda$ is the Liouville function,
\item $g(x,y)$ is of the form
\begin{equation}\label{eq:prodi}g(x,y) = g_1(P_1(x,y),Q_1(x,y)) \cdot g_2(P_2(x,y),Q_2(x,y))\cdot \dotsb \cdot
g_{k}(P_{k}(x,y),Q_{k}(x,y)),
\end{equation}
where each $v_j$ is a place of $K$,
$g_j$ is a function from 
$K_{v_j}^*/(K_{v_j}^*)^2 \times K_{v_j}^*/(K_{v_j}^*)^2$ to
$\{-1,1\}$, and
$P_j, Q_j \in \mathscr{O}_{K,V}\lbrack x,y\rbrack$ are non-zero
and homogeneous,
 and
\item 
$h(x,y)$ is of the form
\begin{equation}\label{eq:squari}h(x,y) = 
\mathop{\prod_{\mathfrak{p}\nmid \mathfrak{d}}}_{\mathfrak{p}^2|P_w(x,y)}
h_{\mathfrak{p}}(c_6(x/y),P_w(x,y)),\end{equation}
where $h_{\mathfrak{p}}$ is a function from 
$K_{\mathfrak{p}}^*/(K_{\mathfrak{p}}^*)^2 \times
K_{\mathfrak{p}}^*/(K_{\mathfrak{p}}^*)^{12}$ to $\{-1,1\}$.
\end{enumerate}
If the reduction of $\mathscr{E}$ at $w$ is half-bad, then $h(x,y)=1$ 
for all values of $x$ and $y$. 
If the reduction of $\mathscr{E}$ at $w$ is good, then
$g(x,y)=1$ and $h(x,y)=1$ for all values of $x$ and $y$.
\end{prop}
\begin{Rem}
Recall that $P_w\in \mathscr{O}_{K,V}\lbrack x,y\rbrack$ is defined to be
a homogeneous polynomial corresponding to a place $w$ of $K(T)$.
For $g_j(P_j(x,y),Q_j(x,y))$ to be well-defined, 
the values $P_j(x,y)$, $Q_j(x,y)$ must lie in $K_{v_j}^*$;
they lie in $K_{v_j}^*$ whenever $x/y$
is outside a finite set of values -- in this case, the set of roots
of $P_j(t,1)$ and $Q_j(t,1)$. The same goes for $c_6(x/y)$ and $P_w(x,y)$;
they lie in $K_{\mathfrak{p}}^*$ whenever $x/y$ is outside a finite set of
values.
\end{Rem}
\begin{proof}
{\em Case 1: $\mathscr{E}$ has good reduction at $w$.}
By the valuative criterion for good reduction (\S \ref{sss:piu}),
we know that $w(c_4) \geq 4 k$, $w(c_6) \geq 6 k$ and 
$w(\Delta) = 12 k$ for some integer $k$, where we see
$c_4$, $c_6$ and $\Delta$ as elements of $K(T)$. 
Let
$\mathfrak{p}|P_w(x,y)$,
 $\mathfrak{p}\nmid \mathfrak{d}$. 
Since 
$\mathfrak{d}_{\mathscr{E}}|\mathfrak{d}$ and $\gcd(x,y)|\mathfrak{d}$,
then, by one of the main properties of the resultant (see \S \ref{subs:husb}),
the ideal $\mathfrak{p}$ cannot divide any $P_{w'}(x,y)$ with 
$w'\in \mathscr{P}$, $w'\ne w$. Thus
\[\begin{aligned}
v_{\mathfrak{p}}(c_4(x/y)) &= w(c_4)\cdot v_{\mathfrak{p}}(P_w(x,y))\geq
4 k \cdot v_{\mathfrak{p}}(P_w(x,y)),\\
v_{\mathfrak{p}}(c_6(x/y)) &= w(c_6)\cdot v_{\mathfrak{p}}(P_w(x,y))\geq
6 k \cdot v_{\mathfrak{p}}(P_w(x,y)),\\
v_{\mathfrak{p}}(\Delta(x/y)) &= w(\Delta)\cdot v_{\mathfrak{p}}(P_w(x,y))
= 12 k \cdot v_{\mathfrak{p}}(P_w(x,y)).
\end{aligned}\]
Thus, again by the criteria in \S \ref{sss:piu}, the elliptic
curve $\mathscr{E}(x/y)$ over $K$ has good reduction at $\mathfrak{p}$.
Hence $W_{\mathfrak{p}}(\mathscr{E}(x/y)) = 1$. We thus have
\[W_{\mathscr{E},\mathfrak{d},w}(x,y)
 = \mathop{\prod_{\mathfrak{p}\nmid \mathfrak{
d}}}_{\mathfrak{p}|P_{w}(x,y)} 
W_{\mathfrak{p}}(\mathscr{E}(x/y)) = 1\]
for all $x,y\in \mathscr{O}_{K,V}$ with $\gcd(x,y)|\mathfrak{d}^{\infty}$.

{\em Case 2: $\mathscr{E}$ has half-bad reduction at $w$.}
By the criterion for half-bad reduction (end of \S \ref{sss:piu}),
we know that $w(c_4) \geq 4 k + 2$, $w(c_6) \geq 6 k + 3$ and 
$w(\Delta) = 12 k + 6$ for some integer $k$.
Let
$\mathfrak{p}|P_w(x,y)$,
 $\mathfrak{p}\nmid \mathfrak{d}$. Since 
$\mathfrak{d}_{\mathscr{E}}|\mathfrak{d}$ and $\gcd(x,y)|\mathfrak{d}$,
$\mathfrak{p}$ cannot divide any $P_{w'}(x,y)$ with $w'\in \mathscr{P}$,
 $w'\ne w$. Thus
\[\begin{aligned}
v_{\mathfrak{p}}(c_4(x/y)) &= w(c_4)\cdot v_{\mathfrak{p}}(P_w(x,y))\geq
(4 k + 2)\cdot v_{\mathfrak{p}}(P_w(x,y)),\\
v_{\mathfrak{p}}(c_6(x/y)) &= w(c_6)\cdot v_{\mathfrak{p}}(P_w(x,y))\geq
(6 k + 3)\cdot v_{\mathfrak{p}}(P_w(x,y)),\\
v_{\mathfrak{p}}(\Delta(x/y)) &= w(\Delta)\cdot v_{\mathfrak{p}}(P_w(x,y))
= (12 k + 6)\cdot v_{\mathfrak{p}}(P_w(x,y)).
\end{aligned}\]
Thus, again by the criteria in \S \ref{sss:piu}, the elliptic
curve $\mathscr{E}(x/y)$ over $K$ has half-bad reduction at $\mathfrak{p}$
if $v_{\mathfrak{p}}(P_w(x,y))$ is odd, and good reduction at $\mathfrak{p}$
if $v_{\mathfrak{p}}(P_w(x,y))$ is even. Hence 
\[W_{\mathfrak{p}}(\mathscr{E}(x/y)) = \begin{cases}
 1 &\text{if $v_{\mathfrak{p}}(P_w(x,y))$ is even,}\\
(-1/\mathfrak{p}) &\text{if $v_{\mathfrak{p}}(P_w(x,y))$ is odd.}\end{cases}
\]
Thereby
\[W_{\mathscr{E},\mathfrak{d},w}(x,y)
 = \mathop{\prod_{\mathfrak{p}\nmid \mathfrak{
d}}}_{\mathfrak{p}|P_{w}(x,y)}
\quadrec{-1}{\mathfrak{p}}^{v_{\mathfrak{p}}(P_{w}(x,y))}
= ( -1 | P_{w}(x,y))_{\mathfrak{d}},.\]
where $(a|b)_{\mathfrak{d}}$ is as in (\ref{eq:brack}). (Recall that 
$2|\mathfrak{d}_{\mathscr{E}}$, and so $2|\mathfrak{d}$.) By
Cor.\ \ref{cor:macizo}, $(-1|P_w(x,y))_{\mathfrak{d}}$ can be written
in the form
\[g_1(P_1(x,y),Q_1(x,y)) \cdot g_2(P_2(x,y),Q_2(x,y))\cdot \dotsb \cdot
g_{k}(P_{k}(x,y),Q_{k}(x,y)),\]
where $g_j$, $P_j$ and $Q_j$ are as in the statement we wish to prove.
Set, then, $g(x,y) = (-1|P_w(x,y))_{\mathfrak{d}}$ and $h(x,y)=1$.

{\em Case 3: $\mathscr{E}$ has multiplicative reduction at $w$.}
From the criteria for reduction type, we know that 
$w(c_4) = 4 k$, $w(c_6) = 6 k$ and 
$w(\Delta) > 12 k$ for some integer $k$.
 Let $\mathfrak{p}|P_w(x,y)$, $\mathfrak{p}\nmid \mathfrak{d}$.
Since $\mathfrak{d}_{\mathscr{E}}|\mathfrak{d}$ and $\gcd(x,y)|\mathfrak{d}$,
the ideal $\mathfrak{p}$ cannot divide any $P_w(x,y)$ with 
$w\in \mathscr{P}$, $w' \ne w$. Hence,
\[\begin{aligned}
v_{\mathfrak{p}}(c_4(x/y)) &= w(c_4)\cdot v_{\mathfrak{p}}(P_w(x,y))
= 4 k \cdot v_{\mathfrak{p}}(P_w(x,y)),\\
v_{\mathfrak{p}}(c_6(x/y)) &= w(c_6)\cdot v_{\mathfrak{p}}(P_w(x,y))
= 6 k \cdot v_{\mathfrak{p}}(P_w(x,y)),\\
v_{\mathfrak{p}}(\Delta(x/y)) &= w(\Delta)\cdot v_{\mathfrak{p}}(P_w(x,y)) >
12 k \cdot v_{\mathfrak{p}}(P_w(x,y)).
\end{aligned}\]
Thus, $\mathscr{E}(x/y)$ has multiplicative reduction at $\mathfrak{p}$. By
Prop.\ \ref{prop:rohr},
\[W_{\mathfrak{p}}(\mathscr{E}(x/y)) = - (-c_6(x/y)/ \mathfrak{p}).\]
Hence (vd.\ (\ref{eq:brack})),
\begin{equation}\label{eq:ciro}\begin{aligned}W_{\mathscr{E},\mathfrak{d},w}(x,y) 
&= \mathop{\prod_{p\nmid \mathfrak{d}}}_{
\mathfrak{p} | P_w(x,y)} (- (-c_6(x,y)/\mathfrak{p})) \\
&= \mathop{\prod_{p\nmid \mathfrak{d}}}_{
\mathfrak{p} | P_w(x,y)} (- (-c_6(x,y)/\mathfrak{p}))^{
v_{\mathfrak{p}}(P_w(x,y))} 
\cdot
\mathop{\prod_{p\nmid \mathfrak{d}}}_{
\mathfrak{p}^2 | P_w(x,y)} (- (-c_6(x,y)/\mathfrak{p}))^{
v_{\mathfrak{p}}(P_w(x,y)) - 1} \\
&= \mathop{\prod_{p\nmid \mathfrak{d}}}_{
\mathfrak{p} | P_w(x,y)} (-1)^{v_{\mathfrak{p}}(P_w(x,y))}
\cdot (-c_6(x,y) | P_w(x,y))_{\mathfrak{d}}
\\ &\cdot
\mathop{\prod_{p\nmid \mathfrak{d}}}_{
\mathfrak{p}^2 | P_w(x,y)} (- (-c_6(x,y)/\mathfrak{p}))^{
v_{\mathfrak{p}}(P_w(x,y)) - 1} .\end{aligned}\end{equation}
Now
\[\mathop{\prod_{p\nmid \mathfrak{d}}}_{
\mathfrak{p} | P_w(x,y)} (-1)^{v_{\mathfrak{p}}(P_w(x,y))} = 
\lambda(P_w(x,y)) \cdot \prod_{\mathfrak{p}|\mathfrak{d}} 
(-1)^{v_{\mathfrak{p}}(P_w(x,y))} .\]
It is clear that $\prod_{\mathfrak{p}|\mathfrak{d}} 
(-1)^{v_{\mathfrak{p}}(P_w(x,y))}$ is of the form (\ref{eq:prodi})
and \[\mathop{\prod_{\mathfrak{p}\nmid
    \mathfrak{d}}}_{\mathfrak{p}^2|P_w(x,y)} (-(-c_6(x,y)/\mathfrak{p}))^{
v_{\mathfrak{p}}(P_w(x,y))-1}\]
is of the form (\ref{eq:squari}). By Cor.\ \ref{cor:macizo},
$(-c_6(x,y)|P_w(x,y))_{\mathfrak{d}}$ is of the form (\ref{eq:prodi})
as well. Thus, 
we obtain from (\ref{eq:ciro})
that $W_{\mathscr{E}, \mathfrak{d},w}(x,y)$ is of the form
(\ref{eq:rocco}).

{\em Case 4: $\mathscr{E}$ has additive, potentially multiplicative reduction
  at $w$.}
$w(c_4) = 4 k + 2$, $w(c_6) = 6 k + 3$ and 
$w(\Delta) > 12 k + 6$ for some integer $k$. Then, for
$\mathfrak{p}|P_w(x,y)$ with
 $\mathfrak{p}\nmid \mathfrak{d}$,
\[\begin{aligned}
v_{\mathfrak{p}}(c_4(x/y)) &= w(c_4)\cdot v_{\mathfrak{p}}(P_w(x,y))=
(4 k + 2)\cdot v_{\mathfrak{p}}(P_w(x,y)),\\
v_{\mathfrak{p}}(c_6(x/y)) &= w(c_6)\cdot v_{\mathfrak{p}}(P_w(x,y))=
(6 k + 3)\cdot v_{\mathfrak{p}}(P_w(x,y)),\\
v_{\mathfrak{p}}(\Delta(x/y)) &= w(\Delta)\cdot v_{\mathfrak{p}}(P_w(x,y))
> (12 k + 6)\cdot v_{\mathfrak{p}}(P_w(x,y)).
\end{aligned}\]
Thus, $\mathscr{E}(x/y)$ has additive, potentially multiplicative
 reduction at $\mathfrak{p}$
if $v_{\mathfrak{p}}(P_w(x,y))$ is odd, and multiplicative reduction
 at $\mathfrak{p}$
if $v_{\mathfrak{p}}(P_w(x,y))$ is even. Hence, by Prop.\ \ref{prop:rohr},
\begin{equation}\label{eq:ostoro}\begin{aligned}W_{\mathscr{E},\mathfrak{d},w}(x,y) &= 
\mathop{\mathop{\prod_{\mathfrak{p}\nmid \mathfrak{d}}}_{\mathfrak{p}|P_{w}(x,y)}}_{\text{$v_{\mathfrak{p}}(P_{w}(x,y))$ odd}}
 (-1/\mathfrak{p})
\cdot
\mathop{\mathop{\prod_{\mathfrak{p}\nmid \mathfrak{d}}}_{\mathfrak{p}|P_{w}(x,y)}
 }_{\text{$v_{\mathfrak{p}}(P_{w}(x,y))$ even}}
(- (-c_6(x,y)/\mathfrak{p}))\\ &=
(-1|P_w(x,y)) \cdot
\mathop{\mathop{\prod_{\mathfrak{p}\nmid \mathfrak{d}}}_{\mathfrak{p}^2|P_{w}(x,y)}
 }_{\text{$v_{\mathfrak{p}}(P_{w}(x,y))$ even}}
(- (-c_6(x,y)/\mathfrak{p})).\end{aligned}\end{equation}
 By Cor.\ \ref{cor:macizo}, $(-1/P_w(x,y))$ is of the form (\ref{eq:prodi}).
It is clear that the product in the last line of (\ref{eq:ostoro}) is
of the form (\ref{eq:squari}).

{\em Case 5: $\mathscr{E}$ has additive, potentially good reduction at $w$
and
$\gcd(w(\Delta),12)=2$.}
 For each 
$\mathfrak{p}\nmid \mathfrak{d}$, we are in one of three different subcases,
depending on $v_{\mathfrak{p}}(P_w(x,y))$:
\[\gcd(v_{\mathfrak{p}}(\Delta(x/y)),12) = \begin{cases}
\text{$2$ or $6$} &\text{if $v_{\mathfrak{p}}(P_w(x,y)) \equiv 1, 3, 5 \mo
  6$,}\\
4 &\text{if $v_{\mathfrak{p}}(P_w(x,y)) \equiv 2, 4 \mo 6$,}\\
12 &\text{if $v_{\mathfrak{p}}(P_w(x,y)) \equiv 0 \mo 6$.}\end{cases}\]
By Prop.\ \ref{prop:rohr}, the three cases result in
$W_{\mathfrak{p}}(\mathscr{E}(x/y)) = (-1/\mathfrak{p})$,
$W_{\mathfrak{p}}(\mathscr{E}(x/y)) = (-3/\mathfrak{p})$ and
$W_{\mathfrak{p}}(\mathscr{E}(x/y)) = 1$, respectively.
Hence
\begin{equation}\label{eq:yotow}\begin{aligned}
W_{\mathscr{E},\mathfrak{d},w}(x,y)
&= \mathop{\prod_{\mathfrak{p}\nmid
\mathfrak{d}}}_{\mathfrak{p}|P_w(x,y)} (-1/\mathfrak{p})^{v_{\mathfrak{p}}(
P_w(x,y))} \cdot \mathop{\prod_{\mathfrak{p}\nmid \mathfrak{d}}}_{
\mathfrak{p}^2 |P_w(x,y)} h_{\mathfrak{p}}(x,y) \\
&= (-1 |P_w(x,y))_{\mathfrak{d}} \cdot 
\mathop{\prod_{\mathfrak{p}\nmid \mathfrak{d}}}_{
\mathfrak{p}^2 |P_w(x,y)} f_{\mathfrak{p}}(v_\mathfrak{p}(P_w(x,y))) ,
\end{aligned}\end{equation}
where $f_{\mathfrak{p}}(a) = (-3 / \mathfrak{p})$ if
$a \equiv 2, 4 \mo 6$, and $f_{\mathfrak{p}}(a) = 1$ otherwise. 
By Cor.\ \ref{cor:macizo}, the term 
$(-1|P_w(x,y))_{\mathfrak{d}}$ is of the form
(\ref{eq:prodi}); it is clear that the product in the last line of
(\ref{eq:yotow}) is of the form (\ref{eq:squari}).

{\em Case 6: $\mathscr{E}$ has additive, potentially good reduction 
at $w$, and $\gcd(w(\Delta),12)=3$.} 
As in case 4, we have three possible subcases for each $\mathfrak{p}\nmid 
\mathfrak{d}$: 
\[\gcd(v_{\mathfrak{p}}(\Delta(x/y)),12) = \begin{cases}
3 &\text{if $v_{\mathfrak{p}}(P_w(x,y)) \equiv 1, 3 \mo 4$,}\\
6 &\text{if $v_{\mathfrak{p}}(P_w(x,y)) \equiv 2 \mo 4$,}\\
12 &\text{if $v_{\mathfrak{p}}(P_w(x,y)) \equiv 0 \mo 4$.}\end{cases}\]
Hence, by Proposition \ref{prop:rohr},
\[W_{\mathfrak{p}}(\mathscr{E}(x/y)) = \begin{cases}
(-2/\mathfrak{p}) &\text{if $v_{\mathfrak{p}}(P_w(x,y)) \equiv 1, 3 \mo 4$,}\\
(-1/\mathfrak{p}) &\text{if $v_{\mathfrak{p}}(P_w(x,y)) \equiv 2 \mo 4$,}\\
1 & \text{if $v_{\mathfrak{p}}(P_w(x,y)) \equiv 0 \mo 4$.}\end{cases}\]
Thus
\begin{equation}\label{eq:sidit}W_{\mathscr{E},\mathfrak{d},w}(x,y) = (-2|P_w(x,y))_{\mathfrak{d}}
\cdot \mathop{\prod_{\mathfrak{p}\nmid \mathfrak{d}}}_{\mathfrak{p}^2 |
P_w(x,y)} h_{\mathfrak{p}}(v_{\mathfrak{p}}(P_w(x,y))),\end{equation}
where $h_{\mathfrak{p}}(a) = (-1/\mathfrak{p})$ if
$a \equiv 2 \mo 4$ and
$f_{\mathfrak{p}}(a) = 1$ otherwise.
The product in (\ref{eq:sidit}) is of the form (\ref{eq:squari}).
By Cor.\ \ref{cor:macizo}, $(-2|P_w(x,y))_{\mathfrak{d}}$ is of the form
(\ref{eq:prodi}).
{\em Case 7: $\mathscr{E}$ has additive, potentially good reduction at
$w$, and
$\gcd(w(\Delta),12)=4$.}
We have two subcases:
Here
\[\gcd(v_{\mathfrak{p}}(\Delta(x/y)),12) = \begin{cases}
4 &\text{if $v_{\mathfrak{p}}(P_w(x,y)) \equiv 1, 2 \mo 3$,}\\
0 &\text{if $v_{\mathfrak{p}}(P_w(x,y)) \equiv 0 \mo 3$.}\end{cases}
\]
Hence $W_{\mathfrak{p}}(\mathscr{E}(x/y)) = (-3/\mathfrak{p})$
if $3\nmid v_{\mathfrak{p}}(P_w(x,y))$ and
$W_{\mathfrak{p}}(\mathscr{E}(x/y)) = 1$ if
$3|v_{\mathfrak{p}}(P_w(x,y))$. Thus
\begin{equation}\label{eq:sirin}
W_{\mathscr{E},\mathfrak{d},w}(x,y) = (-3|P_w(x,y))_{\mathfrak{d}}
\cdot \mathop{\prod_{\mathfrak{p}\nmid \mathfrak{d}}}_{\mathfrak{p}^2 |
P_w(x,y)} f_{\mathfrak{p}}(v_{\mathfrak{p}}(P_w(x,y))),\end{equation}
where $f_{\mathfrak{p}}(a) = (-3/\mathfrak{p})$ if 
$a \equiv 2, 3, 4 \mo 6$, and
$f_{\mathfrak{p}}(a) = 1$ otherwise.
The product in (\ref{eq:sirin}) is of the form (\ref{eq:squari}).
By Cor.\ \ref{cor:macizo}, $(-3|P_w(x,y))_{\mathfrak{d}}$ is of the form
(\ref{eq:prodi}).
\end{proof}

\subsection{The global root number}
Let us first prove a lemma that will come useful soon.
\begin{lem}\label{lem:cagou}
Let $K$ be a global field. Let $n>1$. Assume either
$\charac(K)=0$ or $\gcd(n,\charac(K))=1$.

Let $f$ be a complex-valued function taking values in $\{-1,1\}$. Suppose
that,
for every $t_0\in \mathbb{P}^1(K)$, there
is a punctured neighbourhood $U_{t_0}$ of $t_0$ such that, for $t\in U_{t_0}$,
the value of $f(t)$ depends only on $(t-t_0) \mo (K^*)^n$. (Read 
$1/t$ instead of $(t-t_0)$ if $t_0=\infty \in \mathbb{P}^1(K)$.)

Then there are rational functions $R_j\in (K(T))^*$ 
and functions $g_j:K^*/(K^*)^n\to \{-1,1\}$
such that
\[f(t) = g_1(R_1(t)) \cdot g_2(R_2(t))\cdot  \dotsb \cdot g_k(R_k(t))\]
for all $t\in K$ outside a finite set.
\end{lem}
This is the converse of Corollary \ref{cor:marneg}.
\begin{proof}
We can assume that each $U_{t_0}$ is a punctured ball around $t_0$.
(Otherwise, simply replace each $U_{t_0}$ by a punctured ball around
$t_0$ contained in $U_{t_0}$.)
The neighbourhoods $U_{t_0} \cup \{t_0\}$ are an open cover of
$\mathbb{P}^1(K)$.
Since $\mathbb{P}^1(K)$ is compact, there is an finite subcover.
Given two balls in a non-archimedean field, either they are disjoint
or one of them contains the other one. Hence, our open subcover is
a partition of $\mathbb{P}^1(K)$. It will thus be enough to construct
a function $g:K^*/(K^*)^n\to \{-1,1\}$ and a rational function
$R\in (K(T))^*$ such that $g(R(t))=f(t)$ for all $t\in U_{t_0}$
and $g(R(t))=1$ for all $t\notin U_{t_0} \cup \{t_0\}$. (We can
then take the product of all such $g$ for all $U_{t_0}$ in the 
finite subcover.) We can assume without loss of generality that
$t_0=0$ and $U_{t_0} = U_0 =  \{t\in K^*: v(t)\geq 0\}$.
We can also assume that $f(t)$ is not identically $1$ on $U_0$, as
otherwise the problem is trivial.

Assume first that $f(t)$ is identically $-1$ on $U_0$. Define
$g(x) = -1$ for $x\notin (K^*)^n$ and $g(x) = 1$ for $x\in (K^*)^n$.
By Hensel's lemma,
$(K^*)^n$ is an open subset of $K^*$; in other words, there is an $\ell>0$
such that every element of the form $1+x$, $v(x)\geq \ell$, lies in
$(K^*)^n$.
Let $a$ be any element of $K^*\setminus (K^*)^n$ with $v(a)<-\ell$.
Then $a + x \in K^*\setminus (K^*)^n$ for all $x\in K^*$ with
$v(x)\geq 0$. Let $k$ be the least integer such that $k n + v(a) > \ell$.
Then $a + t^{k n} \in (K^*)^n$ for all $t\in K^*$ with $v(t)<0$.
Define $R(t) = a + t^{k n}$. We have shown that
$f(t) = g(R(t))$ for all $t\in K^*$.

Assume now that $f(t)$ is not identically $-1$ on $U_0$. By hypothesis,
for all $t\in U_0$, $f(t)=h(t)$, where $h$ is a function from $K^*/(K^*)^n$
to $\{-1,1\}$. Let $g(t) = h(t^{-1})$; since $f(t)$ is not identically
$-1$ on $U_0$, $g(t)$ is not identically $-1$. Let $\ell$ be as before,
i.e., a positive integer such that every element of the form $1+x$,
$v(x)\geq \ell$, lies in $(K^*)^n$. Let
$a$ be any element of $K^*$ with $g(a)=1$ and $v(a)> \ell$.
Let $k$ be the least integer such that $k n + 1 > v(a) + \ell$. Then
$g(a+ 1/t^{k n + 1}) = g(a) = 1$ for all $t$ with $v(t)<0$ and
$g(a + 1/t^{k n + 1}) = g(1/t^{k n + 1}) = h(t^{k n + 1}) = h(t)$ 
for all $t$ with $v(t)\geq 0$. 
Define $R(t) = a + 1/t^{k n + 1}$. We have shown that $f(t) = g(R(t))$ for
all $t\in K^*$.
\end{proof}

We finally arrive at one of our main results. It states that the global
root number $W(\mathscr{E}(x/y))$ 
can be expressed as something close to a finite product
times $\lambda(M_{\mathscr{E}}(x,y))$. 

\begin{prop}\label{prop:ahem}
Let $K$ be a global field of characteristic $\ne 2,3$.
Let $\mathscr{E}$ be an elliptic curve over $K(T)$.

Fix a finite, non-empty set $V$ of places of $K$ containing all 
archimedean places;
let $\mathfrak{d}_0$ be any non-zero ideal of $\mathscr{O}_{K,V}$.
Then, for all $x,y \in \mathscr{O}_{K,V}$ with 
$x/y$ outside a finite set of values
and $\gcd(x,y)|\mathfrak{d}_0^{\infty}$,
\begin{equation}\label{eq:exer}
W(\mathscr{E}(x/y)) = g(x,y) \cdot h(x,y)\cdot 
\lambda(M_{\mathscr{E}}(x,y)),
\end{equation}
where 
\begin{enumerate}
\item $\lambda$ is the Liouville function, 
\item $M_{\mathscr{E}}$ and $B_{\mathscr{E}}$, the polynomials corresponding
to the places of multiplicative and quite bad reduction, are as in
(\ref{eq:M}),
\item $g(x,y)$ is of the form
\begin{equation}\label{eq:prodier}g(x,y) = g_1(P_1(x,y),Q_1(x,y)) \cdot g_2(P_2(x,y),Q_2(x,y))\cdot \dotsb \cdot
g_{k}(P_k(x,y),Q_k(x,y)),
\end{equation}
where $v_j$ is a place of $K$,
$g_j$ is a function from 
$K_{v_j}^*/(K_{v_j}^*)^2 \times K_{v_j}^*/(K_{v_j}^*)^2$ to
$\{-1,1\}$, and
$P_j, Q_j \in \mathscr{O}_{K,V}\lbrack x,y\rbrack$ are non-zero
and homogeneous,
 and
\item 
$h(x,y)$ is of the form
\begin{equation}\label{eq:squarier}h(x,y) = 
\mathop{\prod_{\mathfrak{p}\nmid \mathfrak{d}}}_{\mathfrak{p}^2|B_{\mathscr{E}}(x,y)}
\mathop{\prod_{w}}_{\text{$\mathscr{E}$ has q.\ bad red.\ at $w$}} h_{w,\mathfrak{p}}(c_6(x/y),P_w(x,y)),\end{equation}
where $h_{w,\mathfrak{p}}$ is a function from 
$K_{\mathfrak{p}}^*/(K_{\mathfrak{p}}^*)^2 \times
K_{\mathfrak{p}}^*/(K_{\mathfrak{p}}^*)^{12}$ to $\{-1,1\}$.
\end{enumerate}
The functions $g$ and $h$ depend on $\mathscr{E}$, $K$, $V$ and
$\mathfrak{d}$.
\end{prop}
\begin{proof}
Let $\mathfrak{d} = \mathfrak{d}_0 \mathfrak{d}_{\mathscr{E}}$, where
$\mathfrak{d}_{\mathscr{E}}$ is as in (\ref{eq:toolong}).
By Prop.\ \ref{prop:rohr},
$W_{\mathfrak{p}}(\mathscr{E}(x/y))=1$ when $\mathscr{E}(x/y)$ has
good reduction at $\mathfrak{p}$.
Let $\mathfrak{p}\nmid \mathfrak{d}$ be a prime at which $\mathscr{E}(x/y)$
has bad reduction. Then $v_{\mathfrak{p}}(\Delta(x/y))\ne 0$, and so,
by (\ref{eq:bahc}), $\mathfrak{p}|P_w(x,y)$ for some place $w$ such that
$w(\Delta)\ne 0$. Since 
$\mathfrak{d}_{\mathscr{E}}|\mathfrak{d}$ and $\gcd(x,y)|\mathfrak{d}$,
and $\mathfrak{d}_{\mathscr{E}}$ is defined as a product of resultants 
(\ref{eq:toolong}),
the ideal $\mathfrak{p}$ cannot divide any $P_{w'}(x,y)$ with 
$w'\in \mathscr{P}$, $w'\ne w$.
Hence
\begin{equation}\label{eq:lauchun}\begin{aligned}
W(\mathscr{E}(x/y)) &= \prod_{v\in V \cup \{\mathfrak{p}: \mathfrak{p}|
\mathfrak{d}\}} W_v(\mathscr{E}(x/y))
\cdot \mathop{\prod_{\mathfrak{p}\nmid \mathfrak{d}}}_{\text{$\mathscr{E}$
has bad red.\ at $\mathfrak{p}$}}
 W_{\mathfrak{p}}(\mathscr{E}(x/y))\\
&= \prod_{v\in V \cup \{\mathfrak{p}: \mathfrak{p}|\mathfrak{d}\}} 
W_v(\mathscr{E}(x/y)) 
\cdot \prod_{w: w(\Delta)\ne 0} 
\mathop{\prod_{\mathfrak{p}\nmid \mathfrak{d}}}_{\mathfrak{p}|P_w(x,y)}
 W_{\mathfrak{p}}(\mathscr{E}(x/y))\\
&= \prod_{v\in V \cup \{\mathfrak{p}: \mathfrak{p}|
\mathfrak{d}\}} W_v(\mathscr{E}(x/y)) \cdot
\prod_{w: w(\Delta)\ne 0} 
W_{\mathscr{E},\mathfrak{d},w}(x,y).\end{aligned}\end{equation}
The statement now follows
 from Thm.\ \ref{thm:krone}, Prop.\ \ref{prop:inflace} and
Lem.\ \ref{lem:cagou} 
 (applied to the first product in the last line
of (\ref{eq:lauchun})) and from
Prop.\ \ref{prop:rijk} (applied to the second product in the last line
of (\ref{eq:lauchun})). Note that, again by the fact that
$\mathfrak{d}_{\mathscr{E}}$ is a product of resultants, the condition
$\mathfrak{p}^2|B_{\mathscr{E}}(x,y)$ holds for some $\mathfrak{p}\nmid
\mathfrak{d}$ if and only if $\mathfrak{p}^2|P_w(x,y)$ for some 
place $w$ where $\mathscr{E}$ has quite bad reduction.
\end{proof}

The following is a somewhat less explicit and more readable restatement
of Prop.\ \ref{prop:ahem}. It will be quite enough for all of our purposes,
including the proofs of the main theorems.
\begin{thm}\label{thm:smet}
Let $K$ be a global field of characteristic $\ne 2,3$.
Let $\mathscr{E}$ be an elliptic curve over $K(T)$.

Fix a finite, non-empty set $V$ of places of $K$ containing all 
archimedean places;
let $\mathfrak{d}_0$ be any non-zero ideal of $\mathscr{O}_{K,V}$.
Then, for all $x,y \in \mathscr{O}_{K,V}$ with 
$x/y$ outside a finite set of values 
and $\gcd(x,y)|\mathfrak{d}_0^{\infty}$,
\begin{equation}\label{eq:exerier}
W(\mathscr{E}(x/y)) = g(x,y) \cdot h(x,y)\cdot 
\lambda(M_{\mathscr{E}}(x,y)),
\end{equation}
where 
\begin{enumerate}
\item $\lambda$ is the Liouville function, 
\item $M_{\mathscr{E}}$ is the polynomial corresponding
to the set of places of multiplicative reduction (see (\ref{eq:M})),
\item\label{eq:gard} $g(x,y)$ is of the form
\[g(x,y) = \prod_{v\in S} g_v(x,y),\]
where $S$ is a finite set of places of $K$ and 
$g_v:K_v\times K_v\to \{-1,1\}$ is locally constant outside a finite
set of lines through the origin in $K_v\times K_v$, and
\item\label{eq:gord} $h(x,y)$ is of the form
\[h(x,y) = \prod_{\mathfrak{p}
\notin S}
h_{\mathfrak{p}}(x,y),\]
where $h_{\mathfrak{p}}:K_{\mathfrak{p}}\times K_{\mathfrak{p}}\to \{-1,1\}$
is locally constant outside a finite
set of lines through the origin in $K_{\mathfrak{p}}\times K_{\mathfrak{p}}$.
The equality $h_{\mathfrak{p}}(x,y)=1$ holds whenever
$\mathfrak{p}^2|B_{\mathscr{E}}(x,y)$, where $B_{\mathscr{E}}$ is the
polynomial corresponding to the set of places of quite bad reduction
(see (\ref{eq:M})).
\end{enumerate}
\end{thm}
\noindent {\bf Remark.}
The lines through the origin mentioned in (\ref{eq:gard})
are defined over the algebraic closure of $K_v$; the same is true
of the lines in (\ref{eq:gord}) and $K_{\mathfrak{p}}$. Thus, for example,
if $K=\mathbb{Q}$ and $v$ is the infinite place of $\mathbb{Q}$,
the lines have algebraic slopes (though one may be vertical).
 
For every $\mathfrak{p}\notin S$ and every ball $U$
in $K_{\mathfrak{p}}$, $\int_U h_{\mathfrak{p}}(x,1) dx$ is rational.
For every finite 
place $v\in S$ and every ball $U$ in $K_v$, $\int_U g_v(x,1) dx$
is rational. 

For every finite place $v$, let
 $R_v = \{(x,y)\in K_v\times K_v: \min(v(x),v(y))=0\}$. Then,
for every $\mathfrak{p}\notin S$ and any balls
$U_1, U_2 \subset K_v$, the integrals 
$\int_{R \cup (U_1\times U_2)} g_{\mathfrak{p}}(x,y) dx dy$ and
$\int_{R \cup (U_1\times U_2)} h_{\mathfrak{p}}(x,y) dx dy$ are rational.

\begin{proof}[Proof of Theorem \ref{thm:smet}]
Immediate from Proposition \ref{prop:ahem} and Lemma \ref{lem:hussite}.
(A finite number of factors has to be taken from the product defining
$h(x,y)$ and included in the product defining $g(x,y)$.)
The statements in the remark follow easily from Proposition \ref{prop:ahem}
and
Corollaries \ref{cor:pathan} and \ref{cor:zobra}.
\end{proof}

\section{The distribution of the global root number}\label{sec:avcor}
The time has come to compute averages. From now on, we shall be working
with families $\mathscr{E}$ over $\mathbb{Q}(t)$. (A great deal of what
follows probably holds just as well
 for $\mathscr{E}/K(t)$, $K\ne \mathbb{Q}$, under
different conventions on how to average functions from $K$ to $\mathbb{C}$.
See the comments in \S \ref{subs:atrof}.)

We now know what $W(\mathscr{E}(x/y))$ looks like 
(Thm.\ \ref{thm:smet}). Its behaviour as $x$ and $y$ vary depends on
the polynomials 
$B_{\mathscr{E}}$ and $M_{\mathscr{E}}$. If $B_{\mathscr{E}} = 1$ identically,
then (\S \ref{subs:dafs}) the function $t\mapsto W(\mathscr{E}(t))$
is a product of finitely many functions each of which is locally constant
almost everywhere. Finding the average of such a function (over $\mathbb{Z}$
or over $\mathbb{Q}$) is not particularly hard.

If $B_{\mathscr{E}}\ne 1$, the behaviour of $W(\mathscr{E}(t))$ will be 
more complex. If $B_{\mathscr{E}} \ne 1$ and $M_{\mathscr{E}}$ is
identically $1$, then (\S \ref{subs:vogel}) the average of 
$W(\mathscr{E}(t))$ is not necessarily zero. We will be able to compute
this average under the assumption of hypothesis 
$\mathscr{A}_1(B_{\mathscr{E}}(t,1))$ (if the average is over $\mathbb{Z}$)
or $\mathscr{A}_2(B_{\mathscr{E}})$ (if the average is over $\mathbb{Q}$).

Suppose now that $B_{\mathscr{E}}\ne 1$, $M_{\mathscr{E}}\ne 1$ and
$K = \mathbb{Q}$.  (This is the general case.)
Then (\S \ref{subs:edi}) the average of 
$W(\mathscr{E}(t))$ is in fact $0$. This is our main result.
It is conditional on hypotheses $\mathscr{A}_1(B_{\mathscr{E}}(t,1))$
and $\mathscr{B}_1(M_{\mathscr{E}}(t,1))$ (if the average is over
$\mathbb{Z}$) or $\mathscr{A}_2(B_{\mathscr{E}})$
and $\mathscr{B}_2(M_{\mathscr{E}})$ (if the average is over
$\mathbb{Q}$). We have already discussed (\S \ref{subs:htour})
the cases in which these hypotheses have already been proved.
We shall see some families $\mathscr{E}$ for which the main result
is, in consequence, unconditional.
\begin{center} 
* * *
\end{center}
The same procedure that we use to compute averages can be used to
compute autocorrelations, viz.,
\begin{equation}\label{eq:sohn}
\lim_{N\to \infty} \frac{1}{N} \sum_{n\leq N} W(\mathscr{E}(n))
\cdot W(\mathscr{E}(n+1))
\end{equation}
and the like. In particular, if $M_{\mathscr{E}}(t,1)$ is not identically
zero, then (\ref{eq:sohn}) is zero, conditionally on 
$\mathscr{A}_1(B_{\mathscr{E}}(t,1))$ and
$\mathscr{B}_1(M_{\mathscr{E}}(t,1) M_{\mathscr{E}}(t+1,1))$). 
It is hard to tell whether expressions
such as (\ref{eq:sohn}) are of any interest; an application might make the
difference. We shall not discuss (\ref{eq:sohn}) further.

\subsection{Averaging finite products}\label{subs:hoat}
The following propositions will be very useful for finding the averages of
the root numbers in certain special families (\S \ref{subs:dafs}). To deal
with other families, we will later prove auxiliary results that are
elaborations of the ones here.
\begin{prop}\label{prop:cake}
Let $S$ be a finite set of places of $\mathbb{Q}$. For every $v\in S$, let
$g_v:\mathbb{Q}_v\to \mathbb{C}$ be a bounded function that is locally constant
almost everywhere. Then
\begin{equation}\label{eq:dodo}\av_{\mathbb{Z}} \prod_{v\in S} g_v(n) =
c_{\infty} \cdot \prod_{p\in S} \int_{\mathbb{Z}_p} g_p(x) dx ,
\end{equation}
where, if $\infty\notin S$, then
$c_{\infty} = 1$, and, if
$\infty\in S$, then 
$c_{\infty}$ is the value $g_{\infty}(x)$ takes for all $x$ positive and
sufficiently large.
\end{prop}
Here, as always, ``almost everywhere'' means ``outside a finite set of
points''.
Thus, in particular, $g_{\infty}(x)$ is constant for all sufficiently large
positive reals $x$.

In (\ref{eq:dodo}) as elsewhere, $p\in S$ denotes a finite place of $S$,
i.e., a prime; in other words, the product $\prod_{p\in S}$ ranges over the primes
in $S$.

\begin{proof}
Let $|S|$ be the number of elements of $S$. For $p\in S$, let 
$S_p\subset \mathbb{Q}_p$ be the finite set of points $x$ in $\mathbb{Q}_p$
at which $g_p(x)$ is not locally constant. We can cover $S_p$ by arbitrarily
small balls. In particular, we can cover $S_p$ by disjoint open balls whose
union $X_p\subset \mathbb{Z}_p$ has measure $\leq \frac{\epsilon}{|S|}$.
Let $R_p = \mathbb{Z}_p \setminus X_p$. (Here $X_p$ stands for
``exceptional'' and $R_p$ stands for ``regular''.) Then, since
$g_p$ is bounded,
\begin{equation}\label{eq:poro}
\left|\prod_{v\in S} \int_{\mathbb{Z}_p} g_p(x) dx -
\prod_{v\in S} \int_{R_p} g_p(x) dx \right| \leq
|S| \cdot \frac{\epsilon}{|S|} \cdot \max_x \prod_{p\in S} g_p(x)
= O(\epsilon).\end{equation}

Let $U_p$ be a ball in the cover of $S_p$. We can write $U_p = 
a + p^k \mathbb{Z}_p$ for some $a\in \mathbb{Z}$, $k\geq 0$. Now
\begin{equation}\label{eq:oreja}
\mathop{\sum_{n\leq N}}_{n\equiv a \mo p^k} \prod_{v\in S} g_v(n)\; \ll
\mathop{\sum_{n\leq N}}_{n\equiv a \mo p^k} 1\;\leq\; \frac{N}{p^k} + 1 = N \cdot \mu_p(U_p) + 1,
\end{equation}
where $\mu_p$ is the measure on $\mathbb{Z}_p$. Let $\iota_p$ be
the inclusion $\mathbb{Z}\to \mathbb{Z}_p$. Let
\[\rho(n) = \begin{cases} 1 &\text{if $\iota_p(n)\in R_p$
for every $p\in S$}\\ 0&\text{otherwise.}\end{cases}\]
Then, thanks to the fact that $X_p$ is the union of the sets $U_p$ in
(\ref{eq:oreja}),
\begin{equation}\label{eq:chori}\begin{aligned}
\left|\sum_{n\leq N} (1 - \rho(n)) \cdot \prod_{v\in S} g_v(n)\right| &\ll
\mathop{\sum_{n\leq N}}_{\exists p\in S:\; \iota_p(n)\in X_p} 1\\ &\leq \sum_{p\in S} (N \cdot \mu_p(X_p)
+ O(1)) \leq \epsilon \cdot N + o(N),
\end{aligned}\end{equation}
where the implied constants do not depend on $N$.

By (\ref{eq:poro}) and (\ref{eq:chori}),
what remains to show is that
\[\av_{\mathbb{Z}} \rho(n) \cdot \prod_{v\in S} g_v(n) = 
c_{\infty} \cdot \prod_{p\in S} \int_{R_p} g_p(x) dx.\]

Since $R_p$ is compact and $g_p$ is locally constant on $R_p$, we can cover
$R_p$ by
finitely many balls $U_{p,j}$ on each of which $g_p$ is constant. Given two
balls in $\mathbb{Z}_p$, either they are disjoint or one contains the other;
hence, we can take the balls $U_{p,j}$ to be disjoint, i.e., we have a
partition $R_p = \cup_j U_{p,j}$. It will be enough to show that
\[\lim_{N\to \infty} \left(\frac{1}{N} \mathop{\sum_{n\leq N}}_{\forall p\in S:
\iota_p(n)\in U_{p,j_p}} \prod_{v\in S} g_v(n)\right) = c_{\infty} \cdot
\prod_{p\in S} \int_{U_{p,j_p}} g_p(x) dx\]
for every choice of $\vec{j} = \{j_p\}_{p\in S}$. Now, $g_p(x)$ is
constant on $U_{p,j_p}$, and $g_{\infty}(x) = c_{\infty}$ for all
$x\in \mathbb{R}$ larger than a constant. It remains only to prove
\begin{equation}\label{eq:wiggl}
\lim_{N\to \infty} \left( \frac{1}{N} \mathop{\sum_{n\leq N}}_{\forall p\in S:
\iota_p(n)\in U_{p,j_p}} 1\right) = \prod_{p\in S} \mu_p(U_{p,j_p}).
\end{equation}

We can write $U_{p,j_p} = a_p + p^{e_p} \mathbb{Z}_p$, $a_p\in \mathbb{Z}$,
$e_p\geq 0$. Then $\mu_p(U_{p,j_p}) = p^{-e_p}$, and so
\[\prod_{p\in S} \mu_p(U_{p,j_p}) = \prod_{p\in S} p^{-e_p} .
\]
On the other hand, by the Chinese remainder theorem, the integers
$n$ such that $\iota_p(n) \in U_{p,j_p}$ for every $p\in S$ form an arithmetic
progression of modulus $m = \prod_{p\in S} p^{e_p}$. Hence
\[\lim_{N\to \infty} \frac{1}{N} \mathop{\sum_{n\leq N}}_{\forall p\in S:
\iota_p(n) \in U_{p,j_p}} 1 = \frac{1}{m} = \prod_{p\in S} p^{-e_p} .
\]
We have shown (\ref{eq:wiggl}), and so we are done.
\end{proof}

{\bf Remark.} 
The assumptions on $g_p$ in the statement of Prop.\ \ref{prop:cake}
are stronger than they need to be for the lemma to be true. It would
be enough, for example, to assume that $g_p:\mathbb{Q}_p\to
\mathbb{C}$ is a bounded function that is
continuous outside a set of measure zero. (The condition of boundedness
can also be relaxed.) 
Prop.\ \ref{prop:cake} will do for our purposes as is,
however.

Let us now examine averages over $\mathbb{Z}\times \mathbb{Z}$. 
For any function $f:\mathbb{Z}\times \mathbb{Z}\to \mathbb{C}$, define
\begin{equation}\label{eq:avz2def}\av_{\mathbb{Z}^2, \text{coprime}} f(x,y) = \lim_{N\to \infty} 
\frac{\sum_{(x,y)\in \lbrack - N,N\rbrack^2:\; \gcd(x,y)=1} f(x,y)}{
|\{(x,y)\in \lbrack - N, N\rbrack^2: \gcd(x,y)=1\}|}.\end{equation}
The condition of coprimality will be needed in applications: we mean
$x/y$ to go over all rational numbers without repetitions.
Given a function $f:\mathbb{Q}\to \mathbb{C}$, we define
\begin{equation}\label{eq:avqdef}\av_{\mathbb{Q}} f(t) = \lim_{N\to \infty} 
\frac{\sum_{(x,y)\in \lbrack - N,N\rbrack^2: \gcd(x,y)=1, y\ne 0} f(x/y)}{
|\{(x,y)\in \lbrack - N, N\rbrack^2: \gcd(x,y)=1\}|}.\end{equation}
This is a special case of the averages $\av_{\mathbb{Q},S\cap L}$
defined in \S \ref{sec:ostrorta}, 
(\ref{eq:qavy}): we have $\av_{\mathbb{Q}} = 
\av_{\mathbb{Q},S\cap L}$ for $S = \mathbb{R}^2$, $L = \mathbb{Z}^2$.

\begin{prop}\label{prop:torte}
Let $S$ be a finite set of places of $\mathbb{Q}$. For every $v\in S$, let
$g_v:\mathbb{Q}_v\times \mathbb{Q}_v \to \mathbb{C}$ be a bounded function
that is locally constant outside a finite set of lines through the origin.
Then
\[\av_{\mathbb{Z}^2, \text{coprime}} \prod_{v\in S} g_v(x,y) = c_{\infty}\cdot
\prod_{p\in S} \frac{1}{1 - p^{-2}} \int_{O_p} g_p(x,y)\; dx dy,\]
where $O_p = (\mathbb{Z}_p \times \mathbb{Z}_p) \setminus
(p \mathbb{Z}_p \times p \mathbb{Z}_p)$, $c_{\infty} = 1$ if $\infty \notin
S$,
and
\begin{equation}\label{eq:gogo}
c_{\infty} = \lim_{N\to \infty} \frac{1}{(2 N)^2} \int_{-N}^N \int_{-N}^N
g_{\infty}(x,y)\; dx dy
\end{equation}
if $\infty \in S$. 
\end{prop}
Here $1-p^{-2}$ is simply the measure of $O_p$.
(The measure $\mu_p$ of $\mathbb{Q}_p$ 
is, as usual, normalised so that $\mathbb{Z}_p$
has measure $\mu_p(\mathbb{Z}_p)=1$; then the measure 
$(\mu_p\times \mu_p) (\mathbb{Z}_p\times \mathbb{Z}_p)$ is $1$, and so
the measure $(\mu_p\times \mu_p)(O_p)$ of $O_p$ is $1-p^{-2}$.)
\begin{proof}[Sketch of proof]
The argument is the same as in the proof of Prop.\ \ref{prop:cake}; let us
just remark on the differences.

The exceptional set $X_v$ ($v$ infinite or finite) will consist
of small neighbourhoods of the lines through the origin on which
$g_v$ is not locally constant. If a line in
$\mathbb{Q}_v\times \mathbb{Q}_v$ is of the form $y = a x$,
the neighbourhood will be of the form $\{(x,y): |y/x - a|_v < \epsilon\}$,
where $\epsilon$ is small. (The neighbourhood of a vertical line will
be of the form $\{(x,y): |x/y|_v < \epsilon\}$.) In other words,
the neighbourhoods of the lines will be small angles (or, in the
nomenclature we introduced in \S \ref{subs:seclat}, ``sectors'') containing the lines.

Instead of the Chinese remainder theorem, we use the fact that, given
lattice cosets $a_1 + L_1$, $a_2 + L_2$, $a_3 + L_3$,\dots , $a_n + L_k$ 
in $\mathbb{Z}^2$ with indices
 $\lbrack \mathbb{Z}^2 : L_1\rbrack = p_1^{e_1}$,
$\lbrack \mathbb{Z}^2 : L_2\rbrack = p_2^{e_2}$,
$\lbrack \mathbb{Z}^2 : L_3\rbrack = p_3^{e_3}$,\dots ,
$\lbrack \mathbb{Z}^2 : L_k\rbrack = p_k^{e_k}$ 
(where $p_1$, $p_2$, $p_3$,\dots , $p_k$ are distinct), the intersection
\[\bigcap_{1\leq j\leq k} (a_j + L_j)\]
is a lattice coset $a + L$ with index $\lbrack \mathbb{Z}^2:L\rbrack = 
\prod_{j\leq k} p_j^{e_j}$.

We also need the (easy) fact that, given a lattice coset 
$a+L\subset \mathbb{Z}^2$ and a sector $S\subset \mathbb{R}^2$,
\begin{equation}\label{eq:crusty}
\lim_{N\to \infty} \frac{|S \cap (a + L)|}{(2 N)^2} = 
\frac{1}{\lbrack \mathbb{Z}^2 : L\rbrack} \cdot \lim_{N\to \infty} 
\frac{|S \cap \lbrack -N, N\|}{(2 N)^2} .
\end{equation}
\end{proof}

\subsection{Families of quadratic twists}\label{subs:dafs}
Let $\mathscr{E}$ be a family with $B_{\mathscr{E}}=1$. This implies that
$M_{\mathscr{E}}=1$ as well.
Theorem \ref{thm:smet} then gives us a rather simple expression for
$W(\mathscr{E}(x/y))$: the factors $h(x,y)$ and
$\lambda(M_{\mathscr{E}}(x,y))$
in (\ref{eq:exerier})
 become identically $1$. One can actually obtain an even simpler
expression
fairly easily by examining this sort of family from scratch. Let us do so.

By the definition of $B_{\mathscr{E}}$ (see (\ref{eq:M}) in
\S \ref{subs:fpp}), saying that $B_{\mathscr{E}}$ is identically $1$
is the same as saying that there is a quadratic twist $\mathscr{E}'$
of $\mathscr{E}$ having good reduction over every place of 
$K\lbrack t \rbrack$. A curve $\mathscr{E}'$ over $K\lbrack t\rbrack$
has good reduction over every place of $K\lbrack t\rbrack$ if and only
if it is in fact a curve over $K$. Write $\mathscr{E}'$ in the form
$y^2 = x^3 + a x + b$, where $a, b\in K$. Then $\mathscr{E}$ is of the form
\[f(t) y^2 = x^3 + a x + b\]
for some polynomial $f\in K\lbrack T\rbrack$. 

Let us start by examining the case $f(t)=t$.

\begin{prop}\label{prop:pliqua}
Let $K$ be a global field of characteristic $\ne 2,3$. Let $E$ be a fixed
elliptic curve given by a Weierstrass equation $y^2 = x^3 + a x + b$,
$a,b \in K$. Define
\[E_t : t y^2 = x^3 + a x + b \]
for $t\in K^*$.
Then the root number $W(E_t)$ can be written in the form
\begin{equation}\label{eq:utor}W(E_t) = \prod_{v\in S} w_v(t),\end{equation}
where $S$ is a finite set of places of $K$, and
$w_v:K_v^* \to \{-1,1\}$ is such that $w_v(x)$ depends only on $x (K_v^*)^2
\in (K_v^*)/(K_v^*)^2$.
\end{prop}
Note that the factors $w_v$ are {\em not} the same as the local root numbers
$W_v$. Even here, we will have to use quadratic reciprocity to obtain
an expression for the global root number $W(E_t)$ as a finite product.
\begin{proof}
For every place $v$ of $K$ and every $s\in K_v^*$, the curves $E_t$
and $E_{s^2 t}$ are isomorphic over $K^*$. Hence, for every place $v$,
the local root number $W_v(E_t)$ of $E_t$ depends only on 
the image of $t$ in $K_v^*/(K_v^*)^2$. 

Let $S$ be the union of the set of infinite places, the set of primes
whose residue fields have char.\ $2$ or $3$, and the
set of primes dividing the
discriminant of $y^2 = x^3 + a x + b$. By what we said, the product
$\prod_{v\in S} W_v(E_t)$ depends only on the image of $t$ in
$(K_v^*/(K_v^*)^2)_{v\in S}$. It remains to examine 
$\prod_{v\notin S} W_v(E_t)$.

By Prop.\ \ref{prop:rohr},
$\prod_{\mathfrak{p}\notin S} W(\mathscr{E}(t)) = 
 \prod_{\mathfrak{p}\notin S} (-1/\mathfrak{p})^{v_{\mathfrak{p}}(t)}$.
Now
\[\prod_{\mathfrak{p}\notin S} (-1/\mathfrak{p})^{v_{\mathfrak{p}}(t)}
= \prod_{\mathfrak{p}\notin S}
\left(\frac{t, -1}{\mathfrak{p}}\right) =
\prod_{\mathfrak{p}\notin S}
\left(\frac{-1, t}{\mathfrak{p}}\right) = \prod_{v\in S}
\left(\frac{t, -1}{v}\right) ,\]
where $\left(\frac{a,b}{v}\right)$ 
is the quadratic Hilbert symbol.
Since $\left(\frac{a,b}{v}\right)$ depends on $a,b\in K_v^*$ only $\mo
{K_v^*}^2$, the statement follows.

\end{proof}
Proposition \ref{prop:pliqua} actually gives us the behaviour of the root
number for any family $\mathscr{E}$
 of the form $f(T) y^2 = x^3 + a x + b$: equation
(\ref{eq:utor}) tells us that
\begin{equation}\label{eq:dodoro}
W(\mathscr{E}(t)) = W(E_{f(t)}) = \prod_{v\in S} w_v(f(t)) .\end{equation}
This enables us to prove the following.

\begin{cor}\label{cor:albri}
Let $\mathscr{E}$ be an elliptic curve over $\mathbb{Q}(T)$ of the form
\begin{equation}\label{eq:qatra}f(t) y^2 = x^3 + a x + b ,\end{equation}
for some $a,b \in K$, $f\in \mathbb{Q}\lbrack T\rbrack$. 
Then
\begin{equation}\label{eq:fosa}\av_{\mathbb{Z}} W(\mathscr{E}(t)) = 
c_\infty \cdot \prod_{p\in S} \int_{\mathbb{Z}_p} w_p(f(t)),\end{equation}
where $S$ and $w_v$ are as in Prop.\ \ref{prop:pliqua},
and $c_\infty$ is the value taken by $w_\infty(f(x))$ for all sufficiently
large $x$. Moreover,
the integrals $\prod_p \int_{\mathbb{Z}_p} w_p(f(t))$ are rational numbers.
\end{cor}
Since $w_\infty(t)$ depends only on $t \mo (\mathbb{R}^*)^2 = t
\mo \mathbb{R}^+$,
it depends only on the sign $\sgn(t)$ of $t$. Because $\sgn(f(x))$ is
constant for $x$ sufficiently large, it follows that $w_{\infty}(f(x))$
is equal to a constant $c_{\infty}$ for all sufficiently large $x$.
\begin{proof}
By Prop.\ \ref{prop:pliqua} and Corollary \ref{cor:marneg}, the functions
$t\mapsto w_p(f(t))$ are locally constant almost everywhere. Equation
(\ref{eq:fosa}) now follows immediately from Prop.\ \ref{prop:cake}
and Prop.\ \ref{prop:pliqua}.
The integrals $\int_{\mathbb{Z}_p} w_p(f(t))$ are rational by Corollary \ref{cor:pathan}.
\end{proof}

\begin{cor}\label{cor:albor}
Let $\mathscr{E}$ be an elliptic curve over $\mathbb{Q}(T)$ of the form
(\ref{eq:qatra}).
Then
\begin{equation}\label{eq:foson}\av_{\mathbb{Q}} W(\mathscr{E}(x/y)) = 
c_{\infty}\cdot \prod_{p\in S} \frac{1}{1 - p^{-2}} 
\int_{O_p} w_p(f(x/y)) dx dy,\end{equation}
where $S$ and $w_v$ are as in Prop.\ \ref{prop:pliqua},
$R_p = \mathbb{Z}_p\times \mathbb{Z}_p \setminus (p\mathbb{Z}_p \times
p \mathbb{Z}_p)$, and
\begin{equation}\label{eq:woqo}c_{\infty} =
\lim_{N\to \infty} \frac{1}{(2 N)^2} \int_{-N}^N \int_{-N}^N w_{\infty}(x/y) 
dx dy.\end{equation} Moreover,
the integrals $\prod_p \int_{O_p} w_p(f(t))$ are rational numbers,
and $c_{\infty}$ is an algebraic number.
\end{cor}
We define $w_\infty(x)=1$ for all
$x$ if $\infty \notin S$.
\begin{proof}
By Prop.\ \ref{prop:pliqua} and Corollary \ref{cor:marst}, the functions
$t\mapsto w_p(f(t))$ are locally constant almost everywhere. 
Hence
$(x,y)\mapsto w_{\infty}(x/y)$ is a locally
constant function with values in $\{-1,1\}$ defined in the complement
of a finite set of lines through the origin. Equation
(\ref{eq:foson}) now follows from Prop.\ \ref{prop:torte} and Prop.\
\ref{prop:pliqua}.

The integrals $\int_{O_p} w_p(f(x/y)) dx dy$
are rational by Corollary \ref{cor:zobra}. 
Note that the lines on whose complement $(x,y)\mapsto w_{\infty}(x/y)$
is defined are lines with algebraic slopes. (The lines are of the form
$x = r y$, where $r$ goes through the roots of $f(r)=0$.)
It follows that the limit (\ref{eq:woqo}) is an algebraic number.
\end{proof}

Let us finish by showing that a certain phenomenon first noted over
$\mathbb{Q}$ (\cite{cs}) appears over general global fields as well:
a family of quadratic twists can have constant root number.

\begin{cor}\label{cor:ascs}
Let $K$ be a global field of characteristic $\ne 2,3$. Let $E$ be an 
elliptic curve given by a Weierstrass equation $y^2 = x^3 + a x + b$,
$a,b \in K$. Let $E_t$ be as in Proposition \ref{prop:pliqua}. Then there
is a polynomial
 $f\in K\lbrack T\rbrack \setminus (K\lbrack T\rbrack)^2$ such that 
$t\mapsto W(E_{f(t)})$
is a constant map on $K\setminus \{\text{zeroes of $f$ in $K$}\}$.
\end{cor}
Here $K\lbrack T\rbrack \setminus (K\lbrack T\rbrack)^2$ is the set of
elements of $K\lbrack T\rbrack$ not in $(K\lbrack T\rbrack)^2$.
\begin{proof}
By Proposition \ref{prop:pliqua}, it
 will be enough to
construct a polynomial 
$f\in K\lbrack T\rbrack \setminus (K\lbrack T\rbrack)^2$ such that
$f(t) \in {K_v^*}^2$ for every $v\in S$ and every $t\in K$ such that 
$f(t)\ne 0$. (Here $S$ is a finite set of
places of $K$, as in Proposition \ref{prop:pliqua}.)

Suppose first that $K$ is a function field. We are assuming already that
the characteristic of $K$ is not $2$. Thus, the residue field characteristic
of every localisation $K_v$ of $K$ will be different from $2$. Hensel's
Lemma then tells us that $f(t)$ is a square in $K_v^*$
if and only if $v(f(t))$ is even
and $f(t) \pi^{-v(f(t))}$ is a square in the residue field of $K_v$. (Here $\pi$
is any uniformiser of $K_v$.) 

Let $p^n$ be $\max_{v\in S} N v$, where $N v$ is the cardinality of the
residue field of $K_v$. Then $t^{p^n-1} \equiv 1 \mo v$ for every
$v\in S$ for which $v(t) = 0$. Assume first that $p\ne 5$. Then
$f(t) = 3^2 t^{p^n-1} + 4^2$ is the polynomial we desire: 
when $v(t) = 0$, the residue of $f(t) \mo v$ is
a square, namely, $3^2 + 4^2 = 5^2$, and $v(f(t))$ is an even number,
namely, $0$; 
when $v(t)>0$, the residue of $f(t)$ is a square, 
namely, $4^2$, and $v(f(t))$ is again $0$.
When $v(t)<0$, $v(f(t))$ is the even number $(p^n-1) v(t)$,
and the residue of $f(t) \pi^{-v(f(t))}$ is equal to the residue of
$3^2 (t\cdot \pi^{-v(t)})^{p^n-1}$, which is a square.
Thus, for every $t\in K_v^*$, Hensel's Lemma tells us that $f(t)$
is a square in $K_v^*$.
Consider now the case $p=5$. We choose
$f(t) = (9^2 - 7^2)^2 t^{p^n-1} + (2 \cdot 7 \cdot 9)^2$, and proceed
just as for $p\ne 5$.

Suppose now that $K$ is a number field. Define 
$m = 2^k 5^4 \cdot\prod_{\mathfrak{p}\in S} (N\mathfrak{p} - 1)$, 
where we choose the positive integer $k$
to be large enough for the following argument to work. Define
$f(t) = 3^2 t^m + 4^2$. Then, for every non-archimedean place $v$ whose residue
field characteristic is not $2$, $3$ or $5$, the value $f(t)$ is a square 
in $K_v^*$ for every $t\in K$, by the argument used in the function-field case.
The same argument works for residue field characteristic $5$, and
-- provided that $k>1$ -- for characteristic $3$ as well.
In the case of $v$ with residue field characteristic
$2$, Hensel's lemma tells us that $x\in K_v^*$ is a square
whenever $v(x)$ is even and $x \pi^{-v(x)}$ is a square modulo
$\mathfrak{p}^3$,
where $\pi$ is a uniformiser of $v$ and $\mathfrak{p}$ is the prime ideal of
$v$. Since $2^k|m$, we can ensure that
$t^m \pi^{-v(t^m)} = (t \pi^{-v(t)})^m \equiv 1 \mo \mathfrak{p}^l$ for $l$ arbitrarily high
by setting $k$ sufficiently large. If $K=\mathbb{Q}$, then $l=7$ is large
enough for $t^m \pi^{-v(t^m)}\equiv 1 \mo \mathfrak{p}^l$ to imply $f(t)\equiv
1 \mo \mathfrak{p}^3$;
in general, $l = 7 \deg(K/\mathbb{Q})$ is certainly sufficiently large.
We set $k$ high enough for us to obtain this $l$, and conclude
that $f(t)$ is then a square in $K_v^*$.
 It remains to check the archimedean places. For
$v$ real, $f(t) = 3^2 t^m + 4^2$ is always positive, and hence a square
in $\mathbb{R}$. Every number is a square in $\mathbb{C}$, and so 
$f(t) \in {K_v^*}^2$ for $v$ complex.
\end{proof}
\subsection{Averaging almost-finite products}\label{subs:pajaro}
The following propositions will be needed when we find the averages of the
root numbers of families that are not simply quadratic twists of a fixed
curve, yet lack multiplicative reduction.

Recall that, when we write {\bf Lemma 0.0}($\mathscr{X}(P)$),
we mean that the lemma is conditional on hypothesis $\mathscr{X}$ for
the polynomial $P(t)$. All the hypotheses we shall ever refer to were
described in \S \ref{subs:htour}.

\begin{prop}[$\mathscr{A}_1(B)$]\label{prop:cook}
Let $S$ be a finite set of places of $\mathbb{Q}$. For every $v\in S$, let
$g_v:\mathbb{Q}_v\to \mathbb{C}$ be a bounded function that is locally constant
almost everywhere. For every $p\notin S$, let
$h_p:\mathbb{Q}_p\to \mathbb{C}$ be a function that (a) is locally constant
almost everywhere, (b) satisfies $|h_p(x)|\leq 1$ for all $x$. 

Let $B(t)\in \mathbb{Z}\lbrack t\rbrack$ be a 
non-zero polynomial. Assume that $h_p(x)=1$ whenever $v_p(B(x))<2$.
Let
\begin{equation}\label{eq:ostrog}
W(n) = \prod_{v\in S} g_v(n) \cdot \prod_{p\notin S}
h_p(n).\end{equation}

 Then
\begin{equation}\label{eq:dada}
\av_{\mathbb{Z}} W(n) =
c_{\infty} \cdot \prod_{p\in S} \int_{\mathbb{Z}_p} g_p(x) dx 
\cdot \prod_{p\notin S} \int_{\mathbb{Z}_p} h_p(x) dx,
\end{equation}
where, if $\infty\notin S$, then
$c_{\infty} = 1$, and, if
$\infty\in S$, then 
$c_{\infty}$ is the value $g_{\infty}(x)$ takes for all $x$ positive and
sufficiently large.
\end{prop}
\begin{proof}
Let $M$ be a large integer.
We know from Prop.\ \ref{prop:cake} that
\begin{equation}\label{eq:kheer}
\av_{\mathbb{Z}} \prod_{v\in S} g_v(n) \cdot \mathop{\prod_{p\notin S}}_{p<M}
h_v(n)
=
c_{\infty} \cdot \prod_{p\in S} \int_{\mathbb{Z}_p} g_p(x) dx 
\cdot \mathop{\prod_{p\notin S}}_{p<M} \int_{\mathbb{Z}_p} h_p(x) dx,\end{equation}
since all products in the expression are finite.
We must show that each side of (\ref{eq:kheer}) tends to the corresponding
side of (\ref{eq:dada}) as $M\to \infty$.

 For every $p$, the number of solutions
$t\mo p^2$ to $B(t)\equiv 0 \mo p^2$ is bounded independently of $p$
(by Hensel's lemma).
Hence, for $p$ fixed,
the number of integers $n\leq N$ such that $p^2|B(n)$ is
\[\ll \frac{N}{p^2} + 1.\]
It follows that the number of integers $n\leq N$ such that $p^2|B(n)$
for some $p$ between $M$ and $\sqrt{N}$ (inclusive) is
\[\ll \sum_{M\leq p\leq \sqrt{N}} \frac{N}{p^2} + \sum_{M\leq p\leq \sqrt{N}}
1 \ll \frac{N}{M} + \sqrt{N} = \frac{N}{M} + o(N).\]
By hypothesis $\mathscr{A}_1(B(t))$, the number of integers $n\leq N$ such
that $p^2|B(n)$ for some $p>\sqrt{N}$ is $o(N)$. We conclude that the number
of integers $n\leq N$ such that $p^2|B(n)$ for some $p>M$ is
\[O\left(\frac{N}{M}\right) + o(N).\]
Since $h_p(n)=1$ for all $p$ such that $p^2\nmid B(n)$, it follows that
\[
\left|\sum_{n\leq N} \prod_{v\in S} g_v(n) \cdot \prod_{p\notin S}
h_v(n) -
\sum_{n\leq N} \prod_{v\in S} g_v(n) \cdot \mathop{\prod_{p\notin S}}_{p<M}
h_v(n)\right| = O\left(\frac{N}{M}\right) + o(N),\]
and so
\[\left|\av_{\mathbb{Z}} W(n) - av_{\mathbb{Z}}
 \left(\prod_{v\in S} g_v(n) \cdot \mathop{\prod_{p\notin S}}_{p<M}
h_v(n)\right)\right| = O\left(\frac{1}{M}\right).\]
The expression $O(1/M)$ tends to $0$ as $M\to \infty$.

Now let us examine the right side of (\ref{eq:kheer}) and compare it
to the right side of (\ref{eq:dada}). What we need to show is that
\begin{equation}\label{eq:yate}\lim_{M\to \infty} 
\mathop{\prod_{p\notin S}}_{p\geq M} \int_{\mathbb{Z}_p} h_p(x) dx = 1.
\end{equation}
Recall that the number of solutions
$t\mo p^2$ to $B(t)\equiv 0 \mo p^2$ is bounded. Recall also that
$|h_p(x)|\leq 1$ for all $x\in \mathbb{Z}_p$ and
 $h_p(x)=1$ when $v(B(x))<2$. Hence 
\[\int_{\mathbb{Z}_p} h_p(x) = 1 - O\left(\frac{1}{p^2}\right).\]
Therefore
\[\begin{aligned}
\mathop{\prod_{p\notin S}}_{p\geq M} \int_{\mathbb{Z}_p} h_p(x) dx &=
\prod_{p\geq M} \left(1 - O\left(\frac{1}{p^2}\right)\right)\\
&= 1 - O\left(\sum_{m\geq M} \frac{1}{m^2}\right) = 1 - O\left(\frac{1}{M}\right).\end{aligned}\]
The expression $1 - O\left(\frac{1}{M}\right)$ tends to $1$ as $M\to \infty$.
We have shown (\ref{eq:yate}), and thus we are done.
\end{proof}

\begin{prop}[$\mathscr{A}_2(B)$]\label{prop:cookie}
Let $S$ be a finite set of places of $\mathbb{Q}$. For every $v\in S$, let
$g_v:\mathbb{Q}_v\times \mathbb{Q}_v \to \mathbb{C}$ be a bounded function
that is locally constant outside a finite set of lines through the origin.
For every $p\notin S$, let
$h_p:\mathbb{Q}_p\times \mathbb{Q}_p \to \mathbb{C}$ be a function that
(a) is locally constant outside a finite set of lines through the origin,
(b) satisfies $|h_p(x,y)|\leq 1$ for all $x,y\in \mathbb{Q}_p$.

Let $B\in \mathbb{Z} \lbrack x,y\rbrack$ be a non-zero homogeneous 
polynomial. Assume that $h_p(x,y)=1$ whenever $v_p(B(x,y))<2$.
Let
\[
W(x,y) = \prod_{v\in S} g_v(x,y) \cdot \prod_{p\notin S} h_p(x,y).\]

Then
\[\av_{\mathbb{Z}^2, \text{coprime}} W(x,y) = c_{\infty}\cdot
\prod_{p\in S} \frac{1}{1-p^{-2}} \int_{O_p} g_p(x,y)\; dx dy \cdot
\prod_{p\notin S} 
\frac{1}{1 - p^{-2}}
\int_{O_p} h_p(x,y)\; dx dy 
,\]
where $O_p = (\mathbb{Z}_p \times \mathbb{Z}_p) \setminus
(p \mathbb{Z}_p \times p \mathbb{Z}_p)$, $c_{\infty} = 1$ if $\infty \notin
S$,
and
\[
c_{\infty} = \lim_{N\to \infty} \frac{1}{(2 N)^2} \int_{-N}^N \int_{-N}^N
g_{\infty}(x,y)\; dx dy
\]
if $\infty \in S$. 
\end{prop}
\begin{proof}[Sketch of proof]
Proceed just as in the proof of Prop.\ \ref{prop:cook} -- using
Prop.\ \ref{prop:torte} instead of Prop.\ \ref{prop:cake}.
Rather than the fact that the number of solutions 
$t\mo p^2$ to $P(t)\equiv 0 \mo p^2$ ($P$ a non-zero polynomial)
is bounded independently of $p$, we need the fact that the solutions
$(x,y)\in \mathbb{Z}^2$ to $P(x,y)\equiv 0 \mo p^2$ ($P$ a non-zero homogeneous
polynomial) lie on a bounded number of lattices of index $p^2$.
(Both facts are easy consequences of Hensel's
lemma.)
\end{proof}

\subsection{Families without multiplicative reduction}\label{subs:vogel}
In this subsection, we shall 
consider a family $\mathscr{E}$ of elliptic curves over $\mathbb{Q}$
with $M_{\mathscr{E}}=1$, i.e., an elliptic curve $\mathscr{E}$ over
$\mathbb{Q}(t)$ without any places of multiplicative reduction.
We wish to compute the averages
\[\av_{\mathbb{Z}} W(\mathscr{E}(t)),\;\; \av_{\mathbb{Q}} W(\mathscr{E}(t)).
\]

The root number $W(\mathscr{E}(t))$ was described in Theorem \ref{thm:smet}.
Since $M_{\mathscr{E}}(x,y)=1$, there is no factor of the form 
$\lambda(M_{\mathscr{E}}(x,y))$ in (\ref{eq:exerier})
to give us cancellation. Thus, the
averages of $W(\mathscr{E}(x/y))$ will usually be non-zero.

\subsubsection{The average of the root number in families without 
multiplicative reduction}

Recall that $B_{\mathscr{E}}(x,y)$ is the product of the polynomials 
corresponding to the places of quite bad reduction, and that
$M_{\mathscr{E}}(x,y)$ is the product of the polynomials 
corresponding to the places of multiplicative reduction. 
(See (\ref{eq:M}).) We will be examining the case when $M_{\mathscr{E}}=1$,
i.e., the case of families without multiplicative reduction. 

If we are
looking at averages over $\mathbb{Z}$, we only need $M_{\mathscr{E}}(x,1)=1$
for all $x$, as opposed to $M_{\mathscr{E}}(x,y)=1$ for all $x$, $y$.
The families with $M_{\mathscr{E}}(x,1)=1$ for all $x$ are the families
that have no places of multiplicative reduction other than, possibly,
the place corresponding to the valuation $\deg(\den)-\deg(\num)$.

\begin{prop}[$\mathscr{A}_1(B_{\mathscr{E}}(x,1))$]\label{prop:ant1}
Let $\mathscr{E}$ be an elliptic curve over $\mathbb{Q}(T)$. Assume that
$\mathscr{E}$ has no places of multiplicative reduction over $\mathbb{Q}(T)$,
other than, possibly, the place of $\mathbb{Q}(T)$
 corresponding to the valuation
$\deg(\den) - \deg(\num)$. 

Then
\begin{equation}\label{eq:brah1}
\av_{\mathbb{Z}} W(\mathscr{E}(t)) = 
c_{\infty} \cdot \prod_{p\in S} \int_{\mathbb{Z}_p} g_p(x,1) dx 
\cdot \prod_{p\notin S} \int_{\mathbb{Z}_p} h_p(x,1) dx,
\end{equation}
where $S$, $g_v$ and $h_p$ are as in Thm.\ \ref{thm:smet} and
$c_\infty$ is the value taken by $g_{\infty}(x,1)$ for all sufficiently large
$x$. Moreover, the integrals 
$\int_{\mathbb{Z}_p} g_p(x,1) dx$ and
$\int_{\mathbb{Z}_p} h_p(x,1) dx$ are rational numbers.
\end{prop}
\begin{Rem}
Here as in Prop.\ \ref{prop:cook}, the infinite product $\prod_{p\notin S} 
\int_{\mathbb{Z}_p} h_p(x,1) dx$
is absolutely convergent. (In other words, it is like $\prod_p (1-1/p^2)$,
not like $\prod_p (1-1/p)$.)

This is so because (as stated in Thm.\ \ref{thm:smet}) $h_p(x,1)=1$
whenever $p^2\nmid B_{\mathscr{E}}(x,1)$. For every prime larger than
a constant, the congruence $B_{\mathscr{E}}(x,1) \equiv 0 \mo p^2$
can hold only for $O(1)$ congruence classes $x \mo p^2$; hence
the integral $\int_{\mathbb{Z}_p} h_p(x,1) dx$ is $1 + O(1/p^2)$. Thus,
the infinite product is absolutely convergent, just like
$\prod_p (1+1/p^2)$ or $\prod_p (1-1/p^2)$.
\end{Rem}
\begin{proof}
By Theorem \ref{thm:smet}, the functions $x\mapsto g_v(x,1)$ and
$x\mapsto h_p(x,1)$ are locally constant almost everywhere. Moreover, Theorem
\ref{thm:smet}
states that $h_p(x)=1$ whenever $v_p(B(x))<2$, and $|h_p(x)|\leq 1$ otherwise.
Apply Prop.\ \ref{prop:cook}. The integrals 
$\int_{\mathbb{Z}_p} g_p(x,1) dx$ and
$\int_{\mathbb{Z}_p} h_p(x,1) dx$ are rational by Cor.\ \ref{cor:pathan}.
\end{proof}

\begin{prop}[$\mathscr{A}_2(B_{\mathscr{E}})$]\label{prop:ant2}
Let $\mathscr{E}$ be an elliptic curve over $\mathbb{Q}(T)$. Assume that
$\mathscr{E}$ has no places of multiplicative reduction over $\mathbb{Q}(T)$.

Then
\begin{equation}\label{eq:brah2}
\av_{\mathbb{Q}} W(x,y) = c_{\infty}\cdot
\prod_{p\in S} \frac{1}{1-p^{-2}} \int_{O_p} g_p(x,y)\; dx dy \cdot
\prod_{p\notin S} \frac{1}{1 - p^{-2}} \int_{O_p} h_p(x,y)\; dx dy 
,\end{equation}
where $S$, $g_v$ and $h_p$ are as in Thm.\ \ref{thm:smet}, 
$O_p = (\mathbb{Z}_p \times \mathbb{Z}_p) \setminus
(p \mathbb{Z}_p \times p \mathbb{Z}_p)$, 
$c_{\infty} = 1$ if $\infty \notin
S$,
and
\begin{equation}\label{eq:findem}
c_{\infty} = \lim_{N\to \infty} \frac{1}{(2 N)^2} \int_{-N}^N \int_{-N}^N
g_{\infty}(x,y)\; dx dy
\end{equation}
if $\infty \in S$. 
Moreover, the integrals $\prod_{p\in S} \int_{O_p} g_p(x,y)\; dx dy$,
$\prod_{p\notin S} \int_{O_p} h_p(x,y)\; dx dy$ are rational,
and $c_{\infty}$ is an algebraic number. 
\end{prop}
\begin{Rem}
Here, much as in Prop.\ \ref{prop:cookie}, the infinite product $\prod_{p\notin S} 
\int_{O_p} h_p(x,y)\; dx dy$
is absolutely convergent. 
\end{Rem}
\begin{proof}
Immediate from Theorem \ref{thm:smet} and Prop.\ \ref{prop:cookie}.

The integrals $\prod_{p\in S} \int_{O_p} g_p(x,y)\; dx dy$,
$\prod_{p\notin S} \int_{O_p} h_p(x,y)\; dx dy$ are rational by Cor.\
\ref{cor:zobra}. The function $g_{\infty}(x,y)$ is locally constant
outside the union of finitely many lines through the origin; by the remark
after Theorem \ref{thm:smet}, these lines have algebraic slopes. Since
$g_{\infty}(x,y)$ takes only the values $\{-1, 1\}$, 
the integral (\ref{eq:findem}) defining $c_{\infty}$ is algebraic.
\end{proof}

\subsubsection{Averages over function fields}
The analogues of $\mathscr{A}_1$ and $\mathscr{A}_2$ over function
fields are known\footnote{The assumption of separability in
\cite{Ra} is made unnecessary by the argument in \cite{Po}, \S 7.
The same argument suffices to show that there is no
contribution from inseparable places to the infinite products below:
given an irreducible, inseparable polynomial $f\in F_q
\lbrack u\rbrack\lbrack t\rbrack$, there are only finitely
many irreducibles $\pi\in F_q\lbrack u\rbrack$ such that
$\pi^2|f(t)$ has solutions in $F_q\lbrack u\rbrack$. Various other
proofs
of the same are possible \cite{Co}.}
(\cite{Ra}, \cite{Po}). Thus, the function-field analogues of Propositions \ref{prop:ant1}
and \ref{prop:ant2} should be unconditional. The most natural ordering
here for the purposes of averaging is probably simply the ordering
by degree: given a field $K = \mathbb{F}_q(T)$ with integer ring
$\mathscr{O} = \mathbb{F}_q\lbrack T\rbrack$, and given
a function $f:K\to \mathbb{C}$,
define
\[\av_{\mathscr{O}} f(x) = \lim_{n\to \infty} \frac{1}{q^{n+1}}
\mathop{\sum_{x \in \mathbb{F}_q\lbrack T\rbrack}}_{\deg(x) \leq n} f(x)\]
and
\[\av_{K} f(t) = \lim_{n\to \infty} \frac{
\sum_{x,y \in \mathbb{F}_q\lbrack T\rbrack:\; 
\deg(x),\; \deg(y)\leq n,\; \text{$x$, $y$ coprime}}
 f(x/y)}{
|\{x,y \in \mathbb{F}_q\lbrack T\rbrack: 
\deg(x),\; \deg(y)\leq n,\; \text{$x$, $y$ coprime}\}|} .
\]

With these definitions, the proofs
of the analogues of the results in \S \ref{subs:pajaro} and \S \ref{subs:vogel}
should then go through much as before.
\subsection{Products with cancellation}\label{subs:proca}
The following propositions will be needed when we find the averages of the
root numbers of families with multiplicative reduction.

The main ideas of the proofs are already contained in Propositions \ref{prop:cake}
and \ref{prop:cook}, which address a somewhat different situation. 
We will work in a little more generality than there,
in that we will determine averages over arithmetic progressions and over
lattices, rather than just over $\mathbb{Z}$ and $\mathbb{Q}$. (We could
have done the same there, though essentially the same goal would be fulfilled
by a reparametrisation.)

\begin{prop}[$\mathscr{A}_1(B)$]\label{prop:paneton}
Let $S$ be a finite set of places of $\mathbb{Q}$. For every $v\in S$, let
$g_v:\mathbb{Q}_v\to \mathbb{C}$ be a bounded function that is locally constant
almost everywhere. For every $p\notin S$, let
$h_p:\mathbb{Q}_p\to \mathbb{C}$ be a function that (a) is locally constant
almost everywhere, (b) satisfies $|h_p(x)|\leq 1$ for all $x$. 
Let $B(t)\in \mathbb{Z}\lbrack t\rbrack$ be a 
non-zero polynomial. Assume that $h_p(x)=1$ whenever $v_p(B(x))<2$.

Let $\alpha:\mathbb{Z}\to \mathbb{C}$ be a function such that
\begin{equation}\label{eq:gambo}\av_{a + m\mathbb{Z}}\; \alpha(n) = 0\end{equation}
for every arithmetic progression $a + m \mathbb{Z}$. Assume moreover that 
$|\alpha(n)|\leq 1$
for all $x$.
Let
\begin{equation}\label{eq:charolast}
W(n) = \prod_{v\in S} g_v(n) \cdot \prod_{p\notin S}
h_p(n) \cdot \alpha(n).\end{equation}
 Then
\begin{equation}\label{eq:arar}
\av_{a + m \mathbb{Z}} W(n) = 0
\end{equation}
for every arithmetic progression $a + m \mathbb{Z}$.
\end{prop}
\begin{proof}
Let $M$ be a large integer. Let $a+m\mathbb{Z}$ be an arithmetic progression. We will show that
\[\av_{a + m \mathbb{Z}} \prod_{v\in S} g_v(n) \cdot \mathop{\prod_{p\notin S}}_{p<M}
h_p(n) \cdot \alpha(n)
= 0.\]
Equation (\ref{eq:arar}) will then follow as $M\to \infty$, by the same
argument as in Prop.\ \ref{prop:cook}.

We can define $S' = \{S \cup \{p\notin S: p<M\}\}$ and $g_p = h_p$ for $p\in
S'\setminus S$. We must prove that
\[\av_{a + m \mathbb{Z}} \prod_{v\in S'} g_v(n) \cdot \alpha(n) = 0.\]

We now proceed as in the proof of Prop.\ \ref{prop:cake}. For each $p\in S'$,
let $S'_p\subset \mathbb{Q}_p$ be the finite set of points of $\mathbb{Q}_p$
at which $g_p(x)$ is not locally constant. We cover $S'_p$ by open balls whose
union $X_p$ has measure $\leq \frac{\epsilon}{|S'|}$. We let $R_p =
\mathbb{Q}_p\setminus X_p$. (Here $R_p$ stands for ``regular set''; it is
a set on which $g_p(x)$ is everywhere locally constant.) We define
 \[\rho(n) = \begin{cases} 1 &\text{if $\iota_p(n)\in R_p$
for every $p\in S'$,}\\ 0&\text{otherwise.}\end{cases}\]
Then, as in the proof of Prop.\ \ref{prop:cake},
\[\av_{a + m \mathbb{Z}} \prod_{v\in S'} g_v(n) \cdot \alpha(n) = \av_{a +
  m\mathbb{Z}} \rho(n) \prod_{v\in S'} g_v(n) \cdot \alpha(n) 
+ O(\epsilon).\]

Now $R_p$ is compact, and thus can be covered by finitely many balls $U_{p,j}$
on which $g_{p}$ is constant. (Since $g_{\infty}(x) = c_{\infty}$ for all $x$
larger
than a constant, we do not need to worry about what happens at the infinite place $v=\infty$.)
It will be enough to show that
\[\lim_{N\to \infty} \left(\frac{1}{N/m} \mathop{\mathop{\sum_{n\leq N}}_{n\in
      a + m \mathbb{Z}}}_{\forall p\in S': \iota_p(n)\in U_{p,j_p}} 
\prod_{v\in S'} g_v(n) \alpha(n) \right) = 0\]
for every choice of $\vec{j} = \{j_p\}_{p\in S'}$. In other words, we are
partitioning the progression $a + m\mathbb{Z}$ into sets of the form
$Z_{\vec{j}} \cap (a + m \mathbb{Z})$, where 
\[Z_{\vec{j}} =  \{n\in \mathbb{Z}: \forall p\in S'\;\;\; \iota_p(n)\in
U_{p,j_p}\}\]
for some choice of indices $\vec{j}$, and we would like to see that the average of
$g_p(n) \cdot \alpha(n)$ on each such set is $0$.
Now, $g_p(x)$ is
constant on $U_{p,j_p}$, and so we must just show that
\begin{equation}\label{eq:astorio}\lim_{N\to \infty} \left(\frac{1}{N/m} \mathop{\sum_{n\leq N}}_{n\in
      Z_{\vec{j}} \cap (a + m \mathbb{Z})} \alpha(n) \right) = 0\end{equation}
for every possible choice of indices $\vec{j}$.
As we saw at the end of the proof of Prop.\ \ref{prop:cake}, the set
$Z_{\vec{j}}$ is an arithmetic progression, and hence $(a + m
\mathbb{Z}) \cap Z_{\vec{j}}$ is either an arithmetic progression or empty.
If it is empty, (\ref{eq:astorio}) is trivially true; if it is an arithmetic
progression,
(\ref{eq:astorio}) is true by the assumption (\ref{eq:gambo}) that $\alpha(n)$
averages to $0$ on every arithmetic progression.
\end{proof}

Given a function $f:\mathbb{Z}^2\to \mathbb{C}$, a lattice coset 
$a + L\subset \mathbb{Z}^2$ and a sector $S\subset \mathbb{R}^2$ (see \S
\ref{subs:seclat}), we define
\begin{equation}\label{eq:gogon}
\av_{S\cap (a + L), \text{coprime}} f = \lim_{N\to \infty}
 \frac{\sum_{(x,y)\in S\cap (a + L)\cap \lbrack -N, N\rbrack^2:\; \gcd(x,y)=1} f(x,y)}{|\{(x,y)\in S\cap (a + L)\cap \lbrack -N, N\rbrack^2: \gcd(x,y)=1\}|}
 .\end{equation}
This is just like definition (\ref{eq:pabellon}), 
with the condition $\gcd(x,y)=1$ added.

\begin{prop}[$\mathscr{A}_2(B)$]\label{prop:castrol}
Let 
$\mathscr{S}$ be a finite set of places of $\mathbb{Q}$. For every $v\in 
\mathscr{S}$, let
$g_v:\mathbb{Q}_v\times \mathbb{Q}_v \to \mathbb{C}$ be a bounded function
that is locally constant outside a finite set of lines through the origin.
For every $p\notin \mathscr{S}$, let
$h_p:\mathbb{Q}_p\times \mathbb{Q}_p \to \mathbb{C}$ be a function that
(a) is locally constant outside a finite set of lines through the origin,
(b) satisfies $|h_p(x,y)|\leq 1$ for all $x,y\in \mathbb{Q}_p$.

Let $B\in \mathbb{Z} \lbrack x,y\rbrack$ be a non-zero homogeneous 
polynomial. Assume that $h_p(x,y)=1$ whenever $v_p(B(x,y))<2$. Let
$\alpha:\mathbb{Z}\times \mathbb{Z}\to \mathbb{C}$ be a function
such that
\begin{equation}\label{eq:thoto}
\av_{S\cap (a + L), \text{coprime}}\; \alpha(x,y)=0\end{equation}
for every sector $S$ and every lattice coset $a + L$ such that
$\{(x,y)\in L : \gcd(x,y)=1\}$ is non-empty. Assume  moreover that 
$|\alpha(x,y)|\leq 1$
for all $x$.
Let
\[
W(x,y) = \prod_{v\in \mathscr{S}} g_v(x,y) \cdot \prod_{p\notin \mathscr{S}} h_p(x,y) \cdot \alpha(x,y).\]

Then
\begin{equation}\label{eq:thata}
\av_{S \cap (a + L), \text{coprime}} W(x,y) = 0\end{equation}
for every sector $S$ and every lattice coset $L$.
\end{prop}
We could remove the word ``coprime'' from both (\ref{eq:thoto}) and
(\ref{eq:thata}); the statement and the proof would still be correct.
However, we need the word ``coprime'' in (\ref{eq:thata}) for our applications.
After we are done with the proof we are about to go through, we will see
how to remove the word ``coprime'' from (\ref{eq:thoto}) without removing
it from (\ref{eq:thata}).
\begin{proof}
We will proceed just as in the proof of Prop.\ \ref{prop:paneton}, with
the only difference that we shall borrow from Prop.\ \ref{prop:torte} rather 
than from Prop.\ \ref{prop:cake}.

Let $M$ be a large integer. It will be enough to show that
\[
\av_{S \cap (a + L), \text{coprime}} \prod_{v\in \mathscr{S}} g_v(x,y) \cdot
\prod_{p\notin \mathscr{S}} h_p(x,y) \cdot \alpha(x,y) = 0.\]
Equation (\ref{eq:thata}) will then follow as $M\to \infty$, by the same
argument as in Prop.\ \ref{prop:cook} or Prop.\ \ref{prop:cookie}.

We can define $\mathscr{S}' = \{\mathscr{S} \cup \{p\notin \mathscr{S}: p<M\}\}$ and $g_p = h_p$ for $p\in
\mathscr{S}'\setminus \mathscr{S}$. We must prove that
\begin{equation}\label{eq:zorar}
\av_{S \cap (a + L), \text{coprime}} \prod_{v\in \mathscr{S}'} g_v(x,y) 
\cdot \alpha(x,y) = 0.\end{equation}

Each function $g_v$ is locally constant outside a finite set of lines
in $\mathbb{Q}_v \times \mathbb{Q}_v$. We cover these lines by a set $X_v$
defined as the union of sufficiently narrow sectors (that is, angles)
around the lines. 
The number of integers $-N\leq x,y\leq N$ in an angle $\epsilon$ in
$\mathbb{R}^2 = \mathbb{Q}_{\infty}^2$
is $O(\epsilon N^2)$. Hence, if we take care
to make each sector of $X_{\infty}$ of width $\leq \epsilon$,
\[|\{-N\leq x,y\leq N: (x,y)\in X_{\infty}\}| =
O(\epsilon N^2).\]
Let us now consider $v$ finite. Then $X_v$ consists of a union of sets of the
form $\{(x,y)\in \mathbb{Q}_v\times \mathbb{Q}_v:x/y\in B_p\}$, where $B_p$ is a ball in $\mathbb{Q}_p$ of
area $\leq \epsilon$. Now, if the ball $B_p$ has area $1/p^k$, the
set $\{(x,y)\in \mathbb{Q}_v \times \mathbb{Q}_v:x/y\in B_p\}$ is a lattice coset
$a + L$, where $L$ has index $p^k$. The number of
$-N\leq x,y\leq N$ in a lattice coset $a + L$ where $L$ has index $p^k$
is $(1/p^k + o(1)) (2 N)^2$. Hence
\[|\{-N\leq x,y\leq N: (x,y)\in X_v\}| = O(\epsilon N^2).\]
Thus, the total contribution of all sets $X_v$ to (\ref{eq:zorar})
is $O(\epsilon)$ (where the constant depends on $S$, $\mathscr{S}'$ and $L$).

Much as in Prop.\ \ref{prop:paneton}, we define $R_v = \mathbb{Q}_v\setminus X_v$,
\[\rho(x,y) = \begin{cases} 1 &\text{if $\iota_v(x,y)\in R_v$ for
every $v\in \mathscr{S}'$},\\ 0 &\text{otherwise.} \end{cases}\]
We have shown that
\[\av_{S \cap (a + L), \text{coprime}} \prod_{v\in \mathscr{S}'} g_v(x,y)  \alpha(x,y) = 
\av_{S \cap (a + L), \text{coprime}} \rho(x,y) \prod_{v\in \mathscr{S}'} g_v(x,y) 
 \alpha(x,y)+O(\epsilon).\]

For each $v\in \mathscr{S}'$, the function $g_v$ is locally constant on 
$R_v$.
The set $R_{\infty} \in \mathbb{R}^2$ is the union of the closure
of finitely many disjoint sectors $\{U_{\infty,j}\}$; 
since sectors are connected,
$g_{\infty}$ is constant on each sector. For each $p\in \mathscr{S}'$,
the set $R_p$ is compact; thus, it can be covered by finitely many $p$-adic
balls $U_{p,j}$ in $\mathbb{Q}_p \cap \mathbb{Q}_p$ on each of which $g_p$ is constant. 
We must now show that,
for every choice of indices $\vec{j}$,
the average of 
\[\prod_{v\in \mathscr{S}'} g_v(x,y) 
\cdot \alpha(x,y)\]
on each set of the form
\begin{equation}\label{eq:ranan}\{(x,y)\in S\cap S_j \cap (a + L) \cap (a_j+L_{\vec{j}})\cap \lbrack -N, N\rbrack^2: \gcd(x,y)=1\}\end{equation}
is $0$. (Here $S_j$ is a sector $S_j=U_{\infty,j}$ as above, whereas $L_{\vec{j}}$ is
the intersection of preimages $\iota_p^{-1}(U_{p,j_p})$ of $p$-adic balls $U_{p,j_p}$ as above
(one per prime $p\in \mathscr{S}'$) under the injection
$\iota_p:\mathbb{Z}\times \mathbb{Z} \to \mathbb{Q}_p \times \mathbb{Q}_p$;
since every preimage $\iota_p^{-1}(U_{p,j_p})$ is a lattice coset, the
intersection $L_{\vec{j}} = \cap_{p\in \mathscr{S}'} \iota_p^{-1}(U_{p,j_p})$ is itself a lattice coset.)
The product $\prod_{v\in \mathscr{S}'} g_v(x,y)$ is constant
on each set (\ref{eq:ranan}), and condition (\ref{eq:thoto}) assures us that
$\alpha(x,y)$  averages to $0$ on every set of the form
(\ref{eq:ranan}). Hence we are done.
\end{proof}

We should now prove that we can remove the word ``coprime'' from
(\ref{eq:thoto}). Let us first prove an easy preparatory lemma.
\begin{lem}\label{lem:hogg}
Let $a + L \subset \mathbb{Z}^2$ be a lattice coset. Assume that $\{(x,y)\in a + L: \gcd(x,y)=1\}$ 
is non-empty. Then, for any sector $S$,
\begin{equation}\label{eq:hotoso}\lim_{N\to \infty} \frac{1}{(2 N)^2}
|\{(x,y)\in S\cap (a + L)\cap \lbrack - N, N\rbrack^2: \gcd(x,y)=1\}|
\end{equation}
exists and is non-zero.
\end{lem}
By the definition of sectors (\ref{subs:seclat}), all sectors are non-empty.
\begin{proof}
We will prove something stronger and more explicit:
\begin{equation}\label{eq:dredon}\begin{aligned}
\lim_{N\to \infty} \frac{1}{(2 N)^2}
&|\{(x,y)\in S\cap (a + L)\cap \lbrack - N, N\rbrack^2: \gcd(x,y)=1\}| \\ &= 
\frac{c_{\infty}}{\lbrack \mathbb{Z}^2 : L\rbrack}
 \prod_{p|
|\lbrack \mathbb{Z}^2:L\rbrack|
} (1 - c_p)
\prod_{p \nmid |\lbrack \mathbb{Z}^2:L\rbrack|}
\left(1 - \frac{1}{p^2}\right),\end{aligned}\end{equation}
where
\begin{equation}\label{eq:comino}\begin{aligned}
c_{\infty} &= \lim_{N\to \infty} \frac{|S \cap \lbrack -N, N\|}{(2 N)^2} ,\\
c_p &= \begin{cases} 0 &\text{if 
$(a + L)\cap p \mathbb{Z}^2 = \emptyset$}\\
\frac{\lbrack \mathbb{Z}^2 : L\rbrack}{\lbrack \mathbb{Z}^2 :
L_p\rbrack} &\text{if $(a + L)\cap p \mathbb{Z}^2$ is a coset of
a lattice $L_p\subset \mathbb{Z}^2$.}\end{cases}
\end{aligned}\end{equation}
Now that we know that we have to prove this, the rest is very easy.

Let $M$ be a large integer. Then the right side of (\ref{eq:dredon})
is
\[
\frac{c_{\infty}}{\lbrack \mathbb{Z}^2 : L\rbrack}
\mathop{\prod_{p|
|\lbrack \mathbb{Z}^2:L\rbrack|
}}_{p\leq M} c_p
\mathop{\prod_{p \nmid |\lbrack \mathbb{Z}^2:L\rbrack|}}_{p\leq M}
\left(1 - \frac{1}{p^2}\right) + O\left(\frac{1}{M}\right),\]
whereas
\[|\{-N\leq x,y\leq N: \text{$\exists p>M$ such that $p|\gcd(x,y)$}\}|
= O\left(\frac{1}{M} + \frac{1}{N} \right) \cdot N^2.\]
and so
\[\begin{aligned}
|\{(x,y)\in &S\cap (a + L)\cap \lbrack - N, N\rbrack^2: \gcd(x,y)=1\}|\\ &=
|\{(x,y)\in S\cap (a + L)\cap \lbrack -N, N\rbrack^2: 
(p\nmid x\; \vee\; p\nmid y)\; \forall p\leq M\}|\\ &+ 
O\left(\frac{1}{M} + \frac{1}{N} \right) \cdot N^2.\end{aligned}
\]
Hence, it will be enough to show that
\begin{equation}\label{eq:toyota}\lim_{N\to \infty} \frac{1}{(2 N)^2}
|\{(x,y)\in S\cap (a + L)\cap \lbrack -N, N\rbrack^2: 
(p\nmid x\; \vee\; p\nmid y)\; \forall p\leq M\}|\end{equation}
equals
\[\frac{c_{\infty}}{\lbrack \mathbb{Z}^2 : L\rbrack}
 \cdot \mathop{\prod_{p|
|\lbrack \mathbb{Z}^2:L\rbrack|
}}_{p\leq M} (1 - c_p)
\mathop{\prod_{p \nmid |\lbrack \mathbb{Z}^2:L\rbrack|}}_{p\leq M}
\left(1 - \frac{1}{p^2}\right).\]

Let 
\[c_m = \begin{cases} 0 &\text{if 
$(a + L)\cap m \mathbb{Z}^2 = \emptyset$}\\
\frac{\lbrack \mathbb{Z}^2 : L\rbrack}{\lbrack \mathbb{Z}^2 :
L_m\rbrack} &\text{if $(a + L)\cap m \mathbb{Z}^2$ is a coset of
a lattice $L_m\subset \mathbb{Z}^2$.}\end{cases}\]
Now
\begin{equation}\label{eq:tgarni}\begin{aligned}
|\{(x,y)\in &S\cap (a +L)\cap \lbrack -N, N\rbrack^2: 
(p\nmid x\; \vee\; p\nmid y)\; \forall p\leq M\}|\\ &=
\mathop{\sum_{m}}_{p|m \Rightarrow p\leq M} \mu(m) \cdot
\{(x,y)\in S\cap (a + L)\cap \lbrack -N, N\rbrack^2: m|x\; \wedge m|y\}\\
&= \mathop{\sum_{m}}_{p|m \Rightarrow p\leq M} \mu(m) \cdot
\frac{c_{\infty}}{\lbrack \mathbb{Z}^2 : L\rbrack}
\cdot (2 N)^2 c_m,
\end{aligned}\end{equation}
where we are using (\ref{eq:crusty}). For $m$ square-free,
\[c_m = \prod_{p|m} c_p.\]
Hence (\ref{eq:toyota}) equals
\[\frac{c_{\infty}}{\lbrack \mathbb{Z}^2 : L\rbrack}
\cdot \prod_{p\leq M} (1 - c_p).\]
Now, given two lattice cosets $a + L$, $a + L'$ with coprime indices
$\lbrack \mathbb{Z}^2:L\rbrack$, $\lbrack \mathbb{Z}^2:L'\rbrack$, 
their intersection is a lattice coset of index
$\lbrack \mathbb{Z}^2:L\rbrack\cdot \lbrack \mathbb{Z}^2:L'\rbrack$.
Hence (by (\ref{eq:comino}))
\[c_p = \frac{1}{p^2}\]
for $p\nmid |\lbrack \mathbb{Z}^2:L\rbrack|$,
and so we are done.
\end{proof}
\begin{Rem}
We do not particularly care about the rate of convergence in
(\ref{eq:hotoso}); otherwise we would have proven Lemma \ref{lem:hogg}
differently. For a discussion on how to do things efficiently (from the 
point of view of rates of convergence), see appendix \ref{sec:ratcon}.
(Actually, proving Lemma \ref{lem:hogg} efficiently is fairly easy; taking care of other similar
issues elsewhere in the paper is a harder problem. See also \cite[\S 3]{He}, 
where this sort of issue is dealt with at length.)
\end{Rem}

\begin{lem}\label{lem:honn}
Let $f:\mathbb{Z}^2\to \mathbb{C}$ be a bounded function.
Let $S\subset \mathbb{R}^2$ be a sector. Suppose that
\begin{equation}\label{eq:jamon}\av_{S\cap (a + L)} f = 0\end{equation}
for every lattice coset $a + L$. Then
\begin{equation}\label{eq:queso}\av_{S\cap (a + L), \text{coprime}} f = 0
\end{equation}
for every lattice coset $a + L$ such that 
$\{(x,y)\in a + L: \gcd(x,y)=1\}$ is non-empty.
\end{lem}
\begin{proof}[Sketch of proof]
By Lemma \ref{lem:hogg}, the denominator implicit in the average
in (\ref{eq:queso}) is $\gg N^2$. Hence, it is enough to show that
\[\lim_{N\to \infty} \frac{1}{N^2}
|\{(x,y)\in S\cap (a + L)\cap \lbrack - N, N\rbrack^2: \gcd(x,y)=1\}| = 0.\]
Using (\ref{eq:jamon}), we can prove this just like we proved 
Lemma \ref{lem:hogg}: first we fix a large integer $M$ (which we will let
go to $\infty$ at the end), then we show (exactly as before)
that the difference between
\[|\{(x,y)\in S\cap (a + L)\cap \lbrack - N, N\rbrack^2: \gcd(x,y)=1\}|\]
and
\begin{equation}\label{eq:ojota}
|\{(x,y)\in S\cap (a + L)\cap \lbrack - N, N\rbrack^2: 
(p\nmid x\; \vee\; p\nmid y)\; \forall p\leq M\}|
\}|\end{equation}
is $O(1/M)$, then we express (\ref{eq:ojota}) as a sum of the form
$\sum_m \mu(m) \dotsc$ as in (\ref{eq:tgarni}) -- and, lastly, we use (\ref{eq:jamon})
instead of (\ref{eq:crusty}).
\end{proof}
\subsection{Families with multiplicative reduction}\label{subs:edi}
It remains to prove the main theorems. Recall that, given
a function $f:\mathbb{Z}\to \mathbb{C}$, ``$f$ has strong zero
average over the integers'' means that
\[\av_{a + m \mathbb{Z}} f(n) = 0\]
for every arithmetic progression $f(n)$. For $f:\mathbb{Q} \to \mathbb{C}$,
``$f$ has strong zero average over
the rationals'' means that
\[\av_{\mathbb{Q}, S\cap (a + L)} g = 0
\]
for every sector $S\subset \mathbb{R}^2$ and every lattice coset
$a + L \subset \mathbb{Z}^2$. Lastly, for $f:\mathbb{Z}^2 \to \mathbb{C}$,
``$f$ has strong zero average over $\mathbb{Z}^2$'' means that
\[\av_{S\cap (a + L)} g = 0
\]
for every sector $S\subset \mathbb{R}^2$ and every lattice coset
$a + L \subset \mathbb{Z}^2$. See also (\ref{eq:gogon}).

Recall as well that, when we write
{\bf Theorem 0.0} ($\mathfrak{X}(P)$, $\mathfrak{Y}(Q)$),
we mean a theorem conditional on hypotheses $\mathfrak{X}$ and $\mathfrak{Y}$
in so far as they concern the objects $P$ and $Q$, respectively. 
The polynomials $B_{\mathscr{E}}$ and $M_{\mathscr{E}}$ were defined in
(\ref{eq:M}); they are the polynomials describing the places at which
a curve $\mathscr{E}$ over $\mathbb{Q}(T)$ has quite bad or multiplicative
reduction, respectively.

\begin{main1}[$\mathscr{A}_1(B_\mathscr{E}(t,1))$,
$\mathscr{B}_1(M_\mathscr{E}(t,1))$]
Let $\mathscr{E}$ be a family of elliptic curves over $\mathbb{Q}$. 
Assume that 
$M_\mathscr{E}(t,1)$ is not constant. Then $t\mapsto W(\mathscr{E}(t))$ 
has strong zero average over the integers.
\end{main1}
\begin{proof}
Immediate from Theorem \ref{thm:smet} and Prop.\ \ref{prop:paneton}
with $\alpha(x) = \lambda(M_{\mathscr{E}}(x,1))$. Condition (\ref{eq:gambo})
is furnished by hypothesis $\mathscr{B}_1(M_{\mathscr{E}}(x,1))$.
\end{proof}

\begin{main2}[$\mathscr{A}_2(B_\mathscr{E})$,
$\mathscr{B}_2(M_\mathscr{E})$]
Let $\mathscr{E}$ be a family of elliptic curves over $\mathbb{Q}$. 
Assume that 
$M_\mathscr{E}$ is not constant\footnote{In other words,
assume that $\mathscr{E}$ has at least one point of multiplicative
reduction over $\mathbb{Q}(T)$.}. 
Then $t\mapsto W(\mathscr{E}(t))$ has
strong zero average over the rationals.
\end{main2}
\begin{proof}
Immediate from Theorem \ref{thm:smet} and Prop.\ \ref{prop:castrol}
with $\alpha(x,y) = \lambda(M_{\mathscr{E}}(x,y))$.
(Use Lemma \ref{lem:honn} to prove condition (\ref{eq:thoto})
using hypothesis $\mathscr{B}_2(M_{\mathscr{E}})$.)
\end{proof}

One can go in the other direction: if the average root number of a family
is zero, then we can prove hypothesis $\mathscr{B}_j$ for the associated
polynomial $M_{\mathscr{E}}$. 

\begin{prop}[$\mathscr{A}_1(B_\mathscr{E}(t,1))$]
Let $\mathscr{E}$ be a family of elliptic curves over $\mathbb{Q}$. 
Assume that 
$M_\mathscr{E}(t,1)$ is not constant. 

Suppose that $t\mapsto W(\mathscr{E}(t))$ 
has strong zero average over the integers. Then
hypothesis $\mathscr{B}_1(M_{\mathscr{E}}(x,1))$ holds, i.e.,
$\lambda(M_{\mathscr{E}}(x,1))$ has strong zero average over the integers.
\end{prop}
\begin{proof}
Equation (\ref{eq:exerier}) in Theorem \ref{thm:smet} states that
\begin{equation}\label{eq:rara}W(\mathscr{E}(x/y)) = g(x,y) \cdot h(x,y)\cdot 
\lambda(M_{\mathscr{E}}(x,y)).\end{equation}
Since the left side is never zero, the expressions on the right side are
non-zero for all but finitely many rationals $x/y$. (Equation
(\ref{eq:exerier}) is valid for all but finitely many rationals.)
This implies that
\begin{equation}\label{eq:roro}\lambda(M_{\mathscr{E}}(x,y)) = 
(g(x,y))^{-1} \cdot (h(x,y))^{-1}\cdot W(\mathscr{E}(x/y)) .\end{equation}
Now apply Prop.\ \ref{prop:paneton}
with $\alpha(x) = W(\mathscr{E}(x))$. 
\end{proof}

\begin{prop}[$\mathscr{A}_2(B_\mathscr{E})$]
Let $\mathscr{E}$ be a family of elliptic curves over $\mathbb{Q}$. 
Assume that 
$M_\mathscr{E}$ is not constant.
Then $\mathscr{B}_2(M_{\mathscr{E}}(x,y))$ holds, i.e.,
$\lambda(M_{\mathscr{E}}(x,y))$ has strong zero average over $\mathbb{Z}^2$.
\end{prop}
\begin{proof}
Again, (\ref{eq:rara}) implies (\ref{eq:roro}). Now apply
\[\lambda(M_{\mathscr{E}}(x,y)) = 
(g(x,y))^{-1} \cdot (h(x,y))^{-1}\cdot W(\mathscr{E}(x/y)).\]
Now apply Prop.\ \ref{prop:castrol}
with $\alpha(x,y) = W(\mathscr{E}(x,y))$.
We obtain that
\begin{equation}\label{eq:vick}
\av_{S\cap L, \text{coprime}} \lambda(M_{\mathscr{E}}(x,y)) = 0.\end{equation}
Since $\lambda$ is completely multiplicative and
$M_{\mathscr{E}}$ is homogeneous,
\[\lambda(M_{\mathscr{E}}(m x,m y)) = \lambda(m^d) \lambda(M_{\mathscr{E}}(x,y))\]
for all $m$, $x$, $y$, where $d = \deg(M_{\mathscr{E}})$. Summing over
all $m\leq M$, $M$ large, we deduce from (\ref{eq:vick}) that
\[\begin{aligned}
\mathop{\sum_{(x,y)\in S\cap L \cap \lbrack - N, N\rbrack^2}}_{\gcd(x,y)\leq
  M}
\lambda(M_{\mathscr{E}}(x,y))
&= 
\sum_{1\leq m\leq M} \mu(m) \lambda(m^d)  \cdot o(N^2/m^2)\\
&= o(N^2).\end{aligned}\]
The number of pairs $(x,y)\in S\cap L \cap \lbrack - N,N\rbrack^2$ such
that $\gcd(x,y)>M$ is $O(N^2/M) + O(N)$. We let $N\to \infty$ and obtain
\[\av_{S\cap L} \lambda(M_{\mathscr{E}}(x,y)) = O(1/M) + o(1).\]
Since $M$ can be set to be arbitrarily large, we are done.
\end{proof}

In general, hypothesis $\mathscr{B}_2$ is very difficult to prove.
The state of knowledge on the subject was static for over a century.
It was then proven for polynomials of degree $3$ by the author
(\cite{He3}, \cite{Heirr}; the latter paper builds on work by
Friedlander and Iwaniec \cite{FI} as well as Heath-Brown
\cite{HB}). Hypothesis
$\mathscr{B}_2$ has recently been proven
 for products of four linear factors by Green and Tao \cite{GT}
(using methods related to their proof of the existence of arbitrarily
long progressions in the primes). It is possible that Green and Tao 
(or their associates) shall soon prove $\mathscr{B}_2$ for products of more
than four linear polynomials in $\mathbb{Z}\lbrack x,y\rbrack$. Further
progress in the near future seems unlikely. 

Hypothesis $\mathscr{B}_1$ is even harder. It has been known for degree 1
for more than a hundred years, but it is extremely improbable that it
will proven for any (square-free) polynomials of degree $\geq 2$ in the near 
future: such a task would be of the same difficulty as proving the twin
prime number conjecture.

To summarise: 
(a) the main theorems
assert the truth of general results on root numbers under the assumption of standard 
conjectures in analytic number theory (hypotheses $\mathscr{A}_j$,
$\mathscr{B}_j$); (b) thus, these general results are
unconditionally true in the cases in which these conjectures are known
to be true; (c) if the results on numbers were proven unconditionally for
some more cases, then one of the harder conjectures ($\mathscr{B}_1$ or
$\mathscr{B}_2$) would follow in a case in which it is not yet known.
(In (c), one of the easier conjectures ($\mathscr{A}_j$) is assumed;
recall that the easier conjectures are known for all polynomials of
degree $\leq 3$ (if $j=1$) or degree $\leq 6$ (if $j=2$).)

It thus seems that the main theorems close the matter of averages of root
numbers on one-parameter families for the while being.

\subsection{Some unconditional examples}
It is a simple matter to construct
a family $\mathscr{E}$ with 
$M_{\mathscr{E}}$ equal to a given square-free homogeneous polynomial in
$\mathscr{O}_{K,V}\lbrack x,y\rbrack$, up to multiplication by a scalar.
For example, if we want $M_{\mathscr{E}} = x^3 + 2 y^3$, we may choose
\[c_4 = 1 - 1728 (t^3+2),\;\;\;
c_6 = (1 - 1728 (t^3+2))^2 .\] 
We then have
\[\Delta = \frac{c_4^3 - c_6^2}{1728} = (t^3 + 2) \cdot (1 - 1728 (t^3 +
2))^3,\]
and so the only place of multiplicative reduction over $\mathbb{Q}(T)$
is the one corresponding to $t^3 + 2\in \mathbb{Q}\lbrack T\rbrack$. In 
other words, $M_{\mathscr{E}} = x^3 + 2 y^3$.

In general, we can construct an elliptic curve $\mathscr{E}$ over $K(T)$ 
with given factors in the denominator of $j$. The factors in the
denominator of $j$ correspond to the places of $K(T)$ where the reduction
is either additive and potentially multiplicative or multiplicative.
We can then apply quadratic
twists to make places of additive, potentially multiplicative
 reduction into places
of multiplicative reduction, and vice versa, until we have multiplicative
reduction exactly at the places we wish to have it, and at no other place.
It is thus that we can find families $\mathscr{E}$ with any given
(square-free, homogeneous) $M_{\mathscr{E}}$.

Families $\mathscr{E}$ with $M_{\mathscr{E}}$ given
and $c_4$, $c_6$ coprime
are slightly harder to find. Here are some examples:
\begin{equation}\label{eq:jojojo}\begin{array}{lll}
c_4 = 2 + 4 t + t^2,  &c_6 = 1 + 9 t + 6 t^2 + t^3,
&M_{\mathscr{E}} = (2 x + 7 y) ( x^2 + 4 x y + y^2),\\
c_4 = 2 - 4 t + t^2, &c_6 = 3 + 9 t - 6 t^2 + t^3,
&M_{\mathscr{E}} = 10 x^3 - 63 x^2 y + 102 x y^2 + y^3,\\
c_4=4, &c_6 = 11 + t, &M_{\mathscr{E}} = y (x + 3 y) (x + 19 y),\\
c_4=3, &c_6 = 2 + 7 t, &M_{\mathscr{E}} = y (49 x^2 + 28 x y - 23 y^2),\\
c_4 = 1+t, &c_6=-1 + 3t, &M_{\mathscr{E}} = x y (x -3 y).\end{array}
\end{equation}
Since $\deg(M_{\mathscr{E}})\leq 3$ and $\deg(B_{\mathscr{E}})\leq 6$
in all of these cases, Main Theorem 2 is unconditional for all of
these families.

\appendix
\section{Rates of convergence}\label{sec:ratcon}

We have computed averages over $\mathbb{Z}$ and $\mathbb{Q}$, but
we have not determined their rates of convergence. How do we show that
-- for families $\mathscr{E}$ for which we can determine the average
$\av_{\mathbb{Z}} W(\mathscr{E}(t)$ and
$\av_{\mathbb{Q}} W(\mathscr{E}(t)$ 
unconditionally -- the expressions
\[\begin{aligned}
\frac{1}{N} \sum_{1\leq n\leq N} W(\mathscr{E}(n)) &= 
\av_{\mathbb{Z}} W(\mathscr{E}(t)) + o(1),\\
\frac{\sum_{(x,y)\in \lbrack - N,N\rbrack^2: \gcd(x,y)=1} 
W(\mathscr{E}(x/y))}{
|\{(x,y)\in \lbrack - N, N\rbrack^2: \gcd(x,y)=1\}|} &=
\av_{\mathbb{Q}} W(\mathscr{E}(t)) + o(1)
\end{aligned}\]
have small error terms $o(1)$?

There are limitations imposed by the inputs: there are error
terms implicit in hypotheses
$\mathscr{A}_1$, $\mathscr{A}_2$, $\mathscr{B}_1$ and $\mathscr{B}_2$,
and the cases of the hypotheses that are known are not always known
with the best error terms they might have. For example, when I showed 
\cite[Thm.\ 3.3]{He3} that $\av_{\mathbb{Z}^2} \lambda(x y (x + y)) = 0$,
I showed, in fact, that
\begin{equation}\label{eq:toroto}
\frac{1}{|\lbrack - N, N\rbrack|}
\sum_{-N\leq x,y\leq N} \lambda(x y (x+y)) = 
O\left(\frac{\log \log N}{\log N}\right) ,
\end{equation}
where the implied constant is absolute.
It is quite probable that (\ref{eq:toroto}) is still true if
$O\left(\frac{\log \log N}{\log N}\right)$ is replaced by
$O(N^{-\beta})$, $\beta>0$ small, but we are far from proving that.

Putting the issue of the inputs to one side, we can see -- if we
go through the proofs in the present paper keeping track of the error
terms - that there is one clearly non-trivial issue. The issue
is the passage to the limit at the beginning of
 the proofs 
of Propositions \ref{prop:cook} and \ref{prop:paneton} (and their
two-dimensional analogues,
Propositions \ref{prop:cookie} and \ref{prop:castrol}). 

This issue is not particular to the study of root numbers. It arises,
for example, in the classical problem of determining the number of integers
$n\leq N$ such that $f(n)$ is square-free (or free of factors that
are $k$th powers, in general), where $f$ is a polynomial. 
Say we know that
\begin{equation}\label{eq:gustu}
|\{n\leq N: \text{$p^2|f(n)$ for some $p>n^{0.1}$}\}| = O(N^{1-\beta}),\end{equation}
where $\beta>0$. We would like to conclude that
\begin{equation}\label{eq:gasta}
|\{n\leq N: \text{$f(n)$ is square-free}\}| = c_f \cdot N + O(N^{1-\beta})
+O(N^{1-\beta'}),\end{equation}
say, where $c_f$ is a constant and $\beta'>0$. 
It turns out that one can in fact conclude this 
(with $\beta' = 1/3 - \epsilon$, $\epsilon>0$ arbitrarily small), but this is not very easy.
If one were to do things in the obvious way, one would end up having an error
term of $O(N/(\log N))$ instead of $O(N^{1-\beta'})$.

The implication (\ref{eq:gustu}) $\Rightarrow$ (\ref{eq:gasta}) was proven
in \cite[Prop.\ 3.4]{He}. The issue is solved in full generality
in \cite[\S 3.2--3.4]{He}. In particular, \cite[Prop.\ 3.8 and 3.9]{He}
are explicit analogues of Propositions \ref{prop:cook}
and \ref{prop:cookie} in the present paper; they give explicit (and rather
good) error terms. The main idea (namely, the optimised alternative
to the passage to the limit at the beginning of the proof of Prop.\
\ref{prop:cook}) is contained in a rather cryptic-looking general result
\cite[Prop.\ 3.2]{He}. (The setup and proof of
\cite[Prop.\ 3.2]{He} were called
a {\em riddle}; the term was meant to suggest a coarse sieve.)

In general, \cite[\S 3]{He} contains a great deal of what we did with
much less work in \S \ref{subs:hoat}, \S \ref{subs:pajaro} and
\S \ref{subs:proca}. The main
difference is that, in \cite[\S 3]{He}, a great deal of effort was put
into making all error terms explicit and good. (An unintended but
predictable effect was the fact that \cite[\S 3]{He} is a little hard to read.)
If the results
in \S \ref{subs:hoat}, \S \ref{subs:pajaro} and \S \ref{subs:proca}
are replaced by those
in \cite[\S 3]{He}, one obtains analogues of Propositions \ref{prop:ant1}
and \ref{prop:ant2} and of the main theorems with
explicit error terms. As the inputs -- namely, hypotheses $\mathscr{A}_j$
and $\mathscr{B}_j$ -- do not always give good error terms themselves,
the gains over what one would get by proceeding naively are not always
large.

\section{A non-constant family with non-zero average root number over $\mathbb{Q}$}\label{sec:longexample}

A {\em non-constant} family $\mathscr{E}$ is one where the
fibres $\mathscr{E}(t)$ are not all isomorphic to each other over
$\mathbb{C}$. Families of quadratic twists, for example, are not non-constant.
There are known examples of (a) families of quadratic twists $\mathscr{E}$
with $\av_{\mathbb{Q}} W(\mathscr{E}(t)) \ne 0$ (\cite{Riz1}) and (b)
non-constant families $\mathscr{E}$ with
$\av_{\mathbb{Z}} W(\mathscr{E}(t))=0$ (\cite{Riz2}). 
In this appendix, we will see
an example of a non-constant family $\mathscr{E}$ with 
$\av_{\mathbb{Q}} W(\mathscr{E}(t)) \ne 0$.

We recall that one of our main results in the present paper was that
$\av_{\mathbb{Q}} W(\mathscr{E}(t))=0$ for all families $\mathscr{E}$ with
places of multiplicative reduction, conditionally on two standard
conjectures in analytic number theory. We also gave an expression
(Prop. \ref{prop:ant2} (\ref{eq:brah2})) for $\av_{\mathbb{Q}}
W(\mathscr{E}(t))$ when $\mathscr{E}$ has no places of multiplicative reduction.
In the light of (\ref{eq:brah2}), it is unsurprising that there
are non-constant families $\mathscr{E}$ with no places of multiplicative reduction
and $\av_{\mathbb{Q}} W(\mathscr{E}(t)) \ne 0$. What follows serves as a
completely explicit rendering of most of the arguments of the present paper
-- or, at least, of those in \S \ref{sec:globroo}; we will also need
 the results in \S \ref{sec:avcor}.

The arguments in \S \ref{sec:locroot}
will not appear here, as we engaged in them
precisely in order to avoid explicit case-work. 
Since the family we will treat has no multiplicative reduction,
the steps corresponding to \S \ref{eq:ostorma} will be exceedingly
simple, though the key element (namely, quadratic reciprocity) will still
be present (vd.\ equation (\ref{eq:otorongo})).

It will be satisfying to see the infinite product
of integrals in (\ref{eq:brah2}) take a very concrete form, namely,
equation (\ref{eq:costro}).
In general, we have shown that the $p$-adic integrals in (\ref{eq:brah2})
are rational, and that the double integral in (\ref{eq:brah2})
is algebraic; this is part of the statement of Prop.\ \ref{prop:ant2}.
The same method that we used to prove this can be used to provide an
algorithm (long and tedious) to give (complicated) explicit expressions
in terms of the definition of $\mathscr{E}$ for all of the integrals in 
(\ref{eq:brah2}).

The example we are about to see is actually fairly representative, though
it is a little simpler than most other curves without multiplicative reduction
would be. (This is why this example was chosen.) 
In our example, the expressions we will get for
the integrals $\prod_{p\notin S} \int_{O_p} h_p(x,y)\; dx dy$ depend only
on the number of solutions in $\mathbb{Z}/p\mathbb{Z}$
to $19+11 t^2 + 19 t^4 \equiv 0 \mo p$ -- the integrals
 do not depend on $p$ otherwise; an exception is made for $p=2,3,7,19$,
which will have to be treated separately.
In general, for a family $\mathscr{E}$ without multiplicative reduction,
the integrals $\prod_{p\notin S} \int_{O_p} h_p(x,y)\; dx dy$ 
(for $p$ outside a finite set of primes that have to be treated separately)
depend only
on the number of solutions in $\mathbb{Z}/p\mathbb{Z}$ to 
each of the equations $f_i(x)\equiv 0 \mo p$, where
$f_1,f_2,\dotsc, f_k \in \mathbb{Z}\lbrack x \rbrack$ are fixed polynomials
depending only on the definition of $\mathscr{E}$.

\subsection{The family $\mathscr{E}$. Preparatory work}

Let us, then, construct a non-constant family $\mathscr{E}/\mathbb{Q}\lbrack T\rbrack$ such that 
$W(\mathscr{E}(t))$ does not average to $0$ over the rationals. 
Set $f_1 = -5 - 2 T^2$, $f_2 = 2 + 5 T^2$. Let $\mathscr{E}$ be the elliptic
curve over $\mathbb{Q}(T)$ given by the parameters
\[c_4 = f_1 f_2 (f_1^3 - f_2^3)^2,\;\;\;
c_6 = \frac{1}{2} (f_1^3 + f_2^3) (f_1^3 - f_2^3)^3,\]
\[\Delta = \frac{c_4^3 - c_6^2}{1728} = - 2^{-8} 3^{-3} (f_1^3 - f_2^3)^8 .\]

Here note that
\begin{equation}\label{eq:oneil}\begin{aligned}
f_1^3(t) - f_2^3(t) &= 
 (f_1^2(t) + f_1(t) f_2(t) + f_2^2(t)) (f_1(t) - f_2(t)) \\
&= -7 (19 + 11 t^2 + 19 t^4) (1 + t^2) .\end{aligned}\end{equation}
and
\begin{equation}\begin{aligned}
f_1^3(t) + f_2^3(t) &= 
 (f_1^2(t) - f_1(t) f_2(t) + f_2^2(t)) (f_1(t) + f_2(t)) \\
&= 3^2 (13 t^4 + 23 t^2 + 13) (t - 1) (t+1).
\end{aligned}\end{equation}

Let us determine the sets and polynomials that were defined
in \S \ref{sss:nota}. The set of interesting places $\mathscr{P}$ consists
of all places $w$ of $K(T)$ such that $w(c_4)\ne 0$, $w(c_6)\ne 0$
or $w(\Delta)\ne 0$. In our example, 
\[\mathscr{P} = \{w_{f_1}, w_{f_2}, w_{t-1}, w_{t+1}, w_{t^2+1},
w_{19 t^4 + 11 t^2 + 19}, w_{13 t^4 + 23 t^2 + 13}, \deg(\den)-\deg(\num)\},\]
where we denote by $w_f$ the place of $K(T)$ associated to an
irreducible polynomial in $K\lbrack T\rbrack$. The homogeneous
polynomials in $\mathbb{Z}\lbrack x,y\rbrack$
 associated to these places (\ref{subs:fpp}) are
\begin{equation}\label{eq:baboo}\begin{aligned}P_1 &= (-2 (x/y)^2 -5)\cdot y^2 = -2 x^2 - 5 y^2,\\
P_2 &= 5 x^2 + 2 y^2,\;\;\;\;
P_3 = x - y,\;\;\;\;
P_4 = x + y,\;\;\;\;
P_5 = x^2 + y^2,\\
P_6 &= 19 x^4 + 11 x^2 y^2 + 19 y^4,
\;\;\;\;
P_7 = 13 x^4 + 23 x^2 y^2 + 13
y^4,\;\;\;\; P_8 &= y.\end{aligned}\end{equation}
The resultants $\Res(P_i,P_j)$, $1\leq i<j\leq 8$, have the following
prime factors:
\begin{equation}\label{eq:gosto}2,\; 3,\; 5,\; 7,\; 13,\; 19.\end{equation}
(For example, $\Res(P_6,P_7) = 2^4 \cdot 3^4 \cdot 7^8$ and
$\Res(P_1,P_8) = -5$.) If we consider only $P_i$, $P_j$ with 
$1\leq i< j\leq 7$ (that is, if we leave out the possibility $P_j=P_8$), the
resultants $\Res(P_i,P_j)$ have only the following prime factors:
\begin{equation}\label{eq:shortli}2,\; 3,\; 7.\end{equation}

We write
\begin{equation}\label{eq:ferrat}\begin{aligned}
c_4(x/y) &= P_1(x,y) P_2(x,y) (P_5(x,y) \cdot P_6(x,y))^3 P_8(x,y)^{-16},\\
c_6(x/y) &= \frac{1}{2} (P_3(x,y) P_4(x,y) P_7(x,y)) \cdot 
(P_5(x,y) P_6(x,y))^3 P_8(x,y)^{-24},\\
\Delta(x/y) &= -2^{-8} 3^{-3} (P_5(x,y) P_6(x,y))^8 P_8(x,y)^{-48}
. \end{aligned}\end{equation}
Thus, the prime factors of the ideal defined in (\ref{eq:toolong})
are $2$, $3$, $5$, $7$, $13$, $19$.

It is also worthwhile to note that the prime factors of the discriminant
of every polynomial $P_j$, $1\leq j\leq 8$, are all contained in the set
$\{2,3,7,19\}$.

We will look at the local root number of $W(\mathscr{E}(t))$ at the
places $p= 2,3,5,7,13,19$, and then at all other places. 

\subsection{Local root numbers}\label{subs:lort}
\subsubsection{The root number at $p=2$} Suppose first that $v_2(t)>0$.
Then $v_2(f_1(t))=0$, $v_2(f_2(t))=1$, $v_2(f_1^3(t) - f_2^3(t)) = 0$,
$v_2(f_1^3(t) + f_2^3(t))=0$. Hence $v_2(c_4(t))=1$, $v_2(c_6(t)) = -1$,
$v_2(\Delta(t)) = -8$. Since $v_2(c_6(t))=-1<0$ and $v_2(\Delta)=-8<0$, this
cannot be an integral model. We set $c_4' = 2^4 c_4$, $c_6' = 2^6 c_6$,
$\Delta' = 2^{12} \Delta$. This model has
\[v_2(c_4') = 5,\;\;\;\;\; v_2(c_6') = 5,\;\;\;\;\; v_2(\Delta') = 4 .\]
This new model can be written as $y^2 = x^3 - 27 c_4 - 54 c_6$;
since all of the coefficients of this equation 
are in $\mathbb{Z}_3$, the new model is integral.
(A model is integral if it can be written down in (short or long) Weierstrass
form with integral coefficients.) Since the model given by $c_4$, $c_6$,
$\Delta$ was not integral, the new model is a minimal integral model.

Before we use \cite[Table 1]{Ha}, we must determine 
$2^{-v_2(c_6)} c_6 = 2 c_6 \mo 8$. Now $f_1(t) \equiv -5\equiv 3 \mo 8$ and
$v(f_2(t))>0$.
Hence $f_1^3(t) - f_2^3(t) \equiv 3 \mo 8$ and $f_1^3(t) + f_2^3(t) 
\equiv 3 \mo 8$. Thus $2 c_6 \equiv 1 \mo 8$.
The second line marked III in \cite[Table 1]{Ha} now gives us 
$W_2(\mathscr{E}(t))=-1$.

Suppose now that $v_2(t)=0$. 
Then $v_2(f_1(t))=0$ (and $f_1(t)\equiv 1 \mo 8$), 
$v_2(f_2(t))=0$ (and $f_2(t)\equiv -1 \mo 8$), $v_2(f_1^3(t) - f_2^3(t)) = 1$,
$v_2(f_1^3(t) + f_2^3(t))\geq 3$. Hence $v_2(c_4(t))=2$, $v_2(c_6(t)) \geq 5$,
$v_2(\Delta(t)) = 0$. Suppose that this model were integral.
By the transformations in \cite[\S 3.1]{Si} (or
any other standard text), $c_4(t) = b_2^2 - 24 b_4$ and $b_2 = a_1^2 + 4 a_2$,
where $a_1$, $a_2$ and $b_4$ are parameters that are integers if the
model is integral. Since $v_2(c_4(t))=2$ and $8|24$, we obtain $v_2(b_2)=1$,
but, in view of $b_2 = a_1^2 + 4 a_2$, this is impossible. The model
cannot be integral. Let us pass to $c_4' = 2^4 c_4(t)$, $c_6' = 2^6 c_6(t)$,
$\Delta' = 2^{12} \Delta(t)$, which is certainly integral (because
it can be written as $y^2 = x^3 - 27 c_4(t) - 54 c_6(t)$) 
and thus minimal integral.
The new model satisfies
\[v_2(c_4')=6,\;\;\;\;\;\; v_2(c_6')\geq 11,\;\;\;\;\;\; v_2(\Delta') = 12.\]
Since $c_4 = f_1 f_2 (f_1^3 - f_2^3)^2$, we have 
\[\frac{1}{2^6} c_4' = \frac{1}{4} c_4(t) \equiv 1 \cdot (-1) \cdot 
\left(\frac{1^3 - (-1)^3}{2}\right)^2 \equiv 3 \mo 4.\]
Now the second line marked $I_2^*$ in \cite[Table 1]{Ha} gives us
$W_2(\mathscr{E}(t)) = -1$.

Lastly, let us consider the case $v_2(t)<0$. Let $k = -v_2(t)$. Then
$v_3(f_1(t))= - 2 k + 1$, $v_3(f_2(t)) = -2 k$, $v_3(f_1^3(t) - f_2^3(t))
= - 6 k$ and $v_3(f_1^3(t) + f_2^3(t)) = - 6 k$. Hence
\[v_2(c_4(t)) = -16 k + 1,\;\;\;\;\;\; v_2(c_6(t)) = -24 k - 1,\;\;\;\;\;\;
v_2(\Delta(t)) = - 48 k - 8.\]
We set $c_4' = (2^{4 k + 1})^4 c_4(t)$, $c_6' = (2^{4 k + 1})^6 c_6(t)$, 
$\Delta' = (2^{4 k + 1})^{12} \Delta(t)$, and obtain that
\begin{equation}\label{eq:dududu}
v_2(c_4') = 5,\;\;\;\;\; v_2(c_6') = 5,\;\;\;\;\;\; v_2(\Delta')
= 4.\end{equation}
This new model is given by the equation
$y^2 = x^3 - 27 \frac{c_4'}{2^4} - 54 \frac{c_6'}{2^6}$, which,
by (\ref{eq:dududu}), has all of its coefficients in $\mathbb{Z}_2$.
Hence the new model is (minimal) integral.

A brief calculation gives $2^{-5} c_6' = 2^{24 k +1} c_6(t) \equiv -1
\mo 4$. Now the second line from the top in \cite[Table 1]{Ha} gives us
$W_2(\mathscr{E}(t))=1$.


We conclude that, for $t\in \mathbb{Q}_2$,
\begin{equation}\label{eq:karma}W_2(\mathscr{E}(t)) = \begin{cases} -1 &\text{if $v_2(t)\geq 0$,}\\
1 &\text{if $v_2(t)<0$.}\end{cases}\end{equation}

\subsubsection{The root number at $p=3$} Let $t\in \mathbb{Q}_3$.
Suppose first that $v_3(t)>0$.
Then $v_3(f_1(t))=0$, $v_3(f_2(t))=0$, and -- since $f_1(t)\equiv 1 \mo 3$
and $f_2(t)\equiv 2 \mo 3$ -- $v_3(f_1(t)^3 - f_2(t)^3)=0$. 
On the other hand, because $(-5)^3 + 2^3 \equiv 0 \mo 9$, we have
$v_3(f_1(t)^3 + f_2(t)^3)\geq 2$. Therefore
\[v_3(c_4(t))=0,\;\;\;\;\; v_3(c_6(t)) \geq 2,\;\;\;\;\; v_3(\Delta(t))= -3.\]
Since $v_3(\Delta)<0$, we see that this model is not integral at $3$.
If we set $c_4' = 3^4 c_4(t)$, $c_6' = 3^6 c_6(t)$, $\Delta' = 3^{12} \Delta(t)$,
we obtain a new model for the same curve. This model has
\begin{equation}\label{eq:hoto}
v_3(c_4')=4,\;\;\;\;\; v_3(c_6') \geq 8,\;\;\;\;\; v_3(\Delta')= 9.\end{equation}
This new model can be written as $y^2 = x^3 - 27 c_4(t) - 54 c_6(t)$;
since all of the coefficients of this equation 
are in $\mathbb{Z}_3$, the new model is integral, and, in particular,
minimal integral.

By \cite[Table 2]{Ha}\footnote{
As was pointed out 
in \cite{Riz2}, the special condition $c_4'\equiv 1 \mo 4$ should be added
to line 3 in table 1 of \cite{Ha}.}, an elliptic curve $E$ with
a minimal integral model given by parameters
$c_4'$, $c_6'$, $\Delta'$ satisfying (\ref{eq:hoto}) has local root number
$W_3(E)=1$. Hence $W_3(\mathscr{E}(t))=1$ for all $t\in \mathbb{Q}_3$
with $v_3(t)>0$ (and, in particular, for all $t\in \mathbb{Q}$ with
$v_3(t)>0$).

Suppose now that $v_3(t)=0$. Then $t^2\equiv 1 \mod 3$. Hence
$v_3(f_1(t))=0$, $v_3(f_2(t))=0$, and -- since $f_1(t)\equiv 2 \mo 3$
and $f_2(t)\equiv 1 \mo 3$ -- $v_3(f_1(t)^3 - f_2(t)^3)=0$. 
On the other hand, because 
$(-5 - 2\cdot 4)^3 + (2 + 5\cdot 4)^3 \equiv 0 \mo 9$ and
$(-5 - 2\cdot (-4))^3 + (2 + 5\cdot (-4))^3 \equiv 0 \mo 9$, we have
$v_3(f_1(t)^3 + f_2(t)^3)\geq 2$. Proceeding as before, we get that
$c_4' = 3^4 c_4(t)$, $c_6' = 3^6 c_6(t)$, $\Delta' = 3^{12} \Delta(t)$ give us
a minimal integral model with
\[v_3(c_4')=4,\;\;\;\;\; v_3(c_6') \geq 8,\;\;\;\;\; v_3(\Delta')= 9,\]
and so, by the same entry in \cite[Table 2]{Ha} as before, the local root
number $W_3(E)$ at $3$ is $1$. We have shown that
 Hence $W_3(\mathscr{E}(t))=1$ for all $t\in \mathbb{Q}_3$
with $v_3(t)=0$.

Suppose, lastly, that $v_3(t)<0$. Let $k = -v_3(t)$. Then
$v_3(f_1(t)) = - 2 k$, $v_3(f_2(t)) = - 2 k$, $v_3(f_1^3(t) - f_2^3(t)) =
-6 k$ (since $3^{2 k} f_1(t) \equiv 1 \mo 3$ and 
$3^{2 k} f_2(t) \equiv 2 \mo 3$) and $v_3(f_1^3(t) + f_2^3(t)) \geq
-6 k + 2$ (because $(-2)^3 + 5^3 \equiv 0 \mo 9$). Hence
\[v_3(c_4(t)) = -16 k,\;\;\;\;\; v_3(c_6(t)) \geq - 24 k + 2,\;\;\;\;\;\;\;
v_3(\Delta(t)) = -48 k - 3.\]
We set $c_4' = (3^{4 k+1})^4 c_4(t)$, $c_6' = (3^{4 k+1})^6 c_6(t)$,
$\Delta' = (3^{4 k+1})^{12} \Delta(t)$, and obtain that
\[v_3(c_4') = 4,\;\;\;\;\; v_3(c_6') \geq 8,\;\;\;\;\;\; v_3(\Delta') = 9.\]
The model given by $c_4'$, $c_6'$ and $\Delta'$ is a minimal integral
model by the same argument as before. Hence, by the same entry
in \cite[Table 2]{Ha} as before, $W_3(\mathscr{E}(t))=1$ for all
$t\in \mathbb{Q}_3$ with $v_3(t)<0$.

We conclude that \begin{equation}\label{eq:mary}
W_3(\mathscr{E}(t))=1\end{equation} for all $t\in \mathbb{Q}_3$.

\vskip 15pt
\begin{center}
* * *
\end{center}

The remaining cases tend to be a little easier -- it is very simple to
tell what the reduction type is at a prime $p\ne 2,3$. Let us start
by considering $p=7$, since, among the cases that remain, it is the
most complicated one.

\subsubsection{The root number at $p=7$}\label{subs:furt} Suppose first that $v_7(t)>0$. Then
$v_7(f_1(t))=0$, $v_7(f_2(t))=0$, $v_7(f_1^3(t) - f_2^3(t)) = 1$,
$v_7(f_1^3(t) + f_2^3(t))=0$. Hence 
\begin{equation}\label{eq:saltor}
v_7(c_4(t)) = 2,\;\;\;\;\;\; v_7(c_6(t)) = 3,\;\;\;\;\;\; 
v_7(\Delta(t)) = 8.\end{equation}
 For valuations $v$ at places of $\mathbb{Q}$ other
than $2$ or $3$, the inequalities $v(c_4)\geq 0$, $v(c_6)\geq 0$, 
$v(\Delta)\geq 0$ guarantee that the model is integral; the inequality
$v(\Delta)<12$ guarantees that, if the model is integral, it is minimal
integral. Hence our model (given by $c_4(t)$, $c_6(t)$, $\Delta(t)$)
is minimal integral at $v_7$. By (\ref{eq:saltor}), 
the reduction type (see \S \ref{sss:piu})
is additive and potentially multiplicative.
We use (\ref{eq:totoru}) and obtain that
\[W_7(\mathscr{E}(t)) = (-1/7) = -1 \]
for all $t\in \mathbb{Q}_7$ with $v_7(t)>0$.

Suppose now that $v_7(t)=0$, $t\not\equiv \pm 1\mo 7$. Then
$v_7(f_1(t))=0$, $v_7(f_2(t))=0$, 
$f_1(t)\equiv f_2(t) \mo 7$, and so
$v_7(f_1^3(t) + f_2^3(t)) = 0$.
By (\ref{eq:oneil}),
$v_7(f_1^3(t) - f_2^3(t)) = 1 + v_7((5 + 4 t^2 + 5 t^4) (1 + t^2))
= 1$. We conclude that
\[v_7(c_4(t)) = 2,\;\;\;\;\;\; v_7(c_6(t)) = 3,\;\;\;\;\;\; 
v_7(\Delta(t)) = 8,\]
and so, as before,
\[W_7(\mathscr{E}(t)) = (-1/7) = -1 \]
for all $t\in \mathbb{Q}_7$ with $v_7(t)=0$, 
$t\not\equiv \pm 1\mo 7$. 

Consider now the case of $v_7(t)=0$, $t\equiv \pm 1 \mo 7$.
We get easily that
$v_7(f_1(t))\geq 1$, $v_7(f_2(t))\geq 1$ and
 $v_7(f_1^3(t) + f_2^3(t))\geq 3$.
Write $t = 7 s \pm 1$, $s\in \mathbb{Z}_p$. Then 
\[(19 + 11 t^2 + 19 t^4) = 7^3 (133 s^4 + 76 s^3 + 18 s^2 + 2 s) + 7^2
(-s^2 + 1).\]

Thus (see (\ref{eq:oneil})), if $s\not\equiv \pm 1 \mo 7$,
\begin{equation}\label{eq:ofrance}v_7(f_1^3(t) - f_2^3(t)) = 3,\end{equation}
and then
\[v_7(c_4(t))\geq 2 + 2\cdot 3 = 8,\;\;\;\;
  v_7(c_6(t))\geq 3 + 3\cdot 3 = 12,\;\;\;\;
  v_7(\Delta(t))\geq 8 \cdot 3 = 24.\]
We set $c_4' = (7^2)^{-4} c_4(t)$, $c_6' = (7^2)^{-6} c_6(t)$,
$\Delta' = (7^2)^{-12} \Delta(t)$.
Then
\[v_7(c_4')\geq 0,\;\;\;\;\;\; v_7(c_6')\geq 0,\;\;\;\;\;\;
v_7(\Delta')=0.\]
This is a minimal integral model. (At a place $p\ne 2,3$, the model
is integral if and only if $v_p(c_4)$ and $v_p(c_6)$ are both non-negative.)
Hence the curve has good reduction at $7$, and so
\[W_7(\mathscr{E}(t)) = 1\]
for all $t\in \mathbb{Q}_7$ with $v_7(t)=0$, 
$t\equiv \pm 1\mo 7$, $t\not\equiv \pm 7 \pm 1 \mo 49$. 

Suppose now that $s\equiv \pm 1 \mo 7$. Then
\[v_7(f_1^3(t) - f_2^3(t)) \geq 4.\]
On the other hand, a brief calculation shows that now
$v_7(f_1(t))=1$, $v_7(f_2(t))=1$ and $v_7(f_1^3(t) + f_1^3(t))=3$.
Define $k = v_7(f_1^3(t) - f_2^3(t))$.
Then
\[
v_7(c_4(t))\geq 2 + 2\cdot k,\;\;\;\;\;\;
v_7(c_6(t))\geq 3 + 3\cdot k,\;\;\;\;\;\;
v_7(\Delta(t))\geq 8\cdot k\]
Hence, if $k$ is even, the reduction of $\mathscr{E}(t)$ is additive
and potentially multiplicative, and so
\[W(\mathscr{E}(t)) = (-1/7) = -1;\]
if $k$ is odd, the reduction of $\mathscr{E}(t)$ is multiplicative, and so,
by (\ref{eq:shackle}),
\[W(\mathscr{E}(t)) = - ((- 7^{-v_7(c_6(t))} c_6(t))/7) .\]
This is all about the particular case $t\equiv \pm 7\pm 1 \mo 49$.

It remains to consider the case 
$v_7(t)<0$. Let $k = - v_7(t)$. Then $v_7(f_1(t)) = - 2 k$,
$v_7(f_2(t)) = - 2 k$, 
$v_7(f_1^3(t) + f_2^3(t)) = -6 k$,
$v_7(f_1^3(t) - f_2^3(t)) = -6 k + 1$. Hence
\[v_7(c_4(t))=  - 16 k + 2,\;\;\;\;\;\;
v_7(c_6(t))=  - 24 k + 3,\;\;\;\;\;\;
v_7(\Delta(t))=  - 48 k + 8,\;\;\;\;\;\;
\]
and so a minimal integral model has parameters $c_4'$, $c_6'$, $\Delta'$ 
with
\[v_7(c_4(t)) = 2,\;\;\;\;\;\;
v_7(c_6(t)) = 3,\;\;\;\;\;\;
v_7(c_6(t)) = 8 .\]
Thus, the reduction of
$\mathscr{E}(t)$ is additive
and potentially multiplicative, and so
\[W(\mathscr{E}(t)) = (-1/7) = -1.\]

We summarise: for $t\in \mathbb{Q}_7$, unless $v_7(t)=0$ and $t\equiv \pm 1
\mo 7$, we have
\begin{equation}\label{eq:cojo}
W_7(\mathscr{E}(t))= (-1/7) = -1.\end{equation}
If $v_7(t)=0$ and $t\equiv \pm 1 \mo 7$,
\begin{equation}\label{eq:manco}W_7(\mathscr{E}(t)) = 
\begin{cases} 
1 &\text{if $v_7(u(t)) = 3$,}\\
(-1/7) = -1 &\text{if $v_7(u(t))\geq 4$ and $v_7(u(t))$ is even,}\\
-(- 7^{-v_7(c_6(t))} c_6(t)/7) &\text{if $v_7(u(t))\geq 5$ and $v_7(u(t))$ is odd,}
\end{cases}\end{equation}
where we define $u(t) = f_1^3(t) - f_2^3(t)$. 


\subsubsection{The root number at $p=5$} Assume first that $v_5(t)>0$. Then
$v_5(f_1(t))=1$, $v_5(f_2(t))=0$, $v_5(f_1^3(t) - f_2^3(t))=0$.
Hence the reduction is good, and so $W_5(\mathscr{E}(t))=1$.

Suppose now that $v_5(t)=0$. Then $v_5(f_1(t))=0$ and $v_5(f_2(t))=0$.
Since $19 + 11 t^2 + 19 t^4 \equiv 0 \mo 5$ has no roots in
$\mathbb{Z}/4\mathbb{Z}$, (\ref{eq:oneil}) implies that
\[v_5(f_1^3(t)-f_2^3(t)) = v_5(1+t^2).\]
Now, $v_5(1+t^2)>0$ implies $t\equiv \pm 2 \mo 5$, and this, in turn,
implies that $v_5(f_1^3(t) + f_2^3(t))=0$.
Thus, the reduction is additive and potentially multiplicative if
$v_5(f_1^3(t)-f_2^3(t))$ is odd, and multiplicative if
$v_5(f_1^3(t)-f_2^3(t))$ is positive and even. By Prop.\ \ref{prop:rohr},
this means that 
\[W_5(\mathscr{E}(t)) = \begin{cases}
1 &\text{if $v_5(f_1^3(t)-f_2^3(t))=0$,}\\
(-1/5)=1 &\text{if $v_5(f_1^3(t)-f_2^3(t))$ is odd and positive,}\\
-(-5^{v_5(c_6(t))} c_6(t)/5) 
&\text{if $v_5(f_1^3(t)-f_2^3(t))$ is even and positive}
\end{cases}\]
for all $t\in \mathbb{Q}_5$ with $v_5(t)=0$.

Lastly, suppose $v_5(t)<0$. Let $k = -v_5(t)$. Then $v(f_1(t))=-2 k$,
$v(f_2(t)) = -2 k+1$, $v(f_1^3(t) + f_2^3(t))= - 6 k$,
$v(f_1^3(t) - f_2^3(t))= - 6 k$. Hence
\[v_5(c_4(t)) = -16 k + 1,\;\;\;\;\;\;
v_5(c_6(t)) = - 24 k,\;\;\;\;\;\;
v_5(\Delta(t)) = - 48 k.\]
Thus, the reduction at $p=5$ is good, and $W(\mathscr{E}(t)) = 1$.

We conclude that, for all $t\in \mathbb{Q}_5$,
\[W_5(\mathscr{E}(t)) = \begin{cases}
-(-5^{-v_5(c_6(t))} c_6(t)/5) 
&\text{if $v_5(f_1^3(t)-f_2^3(t))$ 
is even and positive,}\\
(-1/5) = 1 
&\text{if $v_5(f_1^3(t)-f_2^3(t))$ 
is odd and positive,}\\
1&\text{if $v_5(f_1^3(t)-f_2^3(t))\leq 0$.}
\end{cases}\]
As we will later see, this is exactly the sort of thing that happens
for every $p$ outside the finite set $\{2,3,5,7,13,19\}$.

\subsubsection{The root number at $p=13$} 
Assume first that $v_{13}(t)\geq 0$. Then $v_{13}(19+11 t^2 + 19 t^4) = 0$,
and so $v(f_1^3(t)-f_2^3(t)) = v(1+t^2)$. If $v(1+t^2)>0$, then
$v(f_1(t))=0$, $v(f_2(t))=0$ and $v(f_1^3(t)+f_2^3(t))=0$. Hence,
the reduction at $p=13$ is good if
$v_{13}(f_1^3(t) - f_2^3(t))=0$,
additive and potentially multiplicative if
$v_{13}(f_1^3(t)-f_2^3(t))$ is odd, and multiplicative if
$v_{13}(f_1^3(t)-f_2^3(t))$ is positive and even. By Prop.\ \ref{prop:rohr},
this means that 
\[W_{13}(\mathscr{E}(t)) = \begin{cases}
1 &\text{if $v_{13}(f_1^3(t)-f_2^3(t))=0$,}\\
(-1/{13})=1 &\text{if $v_{13}(f_1^3(t)-f_2^3(t))$ is odd and positive,}\\
-(-{13}^{v_{13}(c_6(t))} c_6(t)/{13}) 
&\text{if $v_{13}(f_1^3(t)-f_2^3(t))$ is even and positive}
\end{cases}\]
for all $t\in \mathbb{Q}_{13}$ with $v_{13}(t)\geq 0$.

Assume, lastly, that $v_{13}(t)<0$. Let $k = - v_{13}(t)$. Then
$v_{13}(f_1(t)) = - 2 k$, $v_{13}(f_2(t)) = -2 k$,
$v_{13}(f_1(t)^3 + f_2(t)^3) = - 6 k + 1$,
$v_{13}(f_1(t)^3 - f_2(t)^3) = - 6 k$. Thus,
\[v(c_4(t)) = -16 k,\;\;\;\;\;\;
v(c_6(t)) = - 24 k + 1,\;\;\;\;\;
v(\Delta(t)) = -48 k.\]
Thus, the reduction at $13$ is good, and so $W_{13}(\mathscr{E}(t))=1$.

We conclude that 
\[W_{13}(\mathscr{E}(t)) = \begin{cases}
-(-{13}^{-v_{13}(c_6(t))} c_6(t)/{13}) 
&\text{if $v_{13}(f_1^3(t)-f_2^3(t))$ is even and positive,}\\
(-1/{13}) =1 &\text{
if $v_{13}(f_1^3(t)-f_2^3(t))$ is odd and positive,}\\
1& \text{if $v_{13}(f_1^3(t)-f_2^3(t))\leq 0$}
\end{cases}\]
for all $t\in \mathbb{Q}_{13}$ with $v_{13}(t)=0$. Again, here as for
$p=5$, the behaviour of the local root number is as it is at a prime
$p\notin \{2,3,5,7,13,19\}$.

\subsubsection{The root number at $p=19$} Assume first that $v_{19}(t)> 0$.
 Then
$v_{19}(f_1(t))=0$, $v_{19}(f_2(t))=0$, $v_{19}(f_1^3(t)+f_2^3(t))=0$
and (by (\ref{eq:oneil})) $v_{19}(f_1^3(t)-f_2^3(t))=1$. Hence
\[v_{19}(c_4(t)) = 2,\;\;\;\;\;
v_{19}(c_6(t)) = 3,\;\;\;\;\;
v_{19}(\Delta(t)) = 8.\]
Thus, the reduction type is additive and potentially multiplicative.
The root number is hence
\[W_{19}(\mathscr{E}(t)) = (-1/19) = -1.\]

Suppose now that $v_{19}(t)=0$. Because $-1$ is not a quadratic residue
$\mo 19$, we know that $v_{19}(f_1^3(t)-f_2^3(t))=0$, and thus
$v_{19}(\Delta(t)) = 0$. At the same time, $v_{19}(c_4(t))\geq 0$
and $v_{19}(f_6(t))\geq 0$. Hence, the reduction at $19$ is good, and so
$W_{19}(\mathscr{E}(t)) = 1$.

Suppose now that $v_{19}(t)<0$. Let $k = - v_{19}(t)$. Then
$v_{19}(f_1(t)) = - 2 k$, $v_{19}(f_2(t))= - 2 k$,
$v_{19}(f_1^3(t) - f_2^3(t)) = - 6 k$, 
$v_{19}(f_1^3(t) - f_2^3(t)) = -6 k + 1$.
Hence
\[v_{19}(c_4(t)) = -16 k + 2,\;\;\;\;\;
v_{19}(c_6(t)) = - 24 k + 3,\;\;\;\;\;
v_{19}(\Delta(t)) = - 48 k + 8.\]
Thus, the reduction type is additive and potentially multiplicative.
We obtain
\[W_{19}(\mathscr{E}(t)) = (-1/19) = -1.\]

In general, we can now conclude, 
\begin{equation}\label{eq:peter}W_{19}(\mathscr{E}(t)) = \begin{cases}
-1 &\text{if $v_{19}(t)\ne 0$,}\\
1 &\text{if $v_{19}(t)=0$}\end{cases}\end{equation}
for all $t\in \mathbb{Q}_p$. 

\subsubsection{The local root number at a prime $p\notin \{2,3,5,7,13,19\}$}
Assume first that $v_p(t)\geq 0$. Let $k = v_p(f_1^3(t) - f_2^3(t))$.
Since only primes in $\{2,3,5,7,13,19\}$ can derive a resultant
$\Res(P_i,P_j)$, $1\leq i<j\leq 8$ (see (\ref{eq:gosto})),
we conclude that, for $k\geq0$, 
\[v_p(c_4(t)) = 2 k,\;\;\;\;\;
v_p(c_6(t)) = 3 k,\;\;\;\;\; v_p(\Delta(t)) = 8 k.\]
Hence,
the reduction at $p$ is good if
$v_{p}(f_1^3(t) - f_2^3(t))=0$,
additive and potentially multiplicative if
$v_{p}(f_1^3(t)-f_2^3(t))$ is odd, and multiplicative if
$v_{p}(f_1^3(t)-f_2^3(t))$ is positive and even. Thus
\[W_p(\mathscr{E}(t)) = \begin{cases}
-(-p^{v_p(c_6(t))} c_6(t)/p) 
&\text{if $v_{p}(f_1^3(t)-f_2^3(t))$ is even and positive,}\\
(-1/p)^{v_{p}(f_1^3(t)-f_2^3(t))} &\text{otherwise}
\end{cases}\] 
for all $t\in \mathbb{Z}_p$.

Assume now that $v_p(t)<0$. Let $k = - v_p(t)$. Then
\[v_{19}(c_4(t)) = -16 k,\;\;\;\;\;
v_{19}(c_6(t)) = - 24 k,\;\;\;\;\;
v_{19}(\Delta(t)) = - 48 k.\]
Thus, the reduction at $p$ is good, and so
\[W(\mathscr{E}(t))=1 \] for all 
$t\in \mathbb{Q}_t$ with $v_p(t)<0$.

We conclude that, for $p\notin \{2,3,5,7,13,19\}$,
\begin{equation}\label{eq:boclob}W_p(\mathscr{E}(t)) = \begin{cases}
-(-p^{-v_p(c_6(t))} c_6(t)/p) 
&\text{if $v_{p}(f_1^3(t)-f_2^3(t))$ is even and positive,}\\
(-1/p) &\text{if
$v_{p}(f_1^3(t)-f_2^3(t))$ is odd and positive,}\\
1 & \text{if $v_{p}(f_1^3(t)-f_2^3(t))\leq 0$}
\end{cases} \end{equation}
for all $t\in \mathbb{Q}_p$. As we saw before, this is also true for
$p=5, 13$.

\subsection{The global root number} 
We must compute the global root number
\[W(\mathscr{E}(t)) = \prod_v W_v(\mathscr{E}(t)) =  - \prod_p W_p(
\mathscr{E}(t)).\]

First, we separate what we may see as the ``main term'', speaking loosely.
This is the contribution that each prime $p\ne 2,3,7,19$ would make if
$v_p(f_1^3(t) - f_2^3(t))$ were $\leq 1$, as is the case for all but finitely
many primes. For each prime $p\ne 2,3,7,19$ with $v_p(f_1^3(t)-f_2^3(t))\leq 1$,
\begin{equation}\label{eq:osorto}W_p(\mathscr{E}(t)) = \begin{cases}
(-1/p) &\text{if $v_p(f_1^3(t) - f_2^3(t))$ is odd and positive,}\\
1 &\text{otherwise},\end{cases}\end{equation}
by our case-work above. Thus, we would like to determine
\[V = \mathop{\prod_{p\ne 2,3,7,19}}_{v_p(f_1^3(t) - f_2^3(t))\geq 0} 
(-1/p)^{v_p(f_1^3(t) - f_2^3(t))}\]
for all $t\in \mathbb{Q}$.

Write $t = \frac{t_0}{t_1}$, where $t_0$ and $t_1$ are coprime
integers.
Then
\[\begin{aligned}
f_1^3(t) - f_2^3(t) &= (-5 - 2 (t_0/t_1)^2)^3 - (2 + 5 (t_0/t_1)^2)^3\\
&= \frac{1}{t_1^6} \cdot -7 (19 t_1^4 + 11 t_1^2 t_0^2 +19 t_0^4) 
(t_1^2 + t_0^2) .
\end{aligned}\]
Let 
\[P(x,y) = -7 (19 y^4 + 11 y^2 x^2 +19 x^4) 
(y^2 + x^2) .\]
If $p|y$ and $p\ne 7,19$, then, since $x$ and $y$ are coprime,
\begin{equation}\label{eq:oportor}v_p(7 (19 y^4 + 11 x^2 y^2 +19 x^4) 
(y^2 + x^2)) = 0.\end{equation}
Hence, if $v_p(P(t_0,t_1))>0$, then $p\nmid t_1$. It follows that,
for $p\ne 7,19$, $v_p(f_1^3(t) - f_2^3(t))\geq 0$ implies
$v_p(f_1^3(t) - f_2^3(t)) = v_p(P(t_0,t_1))$,
and, again for $p\ne 7,19$, $v_p(P(t_0,t_1))>0$ implies
$v_p(f_1^3(t)-f_2^3(t)) \geq 0$ and $v_p(f_1^3(t) - f_2^3(t)) = v_p(P(t_0,t_1))$.
Hence
\begin{equation}\label{eq:tata}V = \prod_{p\ne 2,3,7,19} 
(-1/p)^{v_p(P(t_0,t_1))}.\end{equation}
(In the notation we introduced in \S  (\ref{eq:brack}), this
can be rewritten as
\[V = (-1|P(t_0,t_1))_d,\]
where $d = 2\cdot 3\cdot 7\cdot 19$.)

Here comes a key step.
By the quadratic reciprocity law,
\begin{equation}\label{eq:otorongo}
\prod_{p\ne 2} (-1/p)^{v_p(P(t_0,t_1))} = \begin{cases}
1 &\text{if $a\equiv 1 \mo 4$,}\\
-1 &\text{if $a\equiv -1 \mo 4$,}\end{cases}\end{equation}
where $a = \prod_{p\ne 2} p^{v_p(P(t_0,t_1))}$. (Actually,
what we need here is not the full law of quadratic reciprocity, 
but only the fact that
$(-1/p) = 1$ if and only if $p\equiv 1 \mo 4$.)

Since $t_0$ and $t_1$ are coprime, they cannot both be even. If one
of them is even, then $2\nmid P(t_0,t_1)$, and so
$a$ is the absolute value $|P(t_0,t_1)|$ of $P(t_0,t_1)$.
In that case,
\[\begin{aligned}
|P(t_0,t_1)| &=  7 (19 t_1^4 + 11 t_0^2 t_1^2 +19 t_0^4) (t_1^2 + t_0^2)\\
&\equiv 3 (3 t_1^4 + 3 t_0^2 t_1^2 + 3 t_0^4) (t_1^2 + t_0^2) \mo 4\\
&\equiv 1 \mo 4.
\end{aligned}\]
Thus, by (\ref{eq:otorongo}),
\[\prod_{p\ne 2} (-1/p)^{v_p(P(t_0,t_1))} = 1\]
when exactly one of $t_0$, $t_1$ is odd and the other one is even.

Suppose now that $t_0$ and $t_1$ are both odd. Then $t_1^2 + t_0^2
\equiv 2 \mo 4$, and so $|P(t_0,t_1)| \equiv 2 \mo 4$. Now
\[\begin{aligned}
\frac{|P(t_0,t_1)|}{2}
 &=  7 (19 t_1^4 + 11 t_0^2 t_1^2 +19 t_0^4) \cdot \frac{t_1^2 + t_0^2}{2}\\
&\equiv 3 (3 t_1^4 + 3 t_0^2 t_1^2 + 3 t_0^4) \cdot 1 \mo 4\\
&\equiv 3 \mo 4.\end{aligned}\]
Thus, by (\ref{eq:otorongo}),
\[\prod_{p\ne 2} (-1/p)^{v_p(P(t_0,t_1))} = -1\]
when $t_0$, $t_1$ are both odd.

Looking back at (\ref{eq:tata}), we see that
\begin{equation}\begin{aligned}V &= \prod_{p\ne 2,3,7,19} 
(-1/p)^{v_p(P(t_0,t_1))} = 
\prod_{p=3,7,19} (-1/p)^{v_p(P(t_0,t_1))} \cdot \prod_{p\ne 2}
(-1/p)^{v_p(P(t_0,t_1))}\\
&= \prod_{p=3,7,19} (-1/p)^{v_p(P(t_0,t_1))} \cdot 
\begin{cases} -1 &\text{if $t_0$, $t_1$ are both odd}\\
1 &\text{otherwise.}\end{cases}
\end{aligned}\end{equation}
We can now write
\begin{equation}\label{eq:gotor}\begin{aligned}
W(\mathscr{E}(t)) &= - \prod_p W_p(\mathscr{E}(t))\\
&= - \prod_{p=2,3,7,19} W_p(\mathscr{E}(t)) \cdot
\prod_{p\ne 2,3,7,19} W_p(\mathscr{E}(t)) \\
&= - \prod_{p=2,3,7,19} W_p(\mathscr{E}(t)) \cdot
\prod_{p\ne 2,3,7,19} \frac{W_p(\mathscr{E}(t))}{(-1/p)^{v_p(P(t_0,t_1)}}
\cdot V \\
&= - W_2(\mathscr{E}(t)) 
\cdot \prod_{p=3,7,19} W_p(\mathscr{E}(t))
\cdot (-1/p)^{v_p(P(t_0,t_1))}\\ &\cdot
\prod_{p\ne 2,3,7,19} \frac{W_p(\mathscr{E}(t))}{(-1/p)^{v_p(P(t_0,t_1))}}
\cdot\begin{cases} -1 &\text{if $t_0$, $t_1$ are both odd}\\
1 &\text{otherwise.}\end{cases}\end{aligned}\end{equation}
We define
\begin{equation}\label{eq:ostor}\begin{aligned}
g_\infty(x,y) &=-1\;\;\;\;\;\;\;\;\; \text{for all $x,y\in \mathbb{R}$},\\
g_2(x,y) &= W_2(\mathscr{E}(x/y)) \cdot
\begin{cases} -1 &\text{if $x$, $y$ are both odd}\\
1 &\text{otherwise.}\end{cases}\\
g_p(x,y) &= W_p(\mathscr{E}(x/y))
\cdot (-1/p)^{v_p(P(x,y))} 
\;\;\;\;\;
\text{for $p=3,7,19$}.\end{aligned}\end{equation}
and
\begin{equation}\label{eq:alhux}h_p(x,y) = 
\frac{W_p(\mathscr{E}(t))}{(-1/p)^{v_p(P(x,y))}}\end{equation}
for $p\neq 2,3,7,19$. 
By (\ref{eq:osorto}) and the discussion after (\ref{eq:oportor}), 
$h_p(x,y)=1$ whenever $p^2\nmid P(x,y)$.
We now rewrite (\ref{eq:gotor}) as follows:
\[W(\mathscr{E}(x/y)) = \prod_{v\in \{\infty,2,3,7,19\}} g_v(x,y)
\cdot \mathop{\prod_{p\notin \{2,3,7,19\}}}_{p^2|P(x,y)} h_p(x,y),\]
for all coprime $x$, $y$, 
where $P(x,y) = -7 (19 y^4 + 11 y^2 x^2 +19 x^4) 
(y^2 + x^2)$.
 This is simply a special case of Thm.\ \ref{thm:smet},
(\ref{eq:exerier}).

Why could we not have simply invoked Thm.\ \ref{thm:smet}? Theorem
\ref{thm:smet} does not give explicit expressions for $g_v$ and $h_p$.
We have found $g_v$ and $h_p$ in our example by, in effect,
 following the proofs of Prop.\ \ref{prop:rijk}, Prop.\ \ref{prop:ahem}
and Thm.\ \ref{thm:smet}.

The functions $g_v$, $h_p$ are given by the expressions in (\ref{eq:ostor})
and (\ref{eq:alhux}), combined with the work we did before on local
root numbers.
Let us now make these expressions more explicit. It is here that the results
of the lengthy casework in \S \ref{subs:lort} finally enter.

Let us first look at $p=2$. For $x$ and $y$ coprime, $v_2(x/y)>0$ if
$x$ is even and $y$ is odd, $v_2(x/y)=0$ if they are both odd and
$v_2(x/y)<0$ if $x$ is odd and $y$ is even.
Hence, by (\ref{eq:karma}) and (\ref{eq:ostor}),
\begin{equation}\label{eq:jada}g_2(x,y) = \begin{cases} -1 &\text{if $x$ is even and $y$ is odd,}\\
1 &\text{otherwise.}\end{cases}\end{equation}

Let us now look at $p=3$. By (\ref{eq:mary}) and (\ref{eq:ostor}),
\[
g_3(x,y) = W_3(\mathscr{E}(x/y))
\cdot (-1/3)^{v_3(P(x,y))} = (-1/3)^{v_3(P(x,y))}.\]
A brief examination suffices to show that $v_3(P(x,y))=0$ for all coprime
$x$, $y$, and so
\begin{equation}\label{eq:jede}g_3(x,y)=1\end{equation}
identically for all coprime $x$, $y$.

Consider now $p=7$. A brief calculation shows that
$v_7(P(x,y)) = v_7(-7 (19 y^4 + 11 y^2 x^2 +19 x^4) 
(y^2 + x^2)) = 1 + v_7(19 y^4 + 11 y^2 x^2 + 19 x^4)$.
If $v_7(19 y^4 + 11 y^2 x^2 + 19 x^4)=0$, then, by (\ref{eq:ostor}) and
(\ref{eq:cojo}),
\begin{equation}\label{eq:fra1}g_7(x,y) = W_7(\mathscr{E}(x/y))
\cdot (-1/7)^{v_{7}(P(x,y))} = (-1) \cdot (-1) = 1.\end{equation}

Suppose that $v_7(19 y^4 + 11 y^2 x^2 + 19 x^4)>0$. Then
$v_7(x/y)=0$ and  $x/y \equiv \pm 1 \mo 7$.
Write $x/y = 7 s \pm 1$, $s\in \mathbb{Z}_p$. 
If $s\not\equiv \pm 1 \mo 7$, we are in the case
of equation (\ref{eq:ofrance}), and so 
\begin{equation}\label{eq:fra2}g_7(x,y) = W_7(\mathscr{E}(x/y))
\cdot (-1/7)^{v_{7}(P(x,y))} = (1) \cdot (-1) = -1.\end{equation}
If $s\equiv \pm 1 \mo 7$, then $v_7(P(x,y))\geq 4$,
and so, by (\ref{eq:manco}), 
\begin{equation}\label{eq:fra3}g_7(x,y) = W_7(\mathscr{E}(x/y))
\cdot (-1/7)^{v_{7}(P(x,y))} = (-1) \cdot (-1)^{v_7(P(x,y))} = -1\end{equation}
if $v_7(P(x,y))$ is even,
\begin{equation}\label{eq:fra4}\begin{aligned}g_7(x,y) &= W_7(\mathscr{E}(x/y))
\cdot (-1/7)^{v_{7}(P(x,y))} = - (-7^{-v_7(c_6(x/y))} c_6(x/y)/7)
 \cdot (-1)^{v_7(P(x,y))}\\ &= (-7^{-v_7(c_6(x/y))} c_6(x/y)/7)\end{aligned}\end{equation}
if $v_7(P(x,y))$ is odd. 

Lastly, let us look at
 $p=19$. Looking at $P(x,y) = -7 (19 y^4 + 11 y^2 x^2 +19 x^4) 
(y^2 + x^2)$, we see that $v_{19}(P(x,y))=1$ if either $19|x$ or $19|y$
(and $x$ and $y$ are coprime) and $v_{19}(x,y)=0$ if $19\nmid x$
and $19\nmid y$. Hence, by (\ref{eq:peter}) and (\ref{eq:ostor}),
\begin{equation}\label{eq:jidi}g_{19}(x,y) = W_{19}(\mathscr{E}(x/y))
\cdot (-1/19)^{v_{19}(P(x,y))} = 1.\end{equation}

Having examined the functions $g_p$, $p\in \{2,3,7,19\}$, we now
look at the functions $h_p$, $p\not\in \{2,3,7,19\}$. By (\ref{eq:boclob})
and (\ref{eq:alhux}),
\begin{equation}\label{eq:rugir}
h_p(x,y) = \frac{W_{p}(\mathscr{E}(x/y))}{(-1/p)^{v_p(P(x,y))}}
= \begin{cases}
-(p^{-v_p(c_6(x/y))} c_6(x/y)/p)&\text{if $v_p(P(x,y))$ is even and $\geq 2$,}\\
1 &\text{otherwise.}
\end{cases}\end{equation}

\subsection{Computing $p$-adic integrals}

By (\ref{eq:jada}), (\ref{eq:jede}) and (\ref{eq:jidi}),
\[\begin{aligned} \int_{O_2} g_2(x,y) dx dy &= \frac{1}{4},\\
\int_{O_3} g_3(x,y) dx dy &= \Area(O_3) = 1 - \frac{1}{9},\;\;\;\;\;\;\;\;\;
\int_{O_{19}} g_{19}(x,y) dx dy = 1 - \frac{1}{19^2}.\end{aligned}\]
Recall that $O_p = (\mathbb{Z}_p \times \mathbb{Z}_p) \setminus
(p\mathbb{Z}_p \times p\mathbb{Z}_p)$.

Before looking at the integral $\int_{O_7} g_7(x,y) dx dy$ of $g_7(x,y)$,
 let us look at the integral of $h_p(x,y)$, $p\ne 2,3,7,19$.

If $p^2\nmid P(x,y)$, where $P(x,y) = -7 (19 y^4 + 11 y^2 x^2 +19 x^4) (x^2 + y^2)$, 
then $h_p(x,y)=1$. Assume $p^2|P(x,y)$. 
We consult the short list (\ref{eq:shortli}) of prime factors of
determinants and conclude that $p\nmid P_j(x,y)$ for every $1\leq j\leq 7$,
$j\ne 5,6$, where the polynomials $P_j$ are as in (\ref{eq:baboo}).
The assumption $p^2|P(x,y)$ also implies that $v_p(x/y) = 0$ 
and so $p\nmid P_8(x,y)=y$. We thus see from (\ref{eq:ferrat}) that,
for $(x,y)\in O_p$ with $p^2|P(x/y)$
and $x, y \mo p$ given, the quantity
$p^{-v_p(c_6(x/y))} c_6(x/y) \mo p$ equals a constant in $(\mathbb{Z}/p\mathbb{Z})^*$
times $p^{-3 v_p(P(x/y,1))} P(x/y,1)^3 \mod p$.  In particular,
$(p^{-v_p(c_6(x/y))} c_6(x/y)/p)$ equal a constant $\in \{-1,1\}$ times
$(p^{-v_p(P(x/y,1))} P(x/y,1) \mo p$.  
Let $f(t)=P(t,1)$.

Since $p\ne 2,3,7,19$, the discriminant of $f$ is not divisible by $p$.
Then $p^2|f(t)$ (or, for that matter, $p|f(t)$) implies $p\nmid f'(t)$.
Hence, for $t$ such that $v_p(t-t_0)>0$ for some root $t_0$ of $f=0$,
\[f(t) = f'(t) (t-t_0) + g(t) (t-t_0)^2,\]
where $v(g(t))\geq 0$, and so
\[p^{-v_p(t)} f(t) \equiv  f'(t) \cdot p^{-v_p(t-t_0)} (t-t_0) \mo p.\]
It follows that, as $t$ ranges over all values of $t\in \mathbb{Q}_p$
with $v_p(t-t_0)$ given, the integral of $(p^{-v_p(t)} f(t)/p)$ is $0$.
(Here we are using the fact that $\sum_{a=1}^{p-1} (a/p) = 0$.)
Since $p\nmid \Disc(f)$, the roots $t_0$ of $f=0$ are not congruent to 
each other modulo $p$, and so the integral of $(p^{-v_p(t)} f(t)/p)$
over all $t$ such that $v_p(f(t))=k$ ($k>0$ arbitrary) is zero.

Hence, the integral of $(p^{-v_p(c_6(x/y))} c_6(x/y) /p)\mo p$ over all
pairs $(x,y)\in O_p$ with $v_p(P(x,y))=k$, $k$ given and positive, is $0$.
Working from (\ref{eq:rugir}), we see that, for $p\ne 2,3,7,19$,
\[\begin{aligned}
\int_{O_p} h_p(x,y) &= \Area(O_p) -
a_p \cdot \mathop{\sum_{k\geq 2}}_{\text{$k$ even}} p^{-k} (1-p^{-1})\\ &=
(1 - p^{-2}) - a_p \frac{p^{-2}(1-p^{-1})}{1 - p^{-2}}  ,
\end{aligned}\]
where $a_p$ is the number of roots of $f(t)=0$ in $\mathbb{Q}_p$. Since
$p\nmid \Disc(f)$, $a_p$ equals the number of roots of the equation
\[(19 t^4 + 11 t^2 +19)(t^2+1) \equiv 0 \mo p\]
in $\mathbb{Z}/p\mathbb{Z}$.

We can now go back to $p=7$. It is much the same story, only that,
since $p$ does divide $\Disc(f)$, we need to treat $v_p(f(t))$ positive
and small as a collection of special cases - and that is exactly
what we did in \S \ref{subs:furt}. (Because $p|\Disc(f)$, some of the
roots $t_i$, $t_j$ of $f(t)=0$ are close to each other: $v(t_i-t_j)>0$.
We saw in \S \ref{subs:furt} that the roots are not very close to each
other: there are four distinct roots $t_1$, $t_2$, $t_3$, $t_4$, and they
satisfy $v(t_i-t_j)$ is $1$ or $2$ for all $1\leq i<j\leq 4$.)
 
Looking at (\ref{eq:fra1})--(\ref{eq:fra4}) and working as we did for $h_p$, we see that
\[\begin{aligned}
\int_{O_7} g_7(x,y) &= 1\cdot (1 - 7^{-2} - 2 (1-7^{-1}) 7^{-1})
+ (-1) \cdot (2 (1-7^{-1}) 7^{-1} (1 - 2\cdot 7^{-1})) \\
&+ (-1) \cdot 4 (1-7^{-1}) 7^{-2}\cdot \mathop{\sum_{k\geq 0}}_{\text{$k$
    even}}
(1 - 7^{-1}) 7^{-k} \\ &+ 
  0 \cdot 4 (1-7^{-1}) 7^{-2}\cdot \mathop{\sum_{k\geq 0}}_{\text{$k$
    odd}} (1 - 7^{-1}) 7^{-k} \\
&=  (1 - 2\cdot 7^{-1} + 7^{-2}) + (-1)\cdot (2\cdot 7^{-1} - 6\cdot 7^{-2} + 4\cdot 7^{-3})\\
&+ (-1) \cdot 4 \frac{(1 - 7^{-1})^2 7^{-2}}{1 - 7^{-2}}\\
&= 1 - 4\cdot 7^{-1} +7\cdot 7^{-2} - 4\cdot 7^{-3} - \frac{4 (1-7^{-1})^2 7^{-2}}{1-7^{-2}}\\
&= \frac{171}{343} = \frac{3^2\cdot 19}{7^3}. 
\end{aligned}\]

Lastly, recall that, by (\ref{eq:ostor}),  $g_{\infty}(x,y) = -1$ identically.
\subsection{Conclusion}

We can now apply Proposition \ref{prop:cookie}.
(Our polynomial \[B_{\mathscr{E}} = (19 y^4 + 11 y^2 x^2 + 19 x^4) (y^2 +
x^2)\]
 is a product of irreducible polynomials of degree
$\leq 4$.  Thus, by \cite{He},
Prop.\ 4.12, hypothesis $\mathscr{A}_2(B_{\mathscr{E}})$ holds, and thus
Proposition \ref{prop:cookie} holds unconditionally.)
We obtain
\[\begin{aligned}
\av_{\mathbb{Q}} W(\mathscr{E}(x/y)) &= -1 \cdot \frac{1/4}{1-2^{-2}}
\cdot \frac{1-1/9}{1-3^{-2}} \cdot \frac{171/343}{1 - 7^{-2}}\cdot
\frac{1 - 1/19^2}{1-19^{-2}} \\ &\cdot
\prod_{p\ne 2,3,7,19} \frac{1 - p^{-2} - a_p
  \frac{p^{-2}}{1-p^{-2}} (1-p^{-1})}{1-p^{-2}},\end{aligned}\]
i.e.,
\begin{equation}\label{eq:costro}
\av_{\mathbb{Q}} W(\mathscr{E}(x/y)) = 
 \frac{-19}{112} \cdot \prod_{p\ne 2,3,7,19} \left(1 - a_p \frac{p^{-2}
 (1-p^{-1})}{(1 - p^{-2})^2}\right),\end{equation}
where $a_p$ is the number of roots of
the equation
\[(19 t^4 + 11 t^2 +19) (t^2+1)\equiv 0 \mo p\]
in $\mathbb{Z}/p\mathbb{Z}$.

Numerically,
\begin{equation}\label{eq:umi}
\av_{\mathbb{Q}} W(\mathscr{E}(x/y)) = -0.15294\dotsc\end{equation}

We can compare this asymptotic result to the average of $W(\mathscr{E}(x/y))$
over all rationals $x/y$, $\gcd(x,y)=1$, $|x|, |y|\leq N$. A program running
on SAGE gives an average of $-0.15529\dotsc$ for $N=500$ and an average of
$-0.155428$ for $N=1500$.

\section{Two-parameter families}\label{sec:goron}
Every elliptic curve $E$ over $\mathbb{Q}$ can be written in the form
$E:y^2 = x^3 + a x + b$, with $a$, $b$ integral. 
Thus, the most natural two-parameter family
of elliptic curves over $\mathbb{Q}$ is $\mathscr{E}:(a,b)\mapsto
(y^2 = x^3 + a x + b)$. Unfortunately, the
problem of averaging the root number over that family $\mathscr{E}$ seems
to necessitate solving an open case of the parity problem; to put
a remark in \cite{Wo}, \S 1, in the language of this paper, we would
in all likelihood need to show that $\lambda(a^2 - b^3)$ averages to 
zero as $a$ and $b$ vary over $\mathbb{Z}$. In the face of \cite{FI}, one 
may think that showing as much is not completely out of reach, but, as of now, 
proving that $\lambda(a^2 - b^3)$ averages to zero presents, at the
very least, some serious technical difficulties, and thus remains open.

We can, however, give unconditional statements for some other 
two-parameter families. Consider the elliptic curves of the form
 $\mathscr{E}(a,b): y^2 = x (x+a) (x+ b)$ for $a,b\in \mathbb{Z}$. 
This family is of some interest, as every
semisimple elliptic curve over $\mathbb{Q}$ with full rational $2$-torsion
$(\mathbb{Z}/2 \mathbb{Z})^2$ is isomorphic over $\mathbb{Q}$ 
to $\mathscr{E}(a,b)$
for some $a,b\in \mathbb{Z}$. We can, in fact, assume that
$\gcd(a,b)$ is square-free, as $\mathscr{E}(a p^2, b p^2)$ is
isomorphic to $\mathscr{E}(a,b)$. 

We will now show that $W(\mathscr{E}(a,b))$ averages to zero. 
Here $a$, $b$ will go through all non-zero integers with
$\gcd(a,b)=1$. The condition $\gcd(a,b)=1$ is there for the sake
of simplicity; the same methods would suffice to compute the average
over all $a,b\in \mathbb{Z}$ with some more work.

The proof will give us another opportunity to use one of the cases
of the parity problem solved in \cite{He3}: for $a$, $b$ varying
over the all integers with $\gcd(a,b)=1$, the expression 
\[(a|b) \cdot \lambda(a b (a - b))\] 
averages to zero.
Here $(a|b) = \prod_{p|b,\; \text{$p$ odd}} (a/p)^{v_p(b)}$, where
$(\cdot/\cdot)$ is the quadratic reciprocity symbol:
$(a/p)=1$ if $a\in (\mathbb{Q}_p^*)^2$,
$(a/p)=1$ if $a\in \mathbb{Q}_p\setminus (\mathbb{Q}_p^*)^2$.

Recall that the Liouville function $\lambda$ is defined by
$\lambda(n) = \prod_p (-1)^{v_p(n)}$. We write $\rad(n)$
for the radical $\rad(n) = \prod_{p|n} p$ of an integer $n$.

\begin{lem}\label{lem:rootab}
For any $a,b\in \mathbb{Z}$ non-zero and coprime, the elliptic curve
$E_{a,b}:y^2 = x(x+a)(x+b)$ satisfies
\begin{equation}\label{eq:basuga}
\begin{aligned}
\frac{W(E_{a,b})}{W_2(E_{a,b})} &= - \prod_{p|b,\; \text{$p$ odd}} (a/p) \cdot
\prod_{p|a,\; \text{$p$ odd}} (b/p) \cdot
\prod_{p|(b-a),\; \text{$p$ odd}} (-a/p) \\
&\cdot \prod_{p|(a b (a-b))} (-1),\end{aligned}\end{equation}
where $W(E)$ is the (global) root number and $W_2$ is the local root number
at $2$.
\end{lem}
\begin{proof}
A model of the form $y^2 = x^3 + C_2 x_2 + C_1 x_1
+ C_0$ is minimal at an odd prime $p$ if $\min_{j=0,1,2}(v_p(C_j))=0$.
(This is so because 
the passage from a model of this form to a long Weierstrass form
(and vice versa)
is integral over $p$ odd; see \cite[\S III.1]{Si}.)
Since $a$ and $b$ are coprime, the model $y^2 = x (x + a) (x +b) =
x^3 + (a + b) x^2 + a b x$ is minimal at every odd prime $p$.

Let $p$ be an odd prime. We see from our minimal model
$y^2 = x (x + a) (x + b)$ that
the curve $E_{a,b}$ has good reduction at $p$ if $p\nmid a b (a-b)$,
and multiplicative reduction at $p$ if $p| a b (a-b)$.
Suppose, then, that $p| a b (a-b)$. 
The reduced curve $\widehat{E_{a,b}}$ over $\mathbb{Z}/p\mathbb{Z}$
 is given by $y^2 = x^2 (x+b)$ for $p|a$,
by $y^2 = x^2 (x+a)$ for $p|b$, and by $y^2 = x (x+a) (x+a)$ for $p|b-a$.
Hence, for $p$ odd, the reduction is split iff $b$ is a square mod $p$
(for $p|a$), iff $a$ is a square mod $p$ (for $p|b$), iff $-a$ is a square
mod $p$ (for $p|b-a$). Since
$W_p(E_{a,b}) = -1$ when the reduction is split, and 
$W_p(E_{a,b})=1$ when it is not,
the statement follows by $W(E) = - \prod_p W_p(E)$.
\end{proof}

Let us simplify (\ref{eq:basuga}) a little. First of all, we notice that
\[\begin{aligned}
\prod_{p|b,\; \text{$p$ odd}} (a/p) &\cdot
\prod_{p|a,\; \text{$p$ odd}} (b/p) \cdot
\prod_{p|(b-a),\; \text{$p$ odd}} (-a/p) \\&=
\prod_{p^2|b,\; \text{$p$ odd}} (a/p)^{v_p(b)-1} \cdot
\prod_{p^2|a,\; \text{$p$ odd}} (b/p)^{v_p(a)-1} \cdot
\prod_{p^2|(b-a),\; \text{$p$ odd}} (-a/p)^{v_p(b-a)-1} \\
&\cdot
(a|b)\cdot (b|a)\cdot (-a|b-a) .
\end{aligned}\]
By (\ref{eq:ednes}) (essentially a form of the law of quadratic reciprocity),
\[\frac{(b-a|-a)}{(-a|b-a)} = \left(\frac{a,b}{\infty}\right)^{-1} \cdot
\left(\frac{a,b}{2}\right)^{-1} ,\]
where $\left(\frac{a,b}{v}\right)$ is the quadratic Hilbert symbol (see
(\ref{eq:utur})). (What will matter  
is that, by (\ref{eq:utur}),
 $(a,b)\mapsto \left(\frac{a,b}{v}\right)$ depends on $a$ and $b$
only modulo $(K_v^*)^2$.)
Now 
\[(b-a|-a) = (b-a|a) = (b|a).\]
Since $(b|a)=\pm 1$, it follows that
\[\begin{aligned}
(a|b)\cdot (b|a)\cdot (-a|b-a) &= (a|b) \cdot (b|a) \cdot (b|a)
\left(\frac{a,b}{\infty}\right)^{-1} \cdot
\left(\frac{a,b}{2}\right)^{-1} \\ &=
\left(\frac{a,b}{\infty}\right)^{-1} \cdot
\left(\frac{a,b}{2}\right)^{-1}\cdot (a|b).\end{aligned}\]

Lastly, since $a$ and $b$ are coprime,
\[\prod_{p|(a b (a-b))} (-1) 
= \lambda(a b (a- b)) \cdot
\prod_{p^2|(a b (a-b))} (-1)^{v_p(a b(a-b))-1}.\]

Therefore, (\ref{eq:basuga}) implies that
\begin{equation}\label{eq:zand}\begin{aligned}
\frac{W(E_{a,b})}{W_2(E_{a,b})} &= 
\left(\frac{a,b}{\infty}\right)^{-1} \cdot
\left(\frac{a,b}{2}\right)^{-1} \cdot
\prod_{p^2|b,\; \text{$p$ odd}} (a/p)^{v_p(b)-1} \cdot
\prod_{p^2|a,\; \text{$p$ odd}} (b/p)^{v_p(a)-1} \\ 
&\cdot
\prod_{p^2|(b-a),\; \text{$p$ odd}} (-a/p)^{v_p(b-a)-1} \cdot
\prod_{p^2|(a b (a-b))} (-1)^{v_p(a b(a-b))-1}\\
&\cdot (a|b) \lambda(a b (a- b)) .\end{aligned}\end{equation}

\begin{prop}\label{prop:avrab}
For $a,b\in \mathbb{Z}$ coprime, let $E_{a,b}$ denote the curve
$y^2 = x (x + a) (x + b)$. Let $S\subset \mathbb{R}^2$ be a sector
and $a+L\subset \mathbb{Z}^2$ a lattice
coset. Then
\[\av_{S\cap (a + L), \text{coprime}} W(E_{a,b}) = 0.\]
\end{prop}
\begin{Rem}
Recall that the average $\av_{S\cap (a + L), \text{coprime}} W(E_{a,b})$ is
defined
to be
\[\lim_{N\to \infty} \frac{\mathop{\sum_{(a,b)\in S \cap (a +  L):\; 
    \gcd(a,b)=1}} 
 W(E_{a,b})}{|\{(a,b)\in S \cap (a +  L): 
    \gcd(a,b)=1\}|}.\]
Here the two terms $W(E_{0,1})$ and $W(E_{1,0})$ are ill-defined
(since $E_{0,1}$ and $E_{1,0}$ are singular curves); this is irrelevant,
as their contribution goes to $0$ as $N\to \infty$.
\end{Rem}
\begin{proof}[Proof of Proposition \ref{prop:avrab}]
The curve $y^2 = x (x+a) (x+b)$ is an elliptic curve whenever $a,b\ne 0$.
Thus,
by Proposition \ref{prop:waggle}, the map $(a,b)\mapsto W_2(E_{a,b})$
is defined and locally constant everywhere on 
$(\mathbb{Q}_2\setminus \{0\}) \times (\mathbb{Q}_2 \setminus \{0\})$.
Hence, (\ref{eq:zand}) can be rewritten in the following form:
for all non-zero, coprime $a,b\in \mathbb{Z}$,
\[W(E_{a,b}) = (g_\infty(a,b) \cdot g_2(a,b)) \cdot \prod_{\text{$p$ odd}} 
h_p(a,b) \cdot (a|b) \lambda(a b (a-b)),
\]
where $g_{\infty}:(\mathbb{R}\setminus \{0\}) \times (\mathbb{R}\setminus \{0\}) 
\to \{-1,1\}$ is everywhere locally constant,
$g_2:(\mathbb{Q}_2\setminus \{0\}) \times (\mathbb{Q}_2\setminus \{0\}) 
\to \{-1,1\}$ is everywhere locally constant, and
$h_p$ is a function from $\mathbb{Q}_p\times \mathbb{Q}_p$ to $\{-1,1\}$
that (a) is defined and locally constant on the complement
of the lines $x=0$, $y=0$, $x=y$ in the plane
$\mathbb{Q}_p\times \mathbb{Q}_p$, (b) satisfies
$h_p(a,b)=1$ when $p^2\nmid a b (a-b)$.

 (Explicitly, $g_{\infty}(a,b) = 
\left(\frac{a,b}{v}\right)^{-1}$ and $g_2(a,b) = 
\left(\frac{a,b}{2}\right)^{-1} W_2(E_{a,b})$.) 

We can now apply Proposition \ref{prop:castrol} with 
$\alpha(x,y) = (a|b) \lambda(a b (a-b))$ and
$B(x,y) = x y (x-y)$. Since $x y (x-y)$ is a product of linear factors,
hypothesis $\mathscr{A}_2(B(x,y))$ is true (and easy to prove). 
Condition (\ref{eq:thoto}) of Proposition \ref{prop:castrol}
requires that
\begin{equation}\label{eq:vane}
\av_{S\cap (a + L), \text{coprime}} (a|b) \lambda(a b (a-b)) = 0\end{equation}
for all sectors $S$ and all lattice cosets $a + L$; this is true by
Proposition 5.1 of \cite{He3}.
\end{proof}


\begin{thebibliography}{CCH04}
\bibitem{cs}
Cassels, J. W. S., and A. Schinzel, Selmer's conjecture and families of
elliptic curves, {\em Bull. London Math. Soc.} \textbf{14} (1982), 345--348.
\bibitem{Ch}
Chowla, S., {\em The Riemann hypothesis and Hilbert's tenth problem}, 
{\em Mathematics and Its Applications}, Vol. 4, Gordon and Breach Science
Publishers, New York--London--Paris, 1965.
\bibitem{Con1}
Connell, I., Calculating root numbers of elliptic curves over $\mathbb{Q}$,
{\em Manuscripta Math.} {\bf 82} (1994), 93--104.
\bibitem{CCG}
Conrad, B., K. Conrad, and R. Gross, Prime specialisation in genus $0$,
{\em Trans.\ Am.\ Math.\ Soc.} {\bf 360} (2008), 2867--2908.
\bibitem{CCH}
Conrad, B., K. Conrad, and H. Helfgott,
Root numbers and ranks in positive characteristic, {\em Adv. Math.}
{\bf 198} (2005), 684--731.
\bibitem{Co}
Conrad, B., and K. Conrad, personal communication.
\bibitem{De} Deligne, P., Les constantes des \'equations fonctionelles
des fonctions $L$, \emph{Modular Functions of One Variable, II}, \emph{\,SLN}
349, Springer-Verlag, New York, 1973, 501--595.
\bibitem{DD}
Dokchitser, T., and V. Dokchitser, Root numbers of elliptic curves in residue
characteristic 2, Bull.\ London Math.\ Soc.\ {\bf 40} (2008), 516--524.
\bibitem{Er}
Erd\H{o}s, P., Arithmetical properties of polynomials, {\em J. London Math. Soc.}
{\bf 28} (1953), 416--425.
\bibitem{Es}
T. Estermann, Einige S\"atze \"uber quadratfreie Zahlen, {\em Math. Ann.}
{\bf 105} (1931), 653--662.
\bibitem{FI} Friedlander, J., and H. Iwaniec, Asymptotic sieve for
primes, \emph{Ann. of Math. (2)} \textbf{148} (1998), 1041--1065.
\bibitem{Gran} Granville, A., $ABC$ allows us to count squarefrees, 
\emph{Internat. Math. Res. Notices} \textbf{1998}, 991-1009.
\bibitem{Gr} Greaves, G., Power-free values of binary forms, \emph{%
Quart. J. Math. Oxford} \textbf{43}(2) (1992), 45-65.
\bibitem{GT} Green, B., and T. Tao, 
Linear equations in the primes, arxiv:math/0606088.
\bibitem{Ha} Halberstadt, E., Signes locaux des courbes elliptiques en 2
et 3, \emph{C. R. Acad. Sci. Paris S\'er. I Math.} \textbf{326} (1998),
1047--1052.
\bibitem{HB} Heath-Brown, D. R.,
Primes represented by $x^3 + 2 y^3$,
{\em Acta Math.} 186 (2001), 1--84.
\bibitem{HBM}
Heath-Brown, D. R., and B. Z. Moroz,
On the representation of primes by cubic 
polynomials in two variables, 
{\em Proc. London Math. Soc. (3)} {\bf 88} (2004), 289--312.
\bibitem{He} Helfgott, H. A.,
On the square-free sieve, {\em Acta Arith.} {\bf 115} (2004), 349--402.
\bibitem{Heirr} Helfgott, H. A., The parity problem for
irreducible cubic forms, submitted.
\bibitem{He3} Helfgott, H. A., The parity problem for
reducible cubic forms, {\em J. London Math. Soc.} {\bf 73}, no.\ 2 
(2006), 415--435.
\bibitem{Ki} Kisin, M., 
Local constancy in $p$--adic families of Galois representations,
{\em Math. Z.} {\bf 230} (1999), 569--593.
\bibitem{Ko} Kobayashi, S., The local root number of elliptic
curves with wild ramification, {\em Math. Ann.} {\bf 323}
(2002), 609--623.
\bibitem{Ma} Manduchi, E., Root numbers of fibers of elliptic surfaces, 
\emph{Compositio Math.} \textbf{99} (1995), 33--58.
\bibitem{Mi} Miller, S. J., One- and two-level densities
for rational families of elliptic curves: evidence for the 
underlying group symmetries, {\em Compositio Math.} {\bf 140} (2004),
 952--992.
\bibitem{Ne} Neukirch, J., Algebraic number theory, Springer--Verlag, Berlin,
1999. 
\bibitem{Po} Poonen, B., Square-free values of multivariable polynomials,
{\em Duke Math. J.} {\bf 118} (2003), 353--373.
\bibitem{Ra} Ramsay, K., Square-free values of polynomials 
in one variable over function fields, {\em Internat. Math. Res. Notices}
{\bf 1992}, 97--102.
\bibitem{Riz1} Rizzo, O. G., Average root numbers in families of
elliptic curves, \emph{Proc. Amer. Math. Soc.} \textbf{127} (1999),
1597--1603.
\bibitem{Riz2} Rizzo, O. G., Average root numbers for a non-constant
family of elliptic curves, {\em Compositio Math.} \textbf{136} (2003), 1--23.
\bibitem{Ro} Rohrlich, D. E., Variation of the root number in families
of elliptic curves, \emph{Compositio Math.} \textbf{87} (1993),
119--151.
\bibitem{Ro2}
Rohrlich, D. E., Elliptic curves and the Weil-Deligne group, {\em Elliptic
curves and related topics}, 125--157, {\em CRM Proc. Lecture Notes} {\bf 4}
AMS, Providence, RI, 1994.
\bibitem{Rog} 
Rohrlich, D. E., Galois theory, elliptic curves, and root numbers,
\emph{Compositio Math.} \textbf{100} (1996), 311--349.
\bibitem{Roc} Rohrlich, D. E., personal communication.
\bibitem{ST} Serre, J.-P., and J. Tate, Good reduction of abelian 
varieties, {\em Ann. of Math. (2)} {\bf 88} (1968), 492--517.
\bibitem{Sisp} Silverman, J. H., Heights and the specialization map for
families of abelian varieties, {\em J. Reine Angew. Math.} {\bf 342} (1983),
197--211.
\bibitem{Si} Silverman, J. H., \emph{The arithmetic of elliptic curves},
Springer-Verlag, New York, 1985.
\bibitem{Wo} Wong, S., On the density of elliptic curves,
{\em Compositio Math.} {\bf 127} (2001), 23--54.
\bibitem{Y1} Young, M., Low-lying zeroes of families of elliptic curves,
{\em J. Amer. Math. Soc.} {\bf 19} (2006), 205--250.
\bibitem{Y2} Young, M., Lower-order terms of the $1$-level density of families
of elliptic curves, {\em Int. Math. Res. Not.} {\bf 2004}, 587--633.
\end{thebibliography}
\end{document}